\documentclass[10pt,reqno]{amsart}
\usepackage{mathrsfs}
\usepackage{amssymb}

\usepackage{setspace}     
\usepackage[colorlinks=true,linkcolor=blue,citecolor=blue]{hyperref}
\DeclareRobustCommand{\SkipTocEntry}[5]{}

\usepackage{enumitem}
\usepackage{bbm}
\usepackage{mleftright}
\renewcommand{\left}{\mleft}
\renewcommand{\right}{\mright}

 \usepackage[abbrev,nobysame]{amsrefs}

\renewcommand{\MR}[1]{}
\usepackage[capitalize,sort]{cleveref}

\crefformat{section}{\S#2#1#3}

\theoremstyle{plain}
\newtheorem{theorem}{Theorem} [section]

\newtheorem{proposition}[theorem]{Proposition} 

\newtheorem{claim}[theorem]{Claim}
\crefname{claim}{Claim}{Claims}

\newtheorem{lemma}[theorem]{Lemma}
\newtheorem{corollary}[theorem]{Corollary}
\theoremstyle{definition}
\newtheorem{definition}[theorem]{Definition}
\theoremstyle{remark}
\newtheorem{remark}[theorem]{Remark}
\newtheorem*{remark*}{Remark}

\renewcommand{\top}{\mathrm{top}}
\def\btm{\mathrm{bot}}

  \newcommand{\restrict}[2]{{#1}{\restriction_{{ #2}}}}

\usepackage{marginnote}

\newcommand{\SL}{\mathrm{SL}}

\DeclareMathOperator \KL{D_{KL}\mkern-1mu}

\newcommand{\id}{\mathrm{Id}}

\newcommand{\sm}{\smallsetminus}
\newcommand{\R}{\mathbb {R}}

\newcommand{\Z}{\mathbb {Z}}
\newcommand{\N}{\mathbb {N}}

\renewcommand {\vec}[1]{\mathbf #1}

\newcommand{\inv}{^{-1}}

\def\vecp{\vec p}

\DeclareMathOperator{\diff}{Diff}
\DeclareMathOperator{\diam}{diam}
\DeclareMathOperator{\Homeo}{Homeo}
\newcommand{\Diff}{\diff}

\DeclareMathOperator{\supp}{supp}

\def\loc{\mathrm{loc}}

\def\calP{\mathcal P}

\def\partP{\mathcal P}
\def\partR{\mathcal R}
\def\partQ{\mathcal Q}

\def\sigA{\mathscr A}
\def\sigB{\mathscr B}
\def\sigC{\mathscr C}
\def\sigF{\mathscr F}
\def\sigG{\mathscr G}

 \newcommand{\oldepsilon}{\mathchar"10F}
\newcommand{\eps}{\oldepsilon_0}

\renewcommand{\emph}[1]{{\bf{#1}}}

\renewcommand{\bar}{\overline}

\title[Entropy formula for Furstenberg entropy] {Entropy formula for Furstenberg entropy: invariance principle and exact dimensionality of contracting stationary   measures} 

\author[A.~Brown]{Aaron Brown}
\address{Northwestern University, Evanston, IL 60208, USA}
\email{awb@northwestern.edu}

\author[{F.~Rodriguez Hertz}]{Federico Rodriguez Hertz}
\address{Pennsylvania State University, State College, PA 16802, USA}
\email{hertz@math.psu.edu}

\newcommand{\hol}[1][\beta]{\,\textnormal{H\"ol}^{#1}}
\def\Hol{\hol}
\def\Holder{H\"older }

\date{\today}

\def\furstent{h^{\mathrm{Furst}}}

\newcommand{\wtd}{\widetilde}
\newcommand{\td}{\tilde}
\newcommand{\what}{\widehat}
\renewcommand{\hat}{\widehat}

\def\ptpart{\mathcal E}
\def\trivpart{\mathcal N}

\def\transent{h^{\mathrm{tr}}}
\def\maxent{h^{\mathrm{ref}}}

\def\maxFent{\maxent_\mu(T\mid \sigF)}

\def\maxFentC{\maxent_\mu(F\mid \sigC)}

\def\siname{invariant}

\allowdisplaybreaks
\begin{document}

\begin{abstract}
We revisit the invariance principle for stationary probability measures for i.i.d.\ random walks and for measures invariant under fibered dynamics admitting stable holono\-mies.  Assuming certain moment and entropy conditions, we establish  sufficient and necessary conditions for stiffness or holonomy invariance of such measures through a Ledrappier--Young-type formula for the Furstenberg entropy.  

To establish such a formula, we show exact dimensionality of stationary (or projected) measures along the filtration by stable manifolds.   This establishes exact dimensionality of the unique stationary probability measure in the setting of (non-affine) contracting iterated function systems.

\end{abstract}

\maketitle
\setcounter{tocdepth}{2}

\tableofcontents

\section{Introduction and overview}
The purpose of this paper is twofold.  First, we revisit the ``invariance principle'' for ergodic stationary probability measures for i.i.d.\ random walks and for measures invariant under fibered dynamics admitting holonomies.  We provide a sufficient and necessary condition for the invariance principle to hold and obtain a Ledrappier--Young-type formula for the Furstenberg entropy of such measures.   Second,  to formulate the Ledrappier--Young-type formula, we must establish the existence of certain pointwise dimensions of stationary measures along contracting foliations.  For contracting iterated function systems (or for stationary probability measures with all Lyapunov exponents negative), this establishes  exact dimensionality of stationary probability measures.

In the remainder of this brief introduction, we further discuss the history of the above two goals and outline of the paper.

\subsection{Invariance principle}
A distribution $\nu$ on the space $C^0(X,X)$ of continuous maps of a compact metric space $X$ generates an i.i.d.\ random walk on $X$.  Such i.i.d.\ random walks always admit stationary probability measures though may not, in general, admit any probability measure invariant under every generator of the random walk.  The ``invariance principle'' provides a  sufficient, though  generally not necessary, condition 
for ``stiffness''---that is, conditions for which a $\nu$-stationary probability measure is invariant under every $f\in C^0(X,X)$ in the support of $\nu$---of a stationary probability measure in the setting of i.i.d.\ random walks.  
In the setting of fibered dynamics admitting stable holonomies,  the invariance principle provides a  sufficient condition for ``holonomy invariance'' of the fiber disintegration of a measure invariant under the total dynamics.  
In both settings, a sufficient and necessary condition for invariance is the vanishing of a certain entropy, the  \emph{Furstenberg} entropy that we designate by $\furstent$ below.    Roughly, the invariance principle states that 
non-negativity of the Lyapunov exponents for the fiberwise derivative cocycle is a sufficient  condition for vanishing of this entropy; see \cref{AV-inv,InvPrinc:stiffness} below for formulation in the setting of i.i.d.\ random walks on manifolds.

For a comprehensive discussion of the invariance principle, see the recent survey \cite{LedrappierIP} of Ledrappier. 
The origin of the invariance principle appears in Furstenberg's proof of the positivity of the top Lyapunov exponent for i.i.d.\ random products of  matrices in $\SL(d,\R)$ (satisfying  the usual moment condition); see \cite[Thm.\  8.6]{MR0163345}.  
For i.i.d.~products of matrices, Furstenberg studies a quantity (now called the Furstenberg entropy) which vanishes precisely when a stationary measure for the projective random walk is invariant; when the stationary measure is in the Lebesgue class,  Furstenberg shows the entropy is given by the top Lyapunov exponent  (up to dimensional normalization); 
see especially \cite[Lem.~8.9]{MR0163345} and its corollary.  
For strongly irreducible products of matrices, an inequality bounding the infimum of the entropy over all stationary measures by the top Lyapunov exponent  is shown by Ledrappier in \cite[Thm.~III.3.1]{MR0876081}; here, Ledrappier approximates a strongly irreducible system by a system admitting stationary measures in the Lebesgue class and uses (semi-)continuity properties of exponents and entropy (see \cite[Thm.~III.2.2, Cor.~III.1.2]{MR0876081}).  
What is often called the invariance principle for i.i.d.~products of  matrices in $\SL(d,\R)$ is a bound on the Furstenberg entropy for a specific stationary measure in terms of the extremal negative Lyapunov exponent for the projective cocycle (and the dimension of the projective space).  Such a formula seems to first appear (at least in print) in Ledrappier's paper \cite{MR850070}; see \cite[Thm.~3]{MR850070}.  
Ledrappier's bound is sufficient to deduce Furstenberg's theorem, but the invariance principle in \cite[Thm.~3]{MR850070}  applies more broadly.  
In the i.i.d.~setting, it requires neither irreducibility of the cocycle nor absolute continuity of the stationary measure; furthermore, the formulation of  \cite[Thm.~3]{MR850070}  gives an invariance principle for measures invariant under the projectivized cocycle for stationary, possibly dependent, random matrix products (subject to the usual moment conditions).

For random walks and for continuous-time stochastic flows of diffeomorphisms on compact manifolds, Baxendale studied the case of  ergodic stationary probability measures that are absolutely continuous with respect to an ambient volume form.  
In this setting, Baxendale showed that the sum of all Lyapunov exponents (counted with multiplicity) is non-positive and vanishes if and only if the Furstenberg entropy vanishes; see \cite[Thm.~4.2]{MR0995809} as well as a proof (which is rigorous if the density is uniformly bounded away from $0$ and $\infty$) in \cite[Thm.\ V.1.3]{MR884892}.  In particular, 
Baxendale's criterion provides a sufficient and necessary condition for stiffness of absolutely continuous stationary probability measures.  (We note that Baxendale's proof in  \cite[Thm.\ 4.2]{MR0995809} requires a mild extra integrability hypothesis on the density of the stationary measure to avoid terms of the form $\infty - \infty$; this extra hypothesis is removed by \cite[Lem.\  2]{MR877387} of Le\,Jan's treatment in \cite[Prop.\  1(b)]{MR877387}.)

For general stationary probability measures in the setting of  random walks or stochastic flows by $C^1$ diffeomorphisms, a bound on Furstenberg entropy in terms of the bottom Lyapunov exponent was established by Crauel; see especially \cite[Thm.\ 5.1]{MR1204858} and \cite[Thm.\ 6.1]{MR1066618}.   Crauel's  result implies the most general versions typically used in the literature, but it seems his results were not well known to many in smooth dynamics; for example, the results in the appendix to \cite{BDZ} by the first author follow directly from  \cite[Thm.\ 6.1]{MR1066618}.  
The invariance principle  for (bounded) $C^1$ cocycles was popularized in smooth dynamics by the more recent paper of Avila and Viana  \cite{MR2651382}, which also establishes versions of the invariance principle stable under perturbations.   For further applications in smooth (especially partially hyperbolic) dynamics, see also  \cite{MR3203216}.

The setting for our results will be i.i.d.\ random walks by  $C^{1+\text{H\"older}}$ diffeomorphisms and more general skew-transformations whose fiber dynamics are $C^{1+\text{H\"older}}$.  We assume the expected integrability conditions and a mild (though perhaps unnecessary) finite entropy condition on the generators for the random dynamics; this setup includes i.i.d.\ random walks with finitely many generators.  

In the presence of negative Lyapunov exponents, the bound provided by Ledrappier, Crauel, or Avila--Viana, does not apply to provide a useful condition to verify  vanishing of Furstenberg entropy.  We take a different (though related) approach and define two distinct notions of the contribution to the total entropy of a system by negative Lyapunov  exponents and two distinct notions of pointwise dimension of the measure corresponding to negative Lyapunov exponents.  These two notions of entropy contribution  and dimension satisfy inequalities which become equalities if and only if  the Furstenberg entropy vanishes; see \cref{eq:invprinc,invprinc:general} below.  Combining these notions of dimension, Lyapunov exponents, and entropy contribution, we obtain a Ledrappier--Young-type formula for the Furstenberg entropy of the system.  See \cref{thm:furstentformula,entformskewfurst} below.

\subsection{Exact dimensionality}
To prove a Ledrappier--Young-type formula for  Furstenberg entropy,  
we establish that certain pointwise dimensions of a stationary probability measure along the filtration by fast stable manifolds  
are well defined.  In particular, while the  existence and properties of stable and unstable pointwise dimensions of the so-called ``sample measures'' (see \cref{sec:sampmsr} below) are well established in the literature, the existence of such pointwise dimension  for stationary probability measures along stable manifolds seems to be established only in the settings of  contracting affine,  projective, or  conformal, dynamics.  
We establish  the existence of pointwise dimension for stationary probability measures along stable manifolds for general $C^{1+\text{H\"older}}$ dynamics; see \cref{prop:exact:iid} below.

Beyond pointwise dimensions along stable manifolds, one may ask about the existence of a well-defined pointwise dimension (that is, exact dimensionality) of an ambient (ergodic) stationary probability measure.  
This  question has been well studied in the  setting of contracting iterated function systems; see discussion in \cref{rem:hyp and compare}  
for a non-exhaustive discussion of related results establishing exact dimensionality and entropy formulas for  stationary probability measures in the setting of  contracting  affine, projective, or conformal  iterated function systems.

A direct consequence of our main theorems is the exact dimensionality of ergodic stationary probability measures with all Lyapunov exponents negative for finitely generated random walks. This follows since  stable manifolds are open in the ambient manifold.  In particular, we obtain exact dimensionality for the unique stationary probability measure  $\mu_0$ for any contracting iterated function system generated by finitely many uniformly $C^{1+\text{H\"older}}$ diffeomorphisms.  See  \cref{thm:exactIFS} below.

\begin{remark*} In a forthcoming work, the first author together with Homin Lee 
extend the 
 exact dimensionality to the setting of ergodic hyperbolic stationary probability measures.  That is, for an iterated function system generated by finitely many uniformly $C^{1+\text{H\"older}}$ diffeomorphisms, if $\mu$ is an ergodic stationary measure with all Lyapunov exponents nonzero, then $\mu$ is exact dimensional.  The dimension of an ergodic  hyperbolic stationary measure $\mu$ is the sum of the stable dimension of the stationary measure and of the unstable dimension of the sample measures (see \cref{sec:sampmsr}.)
\end{remark*}

\subsection{Brief outline of paper}
In \cref{sec:iid}, we formulate our results for stationary probability measures generated by i.i.d.\ random walks on compact manifolds.  We explain in \cref{sss:transfromIFSto RW} how to translate results from the setting of random walks on compact manifolds to the setting of contracting iterated function systems (on the  non-compact manifold $\R^d$).  
In \cref{sec:abserg}, we recall basic facts about entropy for measure-preserving transformations.  We also define a number of notions of entropy relative to a sub-invariant sub-$\sigma$-algebra.  In general, these various notions of entropy relative to a sub-$\sigma$-algebra
are distinct and the remainder of the paper is devoted to the study of these discrepancies in  smooth settings.  
In \cref{formulation and results}, we present our general formulation of fibered systems admitting holonomies and state the main theorems of the paper; we reduce all announced results to the main technical results of the paper,  \cref{thm:realthm}.  
In \cref{sss:preparatory} we formulate a number of preparatory lemmas, mimicking those in \cite{MR819556,MR819557}; in \cref{LYagain}, we revisit the Ledrappier--Young machinery, following \cite{MR819556,MR819557}, to establish our main technical result,  \cref{thm:realthm}.  

\subsection{Origins of this paper}
This paper originates from the authors' previous paper \cite{MR3671937}, specifically the proof of \cite[Prop.~11.1]{MR3671937}.  Here, the authors needed to establish a certain invariance principle for a measure on a two-sided skew product admitting both positive and negative fiber Lyapunov exponents.  Because of the negative exponent, the standard invariance principle did not apply.  However, certain geometric assumptions on the support of the stable conditional measures imply, in the language of this paper, that both the ``transverse'' conditional entropy and ``standard'' conditional entropy vanish.  The arguments in this paper then show the Furstenberg entropy vanishes and  one obtains invariance of the fiber conditional  measures  under stable holonomies.  

In \cite{MR3671937}, the assumption on the support of the stable conditional measures allowed us to establish invariance of fiber conditional  measures  under stable holonomies using more geometric arguments.   However, the heuristic of the argument suggested a general framework: the negative fiberwise Lyapunov exponents contribute to two distinct notions of conditional entropy, the difference of the two entropies should be the Furstenberg entropy, and that Furstenberg entropy satisfies a Ledrappier--Young-type entropy formula.   This paper formalizes this heuristic.  While the ``standard" conditional entropy is well studied, we introduce the second notion of conditional entropy, which we refer to as ``transverse'' conditional entropy.  We show their difference is the Furstenberg entropy and, via a third notion of entropy which we call the ``refinement'' entropy, establish Ledrappier--Young-type entropy formulae for all these entropies.

\addtocontents{toc}{\SkipTocEntry}
\subsection*{Acknowledgements}
The arguments  in \cref{sss:preparatory,LYagain} of this paper naturally follow arguments from the papers of Ledrappier and Young \cite{MR819556,MR819557} and their many generalizations  including \cite{MR968818,MR953818,MR2530169,MR2448134,MR2218998,MR1944407,MR1369243,MR2032494}.   The authors benefited greatly from discussions with Francois Ledrappier and Homin Lee while writing  this paper. 

\addtocontents{toc}{\SkipTocEntry}
\subsection{Use of LLMs}
The authors used a paid subscription to a LLM to audit and referee early drafts of this manuscript, especially to identify missing hypotheses, incomplete proofs, and errors in mathematics and prose.  It also assisted in literature search and summarizing existing literature.  

The LLM also suggested the statement and proof of \cref{lem:superbconvexity} (the ``Data processing inequality'') in order to prove  \cref{subaddfurstent}, simplifying (and correcting in the infinite entropy case) an existing proof of   \cref{subaddfurstent} by the authors.  All other mathematical arguments are due to the authors (and the cited literature).  

All writing was done by the authors with LLM assisting only for copyediting.     

\section{Formulation of results for i.i.d.\ random walks and contracting IFSs} \label{sec:iid}
Although we will consider a more general setup 
in \cref{formulation and results} below, it seems useful to first formulate our results in the setting of i.i.d.\ random walks and contracting iterated function systems.    
\subsection{Random walks and stationary   measures} 
Let $M$ be a compact metric space and let $\nu_0$ be a probability measure on $C^{0}(M)$, the space of $C^0$ self-maps of $M$.  Starting in \cref{sec:lyap}, we will restrict to the case that $M$ is a compact manifold and that $\nu_0$ is supported on $\diff^{r}(M)$ for $r>1$.
We say a Borel probability measure $\mu_0$ on $M$ is \emph{$\nu_0$-stationary} if $$\int f_*\mu_0 \, d \nu_0(f) = \mu_0.$$ 
We sometimes refer to $\nu_0$ as the \emph{driving measure}.

Given $n\ge 1$, let $\nu_0^{\otimes n}$ denote the $n$th product measure on the semi-group $C^0(M)$; that is, $\nu_0^{\otimes n}$ is the image of $\nu_0\times \nu_0\times \cdots \times \nu_0$ under the map $(f_1,\dots, f_n)\mapsto f_n\circ \cdots \circ f_1$.    We note that if $\mu_0$ is $\nu_0$-stationary then $\mu_0$ is $\nu_0^{\otimes n}$-stationary for all $n\ge 1$.  
\subsubsection{Associated skew product}\label{skewext}
Let $C= \supp(\nu_0)$ and let $\widehat \Sigma= C^\N.$
Equip $\widehat \Sigma$ with the measure $\hat \nu :=\nu_0^\N$ and the $\hat \nu$-preserving (left-)shift $\sigma\colon \widehat \Sigma\to\widehat \Sigma$.  

We write $\what F\colon \widehat \Sigma\times M\to  \widehat \Sigma\times M$ for the  skew extension, 
$$\what F\left((f_0,f_1,f_2,\dots),x\right) = \left((f_1,f_2,\dots),f_0(x)\right) 
.$$
Let $$\hat \mu:=\hat \nu\times \mu_0= \nu_0^\N\times \mu_0.$$
We observe  the following well-known fact (which seems to  appear first in \cite[Lem.\ 2.1]{MR0700941}; see also  proof in \cite[Lem.\ I.2.3]{MR884892}): 
 \begin{lemma}{$\mu_0$ is $\nu_0$-stationary if and only if $\hat \mu$ is $\what F$-invariant.}\end{lemma}
 Moreover, for $n\ge 1$, $\mu_0$ is $\nu_0^{\otimes n}$-stationary  if and only if $\hat \mu$ is $\what F^n$-invariant.
 
Using the skew-extension formulation, we say $\mu_0$ is  an \emph{ergodic} $\nu_0$-stationary probability measure if $\hat \mu$ is ergodic as an $\what F$-invariant probability measure.

\subsection{Entropies of stationary probability measures for i.i.d.\ random walks}
We define various notions of entropy for  stationary probability measures for i.i.d.\ random walks.  
For the reader unfamiliar with notions of entropy, all notation and definitions are presented  in a more general setting in \cref{sec:abserg}.
For this subsection, we may take $M$ to be a compact metric space,
$\nu_0$ be a probability measure on $C^{0}(M)$,  and $\mu_0$ a $\nu_0$-stationary probability measure.  

\subsubsection{Standard entropy}
Let $\xi$ be a measurable partition of $M$ with $H_{\mu_0}(\xi)<\infty$.
The (standard) \emph{entropy} of the random dynamical system determined by $\nu_0$ with respect to the $\nu_0$-stationary measure ${\mu_0}$ and partition $\xi$ is (cf.\ \cite[Thm.\ II.1.3]{MR884892} or \cite[Thm.\ 2.1, Def.\ 2.1]{MR1369243})
$$h_{\mu_0}({\nu_0}, \xi) := \lim _{n\to \infty} \frac 1 n\int  H_{\mu_0} \left( \bigvee _{i=0}^{n-1} \bigl(f_\omega ^{(i)}\bigr)\inv \xi \right) \, d {\nu_0^\N}(\omega).$$
The (standard) \emph{entropy} of the random dynamical system determined by $\nu_0$ with respect to the $\nu_0$-stationary measure ${\mu_0}$
$$h_{\mu_0}({\nu_0})=\sup \{ h_{\mu_0}({\nu_0}, \xi)\}$$
where the supremum is over all measurable partitions $\xi$ of $M$ with $H_{\mu_0}(\xi)<\infty$.  

\subsubsection{Transverse entropy}
Let $\xi$ be an arbitrary measurable partition of $M$.  
Given $\omega\in \widehat \Sigma$, write
 $$\xi_\omega^{-} = \bigvee _{i=0}^\infty \bigl(f_\omega ^{(i)}\bigr)\inv \xi.$$
The \emph{transverse entropy} of the random dynamical system determined by $\nu_0$ with respect to the $\nu_0$-stationary measure ${\mu_0}$ and partition $\xi$ is 
\begin{align*}\transent_{\mu_0}({\nu_0}, \xi) 
&= \int H_{\mu_0} \bigl(\xi |f_\omega \inv \xi_{\sigma(\omega)}^-\bigr) \, d {\nu_0^\N} (\omega).
\end{align*}
We have $$\transent_{\mu_0}({\nu_0}, \xi) = \lim _{n\to \infty} \int  H_{\mu_0} \Bigl( \xi \mid \bigvee _{i=1}^n \bigl(f_\omega ^{(i)}\bigr)\inv \xi \Bigr) \, d {\nu_0^\N}(\omega)$$
assuming any term appearing in limit is finite; in particular, this holds if $H_{\mu_0}(\xi)<\infty$.

The \emph{transverse entropy} of the random dynamical system determined by $\nu_0$ with respect to the $\nu_0$-stationary measure ${\mu_0}$ is 
$$\transent_{\mu_0}({\nu_0}):= \sup\{ \transent_{\mu_0}({\nu_0}, \xi)\}$$
where the supremum is over all measurable partitions $\xi$ of $(M,\mu_0)$.

Let $\what \sigF$ denote the sigma algebra on $\what X:= \widehat \Sigma\times M$ given by preimages of the map $\what X\to \widehat \Sigma$. 
In terms of the notation presented in \cref{fiberentdef} below, we have (see \cref{computefiberent})
$$h_{\mu_0}({\nu_0}) = h_{\hat \mu}(\what F\mid \what  \sigF)$$
and 
$$\transent_{\mu_0}({\nu_0}) = \transent_{\hat \mu}(\what F\mid \what \sigF).$$
We note the following inequality always holds (see \cref{prop:basicineq} below).
\begin{claim}$h_{\mu_0}({\nu_0})\le\transent_{\mu_0}({\nu_0}).$\end{claim}

\subsubsection{Furstenberg entropy}
Given a  continuous self-map $f\colon M\to M$, let  $$J^*_{\mu_0}(f) := \dfrac {d {\mu_0}} {df_*{\mu_0}}$$ denote the Radon--Nikodym derivative of the measure ${\mu_0}$ with respect to $f_*{\mu_0}$; that is, 
$J^*_{\mu_0}(f)\in L^1({f_*\mu_0})$ and $${\mu_0}
 = J^*_{\mu_0}(f) f_*{\mu_0} + \eta_f^*$$ where $\eta_f^*$ is a Borel measure with $\eta_f^*\perp f_*{\mu_0}$.  
We define the    \emph{Furstenberg entropy} of a ${\nu_0}$-stationary probability measure ${\mu_0}$ with respect to the random walk determined by ${\nu_0}$ to be $$\furstent_{\mu_0}({\nu_0}):= - \int \int \log (J^*_{\mu_0}(f)(f(x)) ) \, d{\mu_0}(x) \, d {\nu_0}(f).$$
We observe that $0\le \furstent_{\mu_0}({\nu_0})\le \infty$ and that the following are equivalent:
\begin{enumerate}
\item $f_*{\mu_0}\ll{\mu_0}$;  
\item $J^*_{\mu_0}(f)(x)>0$ for $f_*{\mu_0}$-a.e.\ $x$.
\end{enumerate}
In particular, if $\furstent_{\mu_0}({\nu_0})<\infty$ then $f_*{\mu_0}\ll{\mu_0}$ for ${\nu_0}$-a.e.\ $f\in C^{0}(M)$.

When $f\colon M\to M$ is a homeomorphism, we similarly let  $$J_{\mu_0}(f) := \dfrac {d f^{-1}_*{\mu_0}} {d{\mu_0}}.$$ 
Then
$J_{\mu_0}(f)\in L^1({\mu_0})$ and $f^{-1}_*{\mu_0}
 = J_{\mu_0}(f) {\mu_0} + \eta_f$ where $\eta_f$ is a Borel measure with $\eta_f\perp {\mu_0}$.  
When $\nu_0$ is supported on $\Homeo(M)$, we have the equivalent definition of the {Furstenberg entropy} of a ${\nu_0}$-stationary probability measure ${\mu_0}$ 
given by $$\furstent_{\mu_0}({\nu_0}):= - \int \int \log (J_{\mu_0}(f)(x) ) \, d{\mu_0}(x) \, d {\nu_0}(f).$$

Of particular interest is the case that the Furstenberg entropy vanishes; this is the setting for the  so-called invariance principle (cf.\ \cite{MR2651382,MR850070,MR3203216,MR1204858,MR0995809,MR1066618}).  

\begin{proposition}\label{InvPrinc:stiffness}
We have $\furstent_{\mu_0}({\nu_0})=0$ if and only if $f_*{\mu_0} = {\mu_0}$ for ${\nu_0}$-a.e.\ $f\in C^{0}(M)$.
\end{proposition}
See \cref{claim:zero-furstenberg-entropy} below for a self-contained proof.  

\subsection{Lyapunov exponents and smooth ergodic theory of random walks}\label{sec:lyap}
We now make the standing assumption that $M$ is a compact manifold and that $\nu_0$ is supported on $\diff^r(M)$ for some $r\ge 1$.

Given $\omega\in \widehat \Sigma$ with $\omega = (f_0, f_1, f_2, \dots)$ and $n\ge 1$, we write $$f_\omega^{(n)}:=f_{n-1}\circ \cdots \circ f_{0}$$
for the $n$-fold iterate given by the word $\omega$.  (In particular, $f_\omega:=f_\omega^{(1)}= f_0$.) 

When $f\mapsto \log \|f\inv\|_{C^1}$ is $L^1({\nu_0})$, for $\mu_0$-a.e.\ $x$, the bottom Lyapunov exponent is defined as the a.s.~limit (which exists by Kingman's subadditive ergodic theorem)
$$\lambda_\btm(x):= \lim_{n\to \infty} \frac 1 n \log \bigl \|(D_xf^{(n)}_\omega )^{-1}\bigr\|^{-1}.$$

When we  further assume  the maps $f\mapsto \log \|f\|_{C^1}$ and $f\mapsto \log \|f\inv\|_{C^1}$ are $L^1({\nu_0})$ we may apply   Oseledec's theorem (cf.\ \cite[Thm.\ III.1.1]{MR884892} for a formulation in the setting of random dynamics): 
given ${\nu_0^\N}$-a.e.\ $\omega$, ${\mu_0}$-a.e.\ $x$, and every $0\neq v\in T_xM$, 
the limit $$\lim_{n\to \infty} \frac 1 n \log \|D_xf^{(n)}_\omega v\|$$ exists and is finite.  As $T_xM$ is finite dimensional, this limit takes finitely many ($p(x)$ possibly depending on $x$) values  
$$\lambda_\btm(x)=\lambda_1(x)<\lambda_2(x)<\dots<\lambda_{p(x)}(x).$$
We let $s(x) = \max \bigl(\{0\} \cup \{i: \lambda_i(x)<0\} \bigr )$.  Throughout, we will almost always assume  for a.e.\ $x$ that  $0<s(x)\le p(x)$. 
Note that if ${\mu_0}$ is ergodic, the numbers $s(x), p(x)$, and each $\lambda_i(x)$ are independent of $x$ a.s.

Assume now that $\mu_0$ is ergodic.  Write $V^0_{ \omega}( x) =\{0\}$.  For $1\le i\le p$, 
let $$V^i_{ \omega}(x) := \left\{v\in T_xM: \limsup_{n\to \infty}  \frac 1 n \log \|D_x f_{ \omega}^{(n)}v\|\le \lambda_i
\right\}.$$
The \emph{multiplicity} of $\lambda_i$ is the almost surely constant value of $m_i:=\dim \left(V^i_{ \omega} (x)/V^{i-1}_{  \omega}(x)\right).$

Let $d_M$ denote the distance on $M$ induced by a choice of Riemannian metric on $M$.  Suppose $s\ge 1$.   Given $(\omega,x)\in  \widehat \Sigma\times M$ and $1\le i\le s$, let 
$$W^i_{\omega}(x):=\left\{y\in   M:  \limsup_{n\to \infty}  \frac 1 n \log d_M( f_{ \omega}^{(n)}(x),f_{ \omega}^{(n)}(y)) \le \lambda_i\right\}$$
be the \emph{$i$th fast fiberwise stable set} through $x$ associated with $\omega$.
Write $W^i(\omega,x) = \{\omega\}\times W^i_{\omega}(x).$
   The collection of sets $\{W^i{(\omega,x)}: (\omega,x)\in \widehat \Sigma\times M\}$ forms a (generally non-measurable) partition of $\widehat \Sigma\times M$.

Fix $r>1$ and assume that ${\nu_0}$ is supported on $\diff^r(M)$.  If we further assume that  $f\mapsto \log \|f\|_{C^r}$ and $\log \|f\inv\|_{C^r}$ are $L^1({\nu_0})$, we may apply a fibered version of Pesin's stable manifold theorem (cf.\ \cite[Thm.\  III.3.2]{MR1369243} for the case of $C^{1+\mathrm{Lip}}$ maps or discussion in \cite{MR4599404} for setting of $C^{1+\text{\Holder}}$ maps)  and   obtain the following: ${\nu_0^\N}$-a.e.\ $\omega$ and  ${\mu_0}$-a.e.\ $x$, the set $W^i_{\omega}(x)$ is an injectively immersed, $C^{r}$ submanifold tangent to $V^i_{\omega}(x)$ at $x$.

\subsection{Kushnirenko-type upper bound on Furstenberg entropy}
Much of the existing literature provides a sufficient---though not necessary---condition for the vanishing of Furstenberg entropy $\furstent_{\mu_0}({\nu_0})$ by a Kushnirenko-type (cf.\ \cite{MR0177405}) upper bound.  

Let $\nu_0$ be supported on $\diff^1(M)$ and let ${\mu_0}$ be a ${\nu_0}$-stationary probability measure.  
The following bound was proved  in \cite[Thm.\  3]{MR850070} in the setting of projective  random walks.  
For general random walks by $C^1$ diffeomorphisms, the bound was established in \cite[Thm.\ 5.1]{MR1204858} and \cite[Thm.\ 6.1]{MR1066618} 
assuming   $\log \|f_\omega\inv \|_{C^1}\in L^1(\nu_0)$; see also the more  recent  \cite[Prop.\  3.1]{MR2651382} where a similar bound is shown assuming  $\log \|f_\omega\inv \|_{C^1}$ is bounded. 
\begin{proposition}[Invariance principle]\label{AV-inv}
Suppose 
$\log \|f\inv\|_{C^1}$ is $L^1({\nu_0})$.  Then 
$$\furstent_{\mu_0}({\nu_0})\le \dim(M)   \int \max\{-\lambda_\btm(x),0\} \, d {\mu_0}(x).$$
\end{proposition}

It is remarked in the literature (see for instance discussion after  \cite[Prop.~3.1]{MR2651382}) that one should be able to replace the  Kushnirenko-type upper bound in \cref{AV-inv} with a  Ruelle-type (see \cite{MR0516310}) upper bound where   $\furstent_{\mu_0}({\nu_0}) $  is bounded above by the negative sum of all negative Lyapunov exponents with their multiplicities.  Though such an inequality should hold in the $C^1$ setting, to the best of our knowledge no such 
inequality appears in the literature.  
 As a minor consequence of our main results, under a  finite entropy assumption on the driving measure ${\nu_0}$ for the random walk (in particular, we only consider the case of discrete $\nu_0$), we establish such a  Ruelle-type upper bound on  $\furstent_{\mu_0}({\nu_0}) $ for random walks by $C^{1+\text{\Holder}}$ diffeomorphisms; see \cref{thm:furstentformulaRuelle}. This follows from our precise  formula for $\furstent_{\mu_0}({\nu_0})$ in terms of the negative Lyapunov exponents and certain pointwise dimensions; see \cref{thm:furstentformula} below.

\subsection{Standing assumptions on i.i.d.\ random walks}\label{standing ass}
To state our results, we make the following standing assumptions on our random walk.  
Let $M$ be a compact manifold, let $r=1+\beta$ for $\beta>0$, and let $\nu_0$ be a probability measure supported on $\diff^r(M)$.  We assume
\begin{enumerate}
\item (Discreteness.) ${\nu_0}$ is a discrete measure; namely there is a countable set $C\subset \diff^{1+\beta}(M)$ with ${\nu_0}(C)=1$.
\item(Finite entropy.)  $H({\nu_0}):= -\sum_{f\in C} {\nu_0}(f) \log {\nu_0}(f) <\infty$.
\item\label{inthyp} (Finite $\log$-moment.) The assignments 
$f\mapsto \log  \|f\|_{C^{1+\beta}}$ and $f\mapsto \log  \|f\inv\|_{C^{1+\beta}}$ are $L^1({\nu_0})$.
\end{enumerate}

Although a $\log$-integrability hypothesis $f\mapsto \log  \|f\inv\|_{C^{1+\beta}}\in L^1({\nu_0})$  in \eqref{inthyp} above is almost certainly necessary, it remains an open question to what extent we may relax the other assumptions.   In particular, the integrability of $f\mapsto \log  \|f\|_{C^{1+\beta}}$ is likely unnecessary.    It remains  open if the  finite entropy and discrete support hypotheses on $\nu_0$ are needed for our results. 
 See \cref{rem:hyp and compare} below for further discussion.

\subsection{Sample measures and two-sided skew extension}\label{sec:sampmsr}
Consider $\Sigma:= C^\Z$.  Given a two-sided word $\omega\in \Sigma$, the limit $$\td \mu_\omega:= \lim_{n\to \infty } \left(f_{\sigma^{-n}(\omega)}^{(n)}\right)_*{\mu_0}$$ exists (in the weak-$*$ topology on Borel probability measures on $M$) for $\nu_0^\Z$-a.e.\ $\omega\in \Sigma$.  
Moreover, $(f_\omega)_*\td \mu_\omega= \td \mu_{\sigma(\omega)}$ for ${\nu_0}^\Z$-a.e.\ $\omega$.  
 In the literature, the measures $\{\td \mu_\omega\}$ are often referred to as the \emph{sample measures} of the random dynamical system.  

We  write $  F\colon   \Sigma\times M\to    \Sigma\times M$,
$$  F(\omega,x) = (\sigma(\omega), f_{0}(x))$$
for the two-sided skew extension.  

 Let $\td \mu$ denote the unique $F$-invariant probability measure  on $\Sigma\times M$ whose projections to $\Sigma$ and $\widehat \Sigma\times M$, respectively, are ${\nu_0}^\Z$ and $\hat\mu={\nu_0^\N}\times {\mu_0}$ (cf.\ \cite[Prop.\  I.1.2]{MR1369243}).  Then 
 \begin{enumerate}
 \item $\td \mu$ is $F$-invariant, and
 \item  $F\colon (\Sigma\times M , \td \mu) \to (\Sigma\times M, \td \mu)$ is the (ergodic-theoretic) natural extension of $\what F\colon (\widehat \Sigma\times M , \hat \mu) \to (\widehat \Sigma\times M , \hat  \mu).$ 
 \end{enumerate}
Moreover,
 \begin{enumerate}[resume]
 \item  the conditional measures of $\td \mu$ relative to the partition into fibers of the projection $\Sigma\times M\to \Sigma$ are $ \{\td \mu_\omega\}_{\omega\in \Sigma}$.
 \end{enumerate}

\subsection{Exact dimensionality along stable manifolds}
Let ${\mu_0}$ be an ergodic, ${\nu_0}$-stationary probability measure on $M$. 
We remark that fast stable sets are defined using only the forwards dynamics and so are defined for every $(\omega,x)\in \widehat \Sigma\times M$ or $(\omega,x)\in  \Sigma\times M$.  

For $\omega\in \Sigma$ or $\omega\in \widehat \Sigma$ and for $1\le i\le s$, the sets $  W^i_{\omega}(x)$ form a (generally non-measurable) partition of a full $ \mu_0$-measure subset of $M$.  
A measurable partition $\xi$ of  $\widehat \Sigma\times M$ (resp.\ $ \Sigma\times M$) is \emph{subordinate} to the partition into fiber $W^i$-manifolds if for $\hat \mu$-a.e.\ $(\omega,x)\in \widehat \Sigma\times M$ (resp.\ $\td \mu$-a.e.\ $(\omega,x)\in  \Sigma\times M$) the following hold:
\begin{enumerate}
\item $\xi(\omega,x)\subset  \{\omega\}\times W^i_{\omega}(x)$; 
\item each atom $\xi(\omega, x)$ contains an open (w.r.t.\ the immersed manifold topology) neighborhood of $(\omega,x)$ in $ \{\omega\}\times W^i_{\omega}(x)$;
\item each atom $\xi(\omega,x)$ is precompact (w.r.t.\ the immersed manifold topology) in $ \{\omega\}\times W^i_{\omega}(x)$.
\end{enumerate}

Fix any  measurable partitions $\hat \xi^i$ and $\xi^i$ of $(\widehat \Sigma\times M,\hat \mu)$ and $(\Sigma\times M,\td \mu)$, respectively, {subordinate} to the partition into fiber $W^i$-manifolds. 
We may define families of conditional measures $\{\hat \mu^{\hat \xi^i }_{(\omega,x)}\}_{(\omega,x)\in \widehat \Sigma\times M}$
 and $\{\td \mu^{\xi^i }_{(\omega,x)}\}_{(\omega,x)\in  \Sigma\times M}$ of $\hat \mu$ and $\td \mu$, respectively, with respect to these partitions.  We may ask to what extent the measure of  ambient balls in $M$ decays relative to these measures.  
 
 \begin{proposition}[Exact dimensionality along stable manifolds] \label{prop:exact:iid}Assume the standing hypothesis in \cref{standing ass}.

For ${\nu_0^\N}\times {\mu_0}$-a.e.\ $(\omega,x)\in \widehat \Sigma\times M$ the limit 
$$  {\dim}^i(\hat \mu, (\omega,x)):= \lim_{r\to 0}\dfrac{\log \left(\hat \mu^{\hat\xi^i }_{(\omega,x)} (B(x,r)) \right)}{\log r}$$ 
exists and for  $\wtd \mu$-a.e.\ $(\omega,x) \in \Sigma\times M$ the limit
$${\dim}^i(\td  \mu, (\omega,x)):=  \lim_{r\to 0}\dfrac{\log \left(\td \mu^{\xi^i }_{(\omega,x)} (B(x,r)) \right)}{\log r}$$ 
exists.  

Moreover, the limiting values are $\what F$-invariant  and   $F$-invariant, respectively, and are independent of the choice of partition $\hat \xi^i$ and $\xi^i$ subordinate to fiber $W^i$-manifolds.  
\end{proposition}

We write 
$ {\dim}^i(\hat \mu)$ and $ {\dim}^i(\td \mu)$ for the almost surely  constant value of ${\dim}^i(\hat \mu, (\omega,x))$ and ${\dim}^i(\td \mu, (\omega,x))$, respectively.  
We note that $ {\dim}^i(\td \mu)$ is the pointwise dimension of the sample measures $\{\td\mu_\omega\}$ along  $W^i$-manifolds; $ {\dim}^i(\hat \mu)$ is the pointwise dimension of the stationary probability measure $\mu_0$ along $W^i$-manifolds.  
The existence of the pointwise dimensions ${\dim}^i( \td \mu)$ for the sample measures is well known in the existing literature (see e.g.\ \cite{MR2218998}); however, to the best of our knowledge the existence of the  pointwise dimensions $ {\dim}^i(\hat \mu)$ for the stationary probability measure $\mu_0$ along $W^i$-manifolds has not previously been established in any generality.  

In the proof of the existence of ${\dim}^i(\td \mu)$ and $ {\dim}^i(\hat \mu)$ we will show the following:
\begin{claim}\label{claim:totaldim}
Assume the standing hypothesis in \cref{standing ass}.
For $1\le i\le s$, 
${\dim}^i(\td \mu)\le {\dim}^i(\hat \mu).$
\end{claim}

We write $\lambda_\top(\mu_0)= \lambda_p$ for the top Lyapunov exponent of the ergodic stationary probability measure $\mu_0$.  In the case that $s=p$ and so $\lambda_\top(\mu_0)<0$, the full stable manifolds $W^s_{\omega}(x)$ are open in $M$ and thus the conditional measures $\hat \mu^{\hat \xi^s }_{(\omega,x)}$ are simply the (normalized) restriction of $\mu_0$ to 
$\hat \xi^s (\omega,x)$.  This gives the following exact dimensionality for uniformly contracting on average stationary   measures.  

\begin{corollary}[Exact dimensionality of contracting stationary  measures] \label{cor:exactdim} Assume the standing hypothesis in \cref{standing ass}.
If $\lambda_\top(\mu_0)<0$ then $\mu_0$ is exact dimensional: there exists $\delta\ge 0$ such that for $\mu_0$-a.e.\ $x\in M$,
$$\delta =\lim_{r\to 0}\dfrac{\log \left( \mu_0(B(x,r)) \right)}{\log r}.$$
\end{corollary}

\subsection{Entropy formula for Furstenberg entropy}
We continue to consider $\mu_0$, an ergodic $\nu_0$-stationary probability measure on $M$.  
For $1\le  i \le s$, write $\dim^0(\hat \mu) = \dim^0(\td \mu) = 0$ and 
$$\gamma^i(\hat \mu ) = \dim^i( \hat\mu )- \dim^{i-1}( \hat \mu), \quad \quad 
\gamma^i( \td \mu ) = \dim^i(\td  \mu )- \dim^{i-1}(\td \mu).$$
Similar to \cref{claim:totaldim}, we have the following. 
\begin{claim}Assume the standing hypothesis in \cref{standing ass}.
For $1\le i\le s$,
$$0\le \gamma^i( \td \mu )\le \gamma^i(\hat \mu )\le m_i.$$
\end{claim}

We have the following Ledrappier--Young-type entropy formula (cf.\ \cite{MR819557}) for the Furstenberg entropy.  
\begin{theorem}\label{thm:furstentformula}
Assume the standing hypothesis in \cref{standing ass}.  
Then  
\begin{align*}
\furstent_{\mu_0}({\nu_0})&= \transent_{\mu_0}({\nu_0})-  h_{\mu_0}({\nu_0})=- \sum_{i=1}^s\bigl(\gamma^i(\hat \mu )- \gamma^i( \td \mu )\bigr)\lambda_i.\end{align*}
\end{theorem}

As a corollary, we obtain the following Ruelle-type estimate.  
\begin{corollary}\label{thm:furstentformulaRuelle}
Assume the standing hypothesis in \cref{standing ass}.  
Then  
\begin{align*}
\furstent_{\mu_0}({\nu_0})&\le - \sum_{i=1}^s m_i \lambda_i.\end{align*}
\end{corollary}

\subsection{Sufficient and necessary conditions for vanishing of Furstenberg entropy}
The invariance principle in  \cref{AV-inv} provides a sufficient criterion for vanishing of Furstenberg entropy.  In our setting 
 we have the following sufficient and necessary conditions for vanishing of Furstenberg entropy, strengthening \cref{AV-inv}.  

\begin{theorem}\label{eq:invprinc}
Assume the standing hypothesis in \cref{standing ass}.  Let  $\mu_0$ be an ergodic $\nu_0$-stationary probability measure on $M$. The following are equivalent:
\begin{enumerate}
\item $\furstent_{\mu_0}({\nu_0})= 0$.
\item $  h_{\mu_0}({\nu_0})=\transent_{\mu_0}({\nu_0}).$
\item $\dim^i(\hat  \mu )= \dim^i( \td \mu )$ for every $1\le i\le s$.  
\item $\gamma^i( \hat \mu )= \gamma^i( \td \mu )$ for every $1\le i\le s$.  
\item $f_*\mu_0 = \mu_0$ for $\nu_0$-a.e.\ $f$.  
\item  for ${\nu_0}^\Z$-a.e.\ $\omega\in \Sigma$, we have $\td \mu_\omega= \mu_0$.
\end{enumerate}
\end{theorem}

\subsection{General invariant measures for   contracting IFSs and compactification}\label{sss:transfromIFSto RW}
A particular case of random walks is the setting of  (contracting) iterated function systems  on Euclidean space.  
We outline a broad class of measures supported on the attractor of a (contracting) iterated function systems associated to arbitrary measures on the associated shift space.
Also, since $\R^d$ is not compact, we outline a simple modification so that results for  random walks on compact manifolds apply to (contracting) iterated function systems on $\R^d$.  

\subsubsection{General invariant measures for contracting IFSs}
Let $f_1,\dots, f_k\colon \R^d\to \R^d$ be finitely many $C^{1+\beta}$ diffeomorphisms of $\R^d$ for some $\beta>0$. In the case of i.i.d.\ random walks,  let $\vecp$ be a probability vector on $S:= \{1, \dots, k\}$ and let $\mu_0$ be a compactly supported $\vecp$-stationary probability measure.  
Although $\mu_0$ need  not  exist for all such systems, if the maps $f_i$ are uniform contractions on $\R^d$ (that is, if $\max _i \sup_x \|D_xf_i\|<1$), then there exists a unique, compactly supported $\vecp$-stationary probability measure $\mu_0$.

More generally, assume each $f_i$ is a contraction. 
Let $\Sigma=S^\Z$.  
Pick any $\vec x_0\in \R^d$.   Given $\omega = (\dots, {\omega_{-2}}, {\omega_{-1}}, {\omega_0}, {\omega_1},\dots )\in \Sigma$, define $\Pi\colon \Sigma\to \R^d$ by 
\begin{equation}\label{IFS projec cantor} \Pi(\omega) = \lim_{n\to \infty}  f_{\omega_{-1}}\circ \cdots \circ f_{\omega_{-n}}(\vec x_0).\end{equation}
Using that each $f_i$ is a contraction, the limit exists for every $\omega\in \Sigma$ (and is independent of the choice of initial point $\vec x_0$).  We refer to $\Pi\colon \Sigma\to \R^d$ as the \emph{coding map}.

  Given any ergodic, Borel probability measure $\nu$ on the two-sided full shift $\Sigma:=S^\Z$, let $\mu_0= \Pi_*\nu$.  
We again have that $\mu_0$ is compactly supported.
The measure $\mu_0$ may also be recovered as follows: let $\mu$ denote the push-forward of the measure $\nu$ on $\Sigma \times \R^d$ under the graph $\omega\mapsto (\omega, \Pi(\omega))$ of $\Pi$.  Then $\mu$ is invariant (and ergodic) under the skew-extension $F\colon(\omega,x)\mapsto(\sigma(\omega), f_{\omega_0}(x)).$
Then $\mu_0= (\pi_2)_* \mu$ is the image of $\mu$ under the canonical projection $\Sigma\times \R^d\to \R^d$.  
We also work with the push forward $\hat \mu$ of $\mu$ to the one-sided skew extension $\what \Sigma\times \R^d$; $\hat \mu$ is invariant (and ergodic)  under the one-sided skew extension $\what F\colon\what \Sigma\times \R^d\to \what\Sigma\times \R^d$.

\subsubsection{Compactification of IFSs}
We will always assume that $\mu_0$ is compactly supported; as discussed above, this holds whenever $\{f_1,\dots, f_k\}$ are contractions.  Although $\R^d$ is not compact, compactness of the support of $\mu_0$ is all that is needed in arguments.  For completeness, we outline a reduction from $C^r$ ($r\ge 1$) contracting iterated functions systems on $\R^d$ to  random walks on compact manifolds.  

Let $B(0,T)$ denote the open ball in $\R^d$ centered at $0$ of radius $T>0$. 
There are $R>1$ and $\ell>0$ for which
\begin{enumerate}
\item $\supp(\mu_0)\subset B(0,R-1)$;
\item 
$f_i(\overline {B(0,R)})\subset {B(0,R+\ell)}$ for $1\le i\le k$.  
\end{enumerate}  
  Let $N= \overline {B(0,R+\ell)}$ and let $M$ be the double of $N$ across $\partial N$; that is, if $N_1$ and $N_2$ are two copies of $N$ then $M= (N_1\cup N_2)/\partial N\simeq S^d$.   View 
   $\overline {B(0,R)}$ as a subset of a common copy of $N$ and view 
  each  map $$\restrict{ f_i}{\overline {B(0,R)}}\colon \overline {B(0,R)}\to  {B(0,R+\ell)}\subset \operatorname{int}(N)\subset M$$ as a $C^{r}$ embedding of  $\overline {B(0,R)}$ into $M$. 
  By the Palais-Cerf extension theorem
(see \cite[Cor.~2]{MR0117741} and \cite{MR0140120}) each such embedding extends to a
$C^{r}$ diffeomorphism
\[
\widetilde f_i\colon M\longrightarrow M.
\]
See also the modern formulation
\cite[Thm.~1.4]{GGH25}.  
Although  \cite{GGH25} is only stated in the $C^r$ category for 
$r\in\mathbb N\cup\{\infty\}$, the  proof also applies in the $C^{1+\beta}$ category.  Also, since $S^d$ admits an orientation reversing diffeomorphism, the theorem also applies when any of the maps $f_i$ are orientation reversing.  
    Since $\supp (\mu_0)\subset B(0,R)$, the dynamics in a neighborhood of $\supp(\mu_0)$ is as prescribed by the original iterated function system.

  \subsection{Exact dimensionality and related results for contracting IFSs} 
  \subsubsection{Exact dimensionality}
From \cref{cor:exactdim} above in the case of $\vecp$-stationary probability measure for i.i.d.\ random walks or from \cref{prop:dimensions are defined} and \cref{thm:exactdimproj}, below, in the more general setting   where $\mu_0= \Pi_*\nu$ for a general shift-invariant probability measure,  we obtain   the following:
  \begin{theorem}\label{thm:exactIFS}
  Let $f_1,\dots, f_k\colon \R^d\to \R^d$ be finitely many $C^{1+\text{\Holder}}$ diffeomorphisms of $\R^d$. Assume each $f_i$ is a uniform contraction: there is $0<\kappa<1$ such that for every $1\le i\le k$ and every $x\in \R^d$,  
  $$\|D_xf_i\|\le \kappa<1.$$
Let $\mu_0$ be either 
  \begin{enumerate}
  \item an ergodic $\vecp$-stationary probability measure on $\R^d$ where   $\vecp$ is a probability vector on $S:= \{1, \dots, k\}$; or 
  \item $\Pi_*\nu$ where $\Pi$ is as in \eqref{IFS projec cantor} and  $\nu$ is an ergodic, shift-invariant  probability measure on $  \Sigma= S^\Z$.
  \end{enumerate}
  Then    $\mu_0$ is exact dimensional:  there exists $0\le \delta\le d$ such that for $\mu_0$-a.e.\ $x\in \R^d$
$$\delta =\lim_{r\to 0}\dfrac{\log \left( \mu_0(B(x,r)) \right)}{\log r}.$$
  \end{theorem}
The dimension of the measure $\mu_0$ in \cref{thm:exactIFS} also satisfies a certain entropy formula (for the Furstenberg entropy for the associated one-sided skew extension).  See \cref{cor:contracting}.

Note that \cref{thm:exactIFS} holds trivially  if $\nu$ is an ergodic measure with finite support.

\subsubsection{Relationship between entropies and dimensions}
Let $f_1,\dots, f_k\colon M\to M$ be $C^{1+\text{\Holder}}$ diffeomorphisms of either $M=\R^d$ or a compact manifold $M$.
Let $\vecp=\langle p_1, \dots, p_k\rangle $ be a probability vector on $S:= \{1, \dots, k\}$ and let $\mu_0$ be an ergodic $\vecp$-stationary probability measure.
We will assume that $f_i$ are contractions or, more generally, that $\lambda_\top(\mu_0)<0$.  

We have that $H(\vecp) = -\sum_{i=1}^k p_i \log p_i$ is the entropy of the Bernoulli shift determined by $\vecp$.  In the case that $f_i$ are conformal or, more generally, that $\lambda:= \lambda_\btm(\mu_0)=\lambda_\top(\mu_0)<0$, the  relationship between $H(\vecp)$, $\lambda$, and $\dim(\mu_0)$ is well-studied.  We recover the following well-known formula for the dimension.    

\begin{claim}\label{IFSconformal}
For a contracting iterated function system, for any $\vecp$, if the unique $\vecp$-stationary probability measure  $\mu_0$ has a single Lyapunov exponent $\lambda:= \lambda_\btm(\mu_0)=\lambda_\top(\mu_0)<0$ then $\mu_0$ is exact dimensional and $$\dim(\mu_0) =- \furstent_{\mu_0}(\vecp)/\lambda \le -H(\vecp)/\lambda.$$
\end{claim}

To understand the difference between  $\furstent_{\mu_0}(\vecp)$ and $ H(\vecp)$, let $X= S\times M$ be equipped with the product measure $\hat \mu:=\vecp\times \mu_0$.  The map $F\colon (X,\hat \mu)\to (M,\mu_0)$, $F(i,x) = f_i(x)$, is measure preserving.  Let $\ptpart_M$ be the point partition of $M$; then $F\inv \ptpart _M$ is a well-defined partition of $X$.  
We similarly let $\ptpart_X$ denote the partition into points on $X$.  

The following characterizes the defect between $ \furstent_{\mu_0}(\vecp)$ and $ H(\vecp)$.
\begin{claim}
$ \furstent_{\mu_0}(\vecp) = H(\vecp)- H_{\hat \mu}(\ptpart_X\mid F\inv \ptpart _M)$
\end{claim}

\begin{corollary}
For an iterated function system with a single Lyapunov exponent $\lambda<0$,  the measure   $\mu_0$ is exact dimensional and $$\dim(\mu_0) =-( H(\vecp)-H_{\hat \mu}(\ptpart_X\mid F\inv \ptpart _M)) /\lambda.$$
\end{corollary}

\subsection{Relationship to prior results and comments on hypotheses and regularity}\label{rem:hyp and compare}
For ergodic stationary probability measures $\mu_0$, and more generally for  projections $\Pi_*\nu$ of ergodic
shift-invariant measures $\nu$ under the coding of a contracting IFS, 
there are a number of previous results concerning the existence of pointwise dimensions and formulas involving entropy, dimension, and Lyapunov exponents.  
Much of the literature focuses on affine,  projective, and conformal  contracting (or contracting on average) random walks.  We discuss a number of  these settings and prior results here.  

\subsubsection{Conformal contracting on average systems} Dimensional properties of the (often unique) stationary measure (or image measure $\Pi_*\nu$ associated to an ergodic shift-invariant measure $\nu$) are well studied, especially in the case of contracting conformal random walks and IFSs.  
More generally, there are many results assuming the random walk is only asymptotically contracting or conformal.  
Given a random walk and an ergodic stationary measure $\mu_0$ or projected measure $\Pi_*\nu$, we will say the  derivative cocycle (or its restriction to certain invariant fibrations)  is \emph{Lyapunov conformal} if the Lyapunov spectrum consists of  a single  value $\chi$ and \emph{contracting on average} if $\lambda_{\top}<0$.  
In the setting of Lyapunov conformal  contracting on average random walks, the measure is often shown to be exact dimensional  and the product of the pointwise dimension and $(-\chi)$ coincides with an appropriate entropy of the system.   We outline several results in the conformal setting or when the Lyapunov spectrum consists of one (negative) value where exact dimensionality is established.  (We note exact dimensionality often fails in the presence of zero exponents; see for instance  \cite{MR2563101,MR0829860}.)  

In \cite[Part.~V]{MR0876081}, Ledrappier studied the Furstenberg stationary measure $\mu_0$ for i.i.d.~projective action of $\SL(2,\R)$ on $\R P^1$.  
Under standard hypotheses (which imply the Furstenberg measure is the unique stationary measure and has negative Lyapunov exponent), in \cite[Part V, Thm.\ 5.1, Cor.~5.5]{MR0876081}, Ledrappier identified the entropy formula in \cref{IFSconformal} using a slightly weaker notion of  dimension, where the convergence of $\frac{\log \left( \mu_0(B(x,r)) \right)}{\log r}$ to the dimension as $r\to 0$ is only convergence in  probability rather than almost sure pointwise convergence.  
The  first explicit  proof of exact dimensionality of the Furstenberg measure was given by  Hochman and Solomyak in  \cite[Thm.~3.4]{MR3735630}.  It is remarked in \cite[Rem.~1.2]{MR3735630} that this proof was known to experts.  
Although \cite[Thm.~3.4]{MR3735630} is stated more generally, the statement of \cite[Thm.~3.4]{MR3735630} omits the necessary moment condition to apply Oseledec's theorem, and its proof uses discreteness and an integrability step that is automatic for finitely supported driving measures. 
Thus the proof may only formally justify exact dimensionality in the case of finitely supported driving measures.

In the setting of conformal contracting on average random dynamics, discreteness of the driving measure (or countability of the IFS) seems unnecessary.  For instance, Ledrappier does not assume discreteness of the driving measure in \cite[Part V, Thm.\ 5.1, Cor.~5.5]{MR0876081} and only requires the necessary moment condition.  
Continuing  in the projective setting, we highlight the result of Lessa in \cite[Thm.~2]{MR4216247} which, in the setting of i.i.d.~projective random walks, established exact dimensionality of certain fiber measures along one-dimensional fibers; 
the product of the dimension and the absolute value of the fiber exponent recovers a certain relative Furstenberg entropy.   
Lessa does not require discreteness of the driving measure and only assumes the standard  moment condition ($\log$-integrability of all singular values), in addition to a number of non-degeneracy hypotheses on the random walk (ensuring the stationary measure is unique and contracting on average among other properties.)
A special case of Lessa's theorem recovers exact dimensionality of the Furstenberg measure on $\R P^1$ 
for arbitrary, possibly nondiscrete, driving measures satisfying the standard dynamical hypotheses and moment condition, justifying the conclusion of Hochman and Solomyak's \cite[Thm.~3.4]{MR3735630} for general driving measures satisfying the $\log$-moment  hypothesis.  
See also \cite[Thm.~2.7]{MR4845880} which provides the analogue of  Lessa's \cite[Thm.~2]{MR4216247} for more  projective actions on more general flag varieties associated to  Oseledec flags.  In \cite[Thm.~2.7]{MR4845880}, the fibers may be higher-dimensional but always have a single (negative) Lyapunov exponent.  Again, discreteness of the driving measure is not required; under appropriate assumptions on the random walk, the fiber measures are shown to be exact dimensional and satisfy an entropy formula.  

In nonlinear settings, there are several recent results establishing exact dimensionality and entropy formulas for ergodic stationary measures for random walks on the circle.  We  highlight   \cite{HJX23} which studies $C^2$ random walks with finitely supported driving measures and  \cite{SB25} which considers  $C^{1+\text{\Holder}}$ random walks with general driving measures satisfying a moment condition (the $C^1$ norms of the random diffeomorphisms  and their inverses are $\log$-integrable).  
Both \cite{HJX23,SB25} impose certain hypotheses on the random walk to ensure all stationary measures have negative exponents and establish the entropy formula \cref{IFSconformal} for Furstenberg entropy.

In the setting of contracting $C^1$ finitely generated IFSs, Feng and Hu 
 studied ergodic stationary measures $\mu_0$ and more general projected measures of the form $\mu_0=\Pi_*\nu$ for any ergodic $\nu$ on the shift space.   In \cite[Thm.~2.8]{MR2560042}, they establish exact dimensionality and an entropy formula (for what they call ``projection entropy'') assuming the derivative cocycle is Lyapunov conformal with respect to $\nu$.  
See also \cite[Thm.~10.1]{MR2560042} which establishes a similar result for certain countably generated similarity IFSs that are contracting on average (and satisfy certain entropy and moment  conditions).

For random uniformly contracting countable conformal IFSs driven by an ergodic base, Mihailescu and Urbański proved, under a finite conditional-entropy assumption, that the projected conditional measures are exact dimensional for almost every  base point and satisfy an entropy formula
in terms of ``random projectional entropy;'' 
see \cite[Thm.~3.13 and Cor.~3.14]{MR3505755}.  In later work, they considered uniformly contracting countable conformal IFSs with place-dependent probabilities and showed the associated stationary measure is exact dimensional and satisfies the analogous entropy formula; see \cite[Thm.~2.1]{MR4396787}.

We highlight the very recent \cite{2609.10538} that establishes exact dimensionality and the dimension formula for ergodic stationary measures for random $C^1$  diffeomorphisms with a single negative Lyapunov exponent (i.e.~with $\lambda_{\btm}=\lambda_{\top}$) under the standard moment condition.  Importantly, no discreteness assumption on the measure is needed.

\subsubsection{Beyond conformal systems} 
Outside the conformal setting, the majority of exact dimensionality and entropy formula results known to the authors are in the setting of random walks whose driving measure has finite or countable support or for IFSs generated by a finite or countable number of maps.

In the setting of affine contracting on average iterated function systems, the dimension theory of stationary and projected  measures $\Pi_*\nu$ has been extensively studied.  The most general  theorem that holds for all systems is Feng's paper \cite{MR4557759}; in \cite{MR4557759}, Feng studies finitely generated affine IFS that are contracting on average with respect to an ergodic shift-invariant measure $\nu$ and shows the projected measure $\Pi_*\nu$ is exact dimensional and satisfies a Ledrappier--Young-type entropy formula.  
We highlight recent precursors to Feng's paper establishing exact dimensionality with an entropy formula including results of 
Feng and Hu  in  \cite[Thm.~2.12]{MR2560042} for contracting affine IFS with pairwise commuting linear parts and results of B\'ar\'any  for planar self-affine measures under strong separation and dominated splitting and of B\'ar\'any and K\"aenm\"aki
 for all planar self-affine measures, with higher-dimensional analogues under additional spectral or domination hypotheses \cite{MR3413884,MR3689737}.  See also the more classical literature on carpets and sponges for affine IFSs including  \cite{MR0771063,MR1183358,MR2298824,MR1394572,MR1389626}.

 Like affine IFSs, there is also a long history of results for projective random walks.  Rapaport  studied projective actions on $\R P^{n-1}$ and established exact dimensionality and an entropy formula for stationary measures for finitely supported driving measures  (under strong irreducibility and proximality) in \cite{MR4273191}.   Ledrappier and Lessa  subsequently studied stationary measures on full and partial  flag varieties in $\R^d$ in \cite{MR4561150} and on more general flag varieties constructed from Oseledec flags in \cite{MR4845880}.
They obtain exact dimensionality and Ledrappier--Young-type entropy formulas assuming their driving measure is countably supported, in addition to the standard moment and appropriate non-degeneracy hypotheses.  To pass from a single Lyapunov exponent to the case of multiple (negative) Lyapunov exponents, discreteness of the driving measures seems necessary in the approach in  \cite{MR4561150,MR4845880}.  
 It seems plausible that the approach in \cite{MR4561150,MR4845880} could be adapted to our settings to show exact dimensionality assuming only the moment condition and discrete support in \cref{standing ass}.  However, outside the conformal setting, it remains an open question to what extent exact dimensionality results hold when the driving measure $\nu_0$ is not discrete.

Beyond affine and projective systems, a few results treat nonlinear nonconformal IFSs with additional geometric structure. 
Feng and Hu studied the direct product of contracting IFSs coded by the same subshift; assuming each component IFS is Lyapunov conformal, the image measure $\mu_0=\Pi_*\nu$ for the direct product  IFS is shown to be exact dimensional and satisfy an entropy formula in \cite[Thm.~2.11]{MR2560042}.   
   Jurga and Lee \cite{MR4381206} considered finite planar $C^{1+\alpha}$ IFSs with triangular structure, certain domination, and a certain separation condition.  They proved exact dimensionality and a Ledrappier--Young-type entropy formula for the projections of ergodic  quasi-Bernoulli measures, including
Bernoulli measures and Gibbs measures for Hölder potentials. 
B\'ar\'any and K\"aenm\"aki \cite{MR4808250} studied parametrized families of planar $C^2$ triangular IFSs with the opposite domination.  Under a transversality assumption, they proved that, for each quasi-Bernoulli symbolic measure and for almost every parameter, the projected measure is exact dimensional and its dimension equals its Lyapunov dimension.  See also Feng and Simon  where parametrized families of finite uniformly contracting $C^1$ IFSs are studied.  Under a transversality condition, they show in \cite[Thm.~1.2]{MR4492627} for each ergodic shift-invariant measure and  for almost every parameter, the dimension of the projected measure is the minimum of the ambient dimension and the  Lyapunov dimension.

To the best of our knowledge,  \cref{thm:exactIFS} is the first exact dimensionality theorem for the full  class of uniformly contracting iterated function systems generated by finitely many $C^{1+\text{\Holder}}$ diffeomorphisms.  In particular, our results do not require any  affine or projective structure or any product, triangular, domination, separation, or parameter-genericity assumptions appearing in  many prior nonlinear results.
Moreover, by \cref{cor:contracting} below, these measures satisfy an entropy formula for the Furstenberg entropy (for the associated one-sided skew extension).

\subsubsection{Critical regularity for exact dimensionality of contracting iterated function systems}
We remark that if the  generators  $f_1,\dots, f_k\colon \R^d\to \R^d$ of the iterated function system are only assumed to be bi-Lipschitz, then the conclusion of \cref{thm:exactIFS} is false.  This is shown  in  \cite[Thm.~6.6]{BBFLX} (see also  \cite[Thm.~1.6]{BBFLX}), where a   bi-Lipschitz iterated function system is constructed with the property that  $\Pi_*\nu$ is not exact dimensional for any ergodic, shift-invariant  probability measure $\nu$ with positive entropy on $  \Sigma= S^\Z$.  In particular, if $\nu =\vecp^\Z$ is a product (i.e.\ Bernoulli) measure, then the associated  $\vecp$-stationary probability measure $\Pi_*\nu$  for the iterated function system is not exact dimensional.

It seems unknown what regularity $C^r$ is needed to ensure all measures $\Pi_*\nu$ are exact dimensional for finitely generated, contracting $C^r$ iterated function systems.  The examples in  \cite{BBFLX} show  bi-Lipschitz regularity does not suffice whereas our \cref{thm:exactIFS} shows it is sufficient for the iterated function system to be  $C^{1+\text{\Holder}}$.  For conformal (or Lyapunov conformal) contracting iterated function systems, $C^1$ regularity is sufficient.  Results for parametrized families such as \cite{MR4492627} suggest $C^1$ regularity could possibly be sufficient.

\section{\texorpdfstring{Entropies in abstract ergodic theory relative to a $\sigma$-algebra}{Entropies in abstract ergodic theory relative to a σ-algebra}}\label{sec:abserg}

\subsection{Notation, conventions, and overview of entropies}
Let $(X,\sigB, \mu)$ be a standard probability space.
It will be useful at times to distinguish between measurable partitions and sub-$\sigma$-algebras of $(X,\mu)$.  Given a measurable partition  $\xi$ of $(X,\mu)$, write $\sigma(\xi)$ for the associated countably generated sub-$\sigma$-algebra of $\xi$-saturated sets; given a sub-$\sigma$-algebra $\sigA$ of $(X,\mu)$, we write $\Xi(\sigA)$ for the (unique mod $\mu$ choice of)  measurable partition into atoms of $\sigA$.  We write $\xi(x)$ for the element of $\xi$ containing $x$ and $\sigA(x)$ for the element of $\Xi(\sigA)$ containing $x$.  
At times, we abuse notation and write $\xi \vee \sigA$ for the measurable partition $\xi\vee \Xi(\sigA)$.

Let $T\colon X\to X$ be a measurable, $\mu$-preserving transformation.  
Let $\sigA$ be a sub-$\sigma$-algebra of $\sigB$.  
We say $\sigA$ is \emph{sub-invariant} if $T\inv \sigA\subset \sigA$ and \emph{\siname} if $T\inv \sigA= \sigA$. 
A measurable partition $\xi$ of $(X,\mu)$ is \emph{sub-invariant}  (resp.~\emph{invariant}) if the corresponding sub-$\sigma$-algebra $\sigma(\xi)$ is.  
If $\sigA$ is sub-invariant, it induces a measurable map $\what T\colon (X,\sigA, \mu)\to (X,\sigA,\mu)$ which we naturally identify with a factor system $\what T \colon (\what X, \hat \mu)\to(\what X, \hat \mu)$.  We have that factor transformations are in 1-1 correspondence with sub-invariant sub-$\sigma$-algebras.

In   the following subsections, given a sub-invariant sub-$\sigma$-algebra $\sigF$ of $(X,\mu)$, we will define {\it six} distinct notions of entropy of $T\colon (X,\mu)\to (X,\mu)$ relative to $\sigF$.  Briefly, these are:
\begin{enumerate}
\item $h_\mu(T, \Xi(\sigF))$, the entropy of $T$ relative to the partition into atoms of $\sigF$.
\item $h_\mu(T, \sigF)$, the entropy of $T$ relative to all partitions by  $\sigF$-saturated sets; this is the metric entropy of the factor induced by $\sigF$.
\item $h_\mu(T\mid \sigF)$, the ``standard'' conditional entropy of $T$ conditioned on (or relative to) $\sigF$.
\item $\transent_\mu(T\mid \sigF)$, the ``transverse'' conditional entropy of $T$ conditioned on (or relative to) $\sigF$.
\item $\maxFent$, the ``refinement entropy,'' i.e.~the maximal entropy of $T$ relative to all partitions $\xi$ with $\sigF\subset \sigma(\xi)$; in most cases of interest, we have  (see \cref{prop:GenAR})
$$\maxFent= h_\mu(T\mid \sigF)+ h_\mu(T, \Xi(\sigF)).$$  
\item $\furstent_\mu(T\mid \sigF)$, the ``Furstenberg''  entropy of $T$ relative to $\sigF$.
\end{enumerate}

When $\sigF$ is \siname, by \cref{equalifinv} below  we have 
\begin{equation}\label{relsthreeents}h_\mu(T\mid \sigF)=\transent_\mu(T\mid \sigF)= \maxent_\mu(T\mid \sigF)\end{equation}
and $$\furstent_\mu(T\mid \sigF)=0.$$
Similarly, when $T$ is an automorphism (and $\sigF$ is sub-invariant) by \cref{hbarhhatsame} below we have $$\transent_\mu(T\mid \sigF)= \maxent_\mu(T\mid \sigF).$$ 
To a large extent, this paper is devoted to the difference between the three quantities in \eqref{relsthreeents}  and their relationship to $\furstent_\mu(T\mid \sigF)$ when $T$ is assumed non-invertible and $\sigF$ is assumed to be only sub-invariant. 

The remainder of this section presents various properties and relationships among these entropies in the setting of abstract ergodic theory.  
See especially \cref{prop:basicineq,prop:GenAR,entropygap} below.  
In \cref{formulation and results} we present further properties of these entropies for fibered systems with smooth invertible fibered dynamics.

\subsection{Conventions and general properties of entropy}
\subsubsection{Conditional entropy}
Throughout, let $\trivpart=\{X\}$ denote the trivial partition of $(X,\mu)$ and let $\ptpart$ denote the partition  of $X$ into points (and so $\sigma(\ptpart)=\sigB$). 

Given a measurable partition $\xi$ on $(X,\mu)$, let $\{\mu_x^\xi\}$ denote a family of conditional measures relative to this partition.  
Given a sub-$\sigma$-algebra $\sigA$ we also write $\{\mu_x ^\sigA\} = \{\mu_x ^{\Xi(\sigA)}\}.$
Given measurable partitions $\xi$ and $\eta$, 
we define the \emph{conditional entropy} of $\xi$ relative to $\eta$   
 to be 
$$H_\mu(\xi\mid \eta):= -\int \log \mu_x^\eta(\xi(x))\, d\mu(x)$$
When $\eta$ is trivial, we write $H_\mu(\xi) := H_\mu(\xi\mid \trivpart)$ for   the \emph{entropy} of the partition $\xi$.  In particular, if $\xi$ is countable then $H_\mu(\xi) =-\sum_{C\in \xi}\mu(C) \log \mu(C)$.  
Given sub-$\sigma$-algebras $\sigA$ and $\sigC$ and measurable partitions $\xi$ and $\eta$ we also write 
$$H_\mu(\xi \vee  \sigC \mid\eta \vee  \sigA):= H_\mu(\xi\vee \Xi(\sigC) \mid \eta \vee  \Xi(\sigA))$$

\subsubsection{Basic properties of conditional entropy}\label{entfactsbasic}
We enumerate   basic properties that are used throughout.  Let $\xi,\eta,$ and $\zeta$ be measurable partitions and let $\sigF$ and $\sigG$ be  sub-$\sigma$-algebras.  
\begin{enumerate}
\item $H_\mu(\xi\vee \zeta\mid \sigF) = H_\mu( \zeta\mid \sigF) + H_\mu(\xi\mid  \zeta\vee\sigF)$.
\item If $\xi\prec \zeta$ then 
$H_\mu(\xi\mid \sigF)\le H_\mu(\zeta \mid \sigF)$.
\item If $\sigF\subset \sigG$ then 
$H_\mu(\xi\mid \sigG)\le H_\mu(\xi \mid \sigF)$.
\item  \label{keyentropy2} $H_\mu(\xi\mid \sigF) = H_\mu (T\inv \xi\mid T\inv \sigF)$ and  thus if $\sigF$ is sub-invariant then $  H_\mu (T\inv \xi\mid  \sigF) \le H_\mu(\xi\mid \sigF).$
\end{enumerate}

\subsubsection{Metric entropy with respect to a partition}
 Given a measurable partition  $\xi$ of $(X,\mu)$ and non-empty $S\subset \N$,  write 
\begin{enumerate}
\item $T\inv \xi:= \{T\inv (C): C\in \xi\}$;
\item  $\xi^{\emptyset} := \trivpart$;
\item $ \xi^S:= \bigvee_{i\in S} T^{-i} \xi$
 and so $\xi^{[0,n)}= \bigvee_{j=0}^{n-1}  T^{-j}(\xi) $; 
\item $\xi^- :=\xi^{[0,\infty)}= \bigvee_{j=0}^\infty  T^{-j}(\xi)$. 
\end{enumerate}
Also write \begin{enumerate}[resume]
\item $\xi^\emptyset =\trivpart$. 
\end{enumerate}
If $T$ is an automorphism and $S\subset \Z$,  we  also write 
\begin{enumerate}[resume]
\item $ \xi^S:= \bigvee_{i\in S} T^{-i} \xi$; 
\item $\xi^\Z= \bigvee_{j=-\infty}^\infty  T^{j}(\xi)$;
\item $\xi^+= \bigvee_{j=0}^\infty  T^{j}(\xi)$. 
\end{enumerate}
We note that for any measurable partition $\xi$, the $\sigma$-algebra  $\sigma(\xi^-)$ is sub-invariant and, if $T$ is an automorphism, we have that $\sigma(\xi^\Z)$ is \siname.  

The \emph{measure theoretic entropy} of $T$ with respect to $\xi$ is  \begin{equation}\label{eq:defofentreal} h_\mu(T, \xi) := H_\mu( \xi\mid T\inv \xi^-).\end{equation}
If $\xi$ satisfies $H_\mu(\xi)<\infty$ (in particular, if $\xi$ is finite) then (cf. \cite[Fact 2.3.4]{MR2809170})
\begin{equation}\label{eq:defofent}h_\mu(T, \xi)= \inf_n \frac 1 n H_\mu(\xi^{[0,n)}).\end{equation}
The \emph{measure theoretic entropy} of $T$  is  $$h_\mu(T):=\sup \{h_\mu(T, \xi)\}$$
where the supremum is taken over all measurable partitions of $(X,\mu)$.

\subsubsection{Metric entropy with respect to a sub-$\sigma$-algebra}
Given a sub-$\sigma$-algebra $\sigA$, there are two reasonable notions of measure theoretic entropy of $T$ with respect to $\sigA$: the entropy restricted to the partition into atoms or the entropy restricted to $\sigA$-saturated sets.  
We write \begin{equation}  h_\mu(T,\Xi(\sigA))\end{equation} to denote the entropy of $T$ with respect to the partition into atoms of $\sigA$ and, following the conventions of \cite{MR2809170}, write 
 \begin{equation} \label{eq:factorent} h_\mu(T,\sigA):= \sup _{\sigma(\xi)\subset \sigA} h_\mu(T,\xi).\end{equation}
Here, the supremum is taken over all $\sigA$-measurable partitions $\xi$ of $(X,\mu)$.  
In particular, if $\sigA$ is sub-invariant, and if $\what T\colon (\what X, \hat \mu)\to (\what X, \hat \mu)$ is the induced factor map of the quotient measure space then $$h_\mu(T, \sigA) = h_{\hat \mu}(\what T).$$

\subsection{\texorpdfstring{Notions of entropy relative to a $\sigma$-algebra}{Notions of entropy relative to a σ-algebra}}\label{fiberentdef}
For the remainder of this section, we let   $\sigF\subset \sigB$ be a  sub-invariant sub-$\sigma$-algebra.    Here, we  define three distinct notions of entropy of the transformation $T\colon X\to X$ relative to $\sigF$.

\begin{remark}[Discussion on notions of conditional entropy]
Given an \siname\ sub-$\sigma$-algebra  $\sigF$, one may naturally generalize (the suprema of expressions appearing in) equations \eqref{eq:defofent} and \eqref{eq:defofentreal}  to define a notion of conditional entropy relative to $\sigF$.  We present these generalized definitions in  \cref{standard conditional entropy} and \cref{sss:transent} below.  For \siname\ sub-$\sigma$-algebra $\sigF$, these definitions coincide; see  \cref{equalifinv} below.  
 In particular, one could  use the appropriate generalization of either \eqref{eq:defofent} or \eqref{eq:defofentreal} to define  conditional entropy for an \siname\ sub-$\sigma$-algebra  $\sigF$.  

However, when $\sigF$ is only assumed sub-invariant, the natural generalizations of the definitions (even when taking suprema) given by equations  \eqref{eq:defofent}  and \eqref{eq:defofentreal} differ.  For consistency with existing literature (especially the notion of conditional entropy appearing in the Abramov-Rohlin formula \cite{MR0140660,MR0476995,MR1179170}) we designate the generalization of  \eqref{eq:defofent} ``standard conditional entropy.''  Based on our geometric picture in fibered systems described in \cref{formulation and results} below, we designate our generalization of \eqref{eq:defofentreal} ``transverse conditional entropy.''  
Although these two definitions coincide when  $\sigF$ is \siname, much of  the sequel will focus on sub-invariant sub-$\sigma$-algebras  $\sigF$ where  equality in \cref{equalifinv} below does not hold.  We thus maintain and emphasize the distinction between ``standard'' and ``transverse'' conditional entropy throughout, even in situations where $\sigF$ is \siname.  
\end{remark}

\subsubsection{Standard conditional entropy}\label{standard conditional entropy}
Let $\xi$ be a measurable partition of $X$ with $H_\mu(\xi\mid \sigF)<\infty$.  
We define the  standard \emph{conditional entropy} of $T$ with respect to $\xi$ conditioned on $\sigF$ to be $$h_\mu(T,\xi \mid \sigF):= \lim_{n\to \infty} \frac 1 n H_\mu\left(\xi^{[0,n)}  \mid \sigF\right)=   \inf_{n} \frac 1 n H_\mu\left(\bigvee_{i=0}^{n-1} T^{-i}\xi \mid \sigF\right) $$
where the infimum in the final equality follows from subadditivity (using properties in \cref{entfactsbasic}).

We define the  standard \emph{conditional entropy} of $T$  conditioned on $\sigF$ to be $$h_\mu(T\mid \sigF):= \sup \{h_\mu(T,\xi \mid \sigF)\}$$
where the supremum is taken over all measurable partitions $\xi$ of $(X,\mu)$  with $H_\mu(\xi\mid \sigF)<\infty$.

\subsubsection{Transverse conditional entropy}\label{sss:transent}
We present an alternative notion of conditional entropy which, in general, is distinct from the definition in \cref{standard conditional entropy}.  To the best of our knowledge, this does not appear to have any name in the literature. Based on our geometric picture in fibered systems described in \cref{formulation and results} below, we designate this concept of entropy \emph{transverse conditional entropy}.

Let $\xi$ be any measurable partition of $(X,\mu)$.  
The 
 \emph{transverse conditional entropy} of $T$  with respect to $\xi$ conditioned on $\sigF$ is $$\transent_\mu(T,\xi\mid \sigF):= H_\mu\bigl(\xi  \mid T\inv \xi ^- \vee \sigF\bigr). $$ The 
 \emph{transverse conditional entropy} of $T$ conditioned on $\sigF$ is  $$\transent_\mu(T\mid \sigF):= \sup\{ \transent_\mu(T,\xi\mid \sigF)\}$$
where the supremum is taken over all measurable partitions $\xi$ of $(X,\mu).$

\subsubsection{Refinement entropy relative to $\sigF$}
We present one final notion of entropy relative to $\sigF$.  Again, to the best of our knowledge, this does not appear to have any name in the literature.   The \emph{refinement entropy} of $T$ relative to  $\sigF$ is 
$$\maxent_\mu(T\mid \sigF):= \sup \{ h_\mu(T,\xi\vee \Xi( \sigF))\}$$
where the supremum is taken over all measurable partitions $\xi$ of $(X,\mu).$ 
Alternatively, we have $$\maxent_\mu(T\mid \sigF)= \sup_{\sigF \subset \sigma(\xi)}  h_\mu(T,\xi).$$

\subsection{\texorpdfstring{Properties of entropies relative to $\sigF$}{Properties of entropies relative to ℱ}}
We present a number of properties of our notions of   entropy introduced above.

\subsubsection{Properties of standard conditional entropy}

We  often use the following, a direct adaptation of \cite[Fact 2.3.7]{MR2809170}.  
\begin{proposition}[{cf.\ \cite[Fact 2.3.7]{MR2809170}}] \label{downAR} 
Let $\sigF$ be a sub-invariant sub-$\sigma$-algebra.  
Let $\xi$ and $\eta$ be measurable partitions with 
$H_\mu(\xi \mid \eta\vee \sigF)<\infty$.  
Then  \[\lim_{n\to \infty} \frac 1 n H_\mu(\xi^{[0,n)} \mid \eta^{[0,n)} \vee \sigF)= h_\mu(T,\xi \mid  \eta^-\vee \sigF). \] 
\end{proposition}
Note that though \cite[Fact 2.3.7]{MR2809170} is formulated only for $\xi$ and $\eta$ countable, the proof  uses only the finiteness of the conditional entropy $H_\mu(\xi \mid \eta\vee \sigF)$.

We have the following elementary facts about standard conditional entropy.  
\begin{proposition}[Properties of standard conditional entropy ({cf.\ \cite[Fact 2.4.2]{MR2809170}})]\label{prop:standard properties}
Let $\sigF$  and $\sigG$ be sub-invariant sub-$\sigma$-algebras.  Let $\xi$ and $\zeta$ be measurable partitions with $H_\mu(\xi\mid \sigF)<\infty$ and 
$H_\mu( \zeta\mid \sigF)<\infty$.  
\begin{enumerate}
\item $h_\mu(T,\xi\mid \sigF)\le H_\mu(\xi\mid \sigF)$.
\item \label{standard2} $ h_\mu(T,\xi\mid \sigF)\le  h_\mu(T,\xi\vee \zeta \mid \sigF)$ and 
$$ h_\mu(T,\xi\vee \zeta \mid \sigF)\le  h_\mu(T,\xi \mid \sigF)+ H_\mu(\zeta \mid \xi \vee \sigF).$$
\item If $\sigF\subset \sigG$ 
then $h_\mu(T,\xi\mid \sigG)\le h_\mu(T,\xi\mid \sigF)$.

\end{enumerate}

\end{proposition}

For conclusion \eqref{standard2}, using \cref{downAR} we have 
\begin{align*} 
\lim_{n\to \infty} \frac 1 n &H_\mu\left(\xi^{[0,n)}  \mid \sigF\right)
\le \lim_{n\to \infty} \frac 1 n H_\mu\left(\xi^{[0,n)} \vee \zeta^{[0,n)} \mid \sigF\right)\\
&= \lim_{n\to \infty} \frac 1 n H_\mu\left(\xi^{[0,n)}  \mid \sigF\right) +\lim_{n\to \infty} \frac 1 n H_\mu\left(\zeta^{[0,n)}  \mid \xi^{[0,n)} \vee \sigF\right) \\
&=h_\mu(T,\xi \mid \sigF)+ h_\mu(T, \zeta\mid \xi^-\vee \sigF)\\
&\le h_\mu(T,\xi \mid \sigF)+ H_\mu(\zeta\mid \xi\vee \sigF).
\end{align*}

Suppose now that $T\colon X\to X$ is   a measurable automorphism.  Then $T\inv$ is measurable and $\mu$-preserving and $T\sigF$ is a well-defined sub-$\sigma$-algebra.    In this case we write $\sigF^+$ for the smallest $\sigma$-algebra containing $\{T^i\sigF: i\ge 0\}.$
If $\sigF$ is sub-invariant then $\sigF^+$ is  \siname.  Moreover, the induced invertible transformation on the quotient $(X,\mu)/\sigF^+$ coincides with the (ergodic-theoretic) natural extension of the induced transformation on the quotient  $(X,\mu)/\sigF$.  
\begin{lemma}[{\cite[Thm.\  0.4.4]{MR1369243}}]\label{natext}
Suppose that $T\colon X\to X$ is a measurable automorphism.  Then for any measurable partition $\xi$ with $H_\mu(\xi \mid \sigF)<\infty$, $$h_\mu(T,\xi \mid \sigF) = h_\mu(T,\xi \mid \sigF^+).$$
\end{lemma}

Note that if $T$ is an automorphism and if 
$H_\mu(\eta \mid \sigF)<\infty$ and $H_\mu(\xi \mid \eta\vee \sigF)<\infty$ 
then \cref{downAR,natext} imply
 \[\lim_{n\to \infty} \frac 1 n H_\mu(\xi^{[0,n)} \mid \eta^{[0,n)} \vee \sigF)= h_\mu(T,\xi \mid  \eta^\Z\vee \sigF^+). \]

\cref{downAR,natext} also give the following additivity formula for standard conditional entropy.  
\begin{lemma}[{cf. \cite[Facts 2.4.2 + 2.3.7]{MR2809170}}]
\label{lem:natext}
Let $\sigF$ be a sub-invariant sub-$\sigma$-algebra.  
Let $\xi$ and $\eta$ be measurable partitions with $H_\mu(\eta \mid \sigF)<\infty$ and $H_\mu(\xi \mid \eta\vee \sigF)<\infty$.  Then 
$$h_\mu(T,\xi \vee \eta \mid \sigF) = h_\mu(T,\eta  \mid \sigF)  + h_\mu(T,\xi \mid  \eta^-\vee \sigF) . $$
Moreover, if  $T$ is an automorphism then 
$$h_\mu(T,\xi \vee \eta \mid \sigF) = h_\mu(T,\eta  \mid \sigF)  + h_\mu(T,\xi \mid  \eta^\Z\vee \sigF^+) . $$

\end{lemma}
\begin{proof}
We have \begin{align*}
\lim_{n\to \infty} \frac 1 n  H_\mu(\xi^{[0,n)}\vee \eta^{[0,n)} \mid \sigF)
&= 
\lim_{n\to \infty} \frac 1 n \left[H_\mu(\eta^{[0,n)}  \mid \sigF )+ H_\mu(\xi^{[0,n)} \mid \eta^{[0,n)} \vee \sigF)\right]\\
&= 
 h_\mu(T,\eta \mid \sigF) + \lim_{n\to \infty} \frac 1 n H_\mu(\xi^{[0,n)} \mid \eta^{[0,n)} \vee \sigF).
\end{align*}
The result then follows from \cref{downAR} and \cref{natext}.  
\end{proof}

\subsubsection{A monotonicity estimate}

\begin{claim}\label{claim:mylifeboilsdowntoinfinityminusinfinity}
Let $\sigG$ be a sub-invariant sub-$\sigma$-algebra.   Let $\xi$ and $\zeta$ be measurable partitions.  Then $$
H_\mu(\zeta \vee\xi \mid T\inv( \zeta ^- \vee \xi^-)\vee \sigG)\le  H_\mu(\zeta \mid T\inv \zeta ^-\vee \sigG) + H_{\mu}(\xi\mid \zeta\vee \sigG).$$  Moreover, 
if $H_{\mu}(\xi\mid \zeta\vee\sigG)<\infty$ and if $ H_\mu(\zeta \mid T\inv \zeta ^-\vee \sigG) <\infty$ then 
$$
H_\mu(\zeta \mid T\inv \zeta ^-\vee \sigG)\le H_\mu(\zeta \vee\xi \mid T\inv( \zeta ^- \vee \xi^-)\vee \sigG).$$
\end{claim}
\begin{proof}
We have 
\begin{align}
H_\mu(\zeta \vee&\xi \mid T\inv( \zeta ^-  \vee \xi^-)\vee \sigG )\notag \\
&=H_\mu(\zeta  \mid T\inv \zeta ^- \vee T\inv \xi^-\vee \sigG )+ H_\mu(\xi \mid  \zeta ^- \vee T\inv \xi^-\vee \sigG ).\label{eq:rewrite}
\end{align}
The upper bound on $H_\mu(\zeta \vee\xi \mid T\inv( \zeta ^- \vee \xi^-)\vee \sigG )$ follows since \eqref{eq:rewrite} is bounded above by 
$H_\mu(\zeta  \mid T\inv \zeta ^- \vee \sigG )+ H_\mu(\xi \mid  \zeta \vee\sigG).$

For the lower bound on $H_\mu(\zeta \vee\xi \mid T\inv( \zeta ^- \vee \xi^-)\vee \sigG)$,
we have $H_\mu( \xi^{[1,m)}   \mid T\inv \zeta ^- \vee \sigG )\le(m-1) H (\xi\mid \zeta\vee \sigG)<\infty$
and using that   $ H_\mu(\zeta \mid T\inv \zeta ^-\vee \sigG) <\infty$ we  estimate \eqref{eq:rewrite} as 
\begin{align*}
H_\mu&(\zeta  \mid T\inv \zeta ^- \vee T\inv \xi^-\vee \sigG )+ H_\mu(\xi \mid  \zeta ^- \vee T\inv \xi^-\vee \sigG )\\
& = \lim_{m\to \infty} \bigg(H_\mu(\zeta  \mid T\inv \zeta ^- \vee \xi^{[1,m)}\vee \sigG )+ H_\mu(\xi \mid  \zeta ^- \vee \xi^{[1,m)}\vee \sigG )\bigg) \\
&=\lim \bigg(H_\mu(\zeta\vee \xi^{[1,m)}   \mid T\inv \zeta ^- \vee \sigG )-H_\mu( \xi^{[1,m)}   \mid T\inv \zeta ^- \vee \sigG )
\\
&\quad \quad + H_\mu(\xi \mid  \zeta ^- \vee \xi^{[1,m)}\vee \sigG )\bigg) \\
&=\lim \bigg(H_\mu(\zeta  \mid T\inv \zeta ^- \vee \sigG )+ H_\mu(\xi^{[1,m)}  \mid \zeta ^- \vee \sigG )\\
&\quad \quad -H_\mu( \xi^{[1,m)}   \mid T\inv \zeta ^- \vee \sigG )
+ H_\mu(\xi \mid  \zeta ^- \vee \xi^{[1,m)}\vee \sigG )\bigg) \\
&=\lim \bigg(H_\mu(\zeta  \mid T\inv \zeta ^- \vee \sigG )+ H_\mu(\xi^{[0,m)}  \mid \zeta ^- \vee \sigG )\\
&\quad \quad -H_\mu( \xi^{[1,m)}   \mid T\inv \zeta ^- \vee \sigG )
\bigg) \\
&=\lim \bigg(H_\mu(\zeta  \mid T\inv \zeta ^- \vee \sigG )+ H_\mu(T\inv \xi^{[0,m)}  \mid T\inv \zeta ^- \vee T\inv \sigG )\\
&\quad \quad -H_\mu( \xi^{[1,m)}   \mid T\inv \zeta ^- \vee \sigG )
\bigg) \\
&=\lim_{m\to \infty}  \bigg(H_\mu(\zeta  \mid T\inv \zeta ^- \vee \sigG )+ H_\mu( \xi^{[1,m+1)}  \mid T\inv \zeta ^- \vee T\inv \sigG )\\
&\quad \quad -H_\mu( \xi^{[1,m)}   \mid T\inv \zeta ^- \vee \sigG )
\bigg) \\
&\ge H_\mu(\zeta  \mid T\inv \zeta ^- \vee \sigG ) 
\end{align*}
since 
\begin{align*}
H_\mu( \xi^{[1,m+1)}  &\mid T\inv \zeta ^- \vee T\inv \sigG )- H_\mu( \xi^{[1,m)}   \mid T\inv \zeta ^- \vee \sigG )\\&\ge 
H_\mu( \xi^{[1,m+1)}  \mid T\inv \zeta ^- \vee  \sigG ) -H_\mu( \xi^{[1,m)}   \mid T\inv \zeta ^- \vee \sigG )\ge 0.\qedhere\end{align*}

\end{proof}

\subsubsection{Properties of transverse conditional entropy}
From definition and from \cref{claim:mylifeboilsdowntoinfinityminusinfinity}, we have the following facts.  
\begin{proposition}[Properties of transverse conditional entropy]
Let $\sigF$ and $\sigG$ be  sub-invariant sub-$\sigma$-algebras and let $\xi$ and $\zeta$ be measurable partitions.  
\begin{enumerate}
\item  $\transent_\mu(T,\xi\mid \sigF)\le H_\mu(\xi\mid \sigF)$.
\item  If $\sigF\subset \sigG$ then  $\transent_\mu(T,\xi\mid \sigG)\le \transent_\mu(T,\xi\mid \sigF).$
\item   If  $H_\mu(\xi\mid \zeta)<\infty$ and if $\transent_\mu(T,\zeta \mid \sigF)<\infty,$ 
then $$\transent_\mu(T,\zeta\mid \sigF)\le \transent_\mu(T,\zeta \vee \xi\mid \sigF)\le 
\transent_\mu(T,\zeta\mid \sigF) + H_{\mu}(\xi\mid \zeta \vee \sigF).$$  
\item   If $H_\mu(\xi\mid \sigF)<\infty$ then $\transent_\mu(T,\xi\mid \sigF) = \lim _{n\to \infty}    H_\mu \left( \xi \mid  \xi^{[1,n)} \vee \sigF\right)$.
\end{enumerate}

\end{proposition}

For automorphisms,   a similar result to \cref{lem:natext} controls (sub)additivity of metric entropy.  
\begin{proposition}[{cf.\  \cite[Lem.\ 6.1]{MR2729332}, \cite[Thm.\ 6.3]{MR0262464}, \cite[(7.7)]{MR0217258}}] \label{Prop:GenARpart}
Suppose $T$ is an automorphism.  
For any measurable partitions $\xi $ and $ \eta$ of $X$ with  $\transent_\mu(T,\xi\mid \eta^{-})<\infty$, we have 
$$ h_\mu(T,\xi\vee\eta)\le \transent_\mu(T, \xi\mid  \eta^\Z)+ h_\mu(T, \eta).$$
Moreover, if   $H_\mu(\xi\vee \eta  \mid T\inv \eta^-)<\infty$  then  
\begin{equation} \label{eq:GAR} 
  h_\mu(T,\xi\vee\eta)=   \transent_\mu(T, \xi\mid  \eta^\Z)+  h_\mu(T, \eta).\end{equation}
\end{proposition}

\begin{proof}Following the proof of \cite[Lem.\ 6.1]{MR2729332}, 
$$ h_\mu(T, \xi\vee \eta)
=   h_\mu(T, \xi\vee T^n\eta)=  H_\mu( \xi\mid T\inv \xi^{- }\vee T^{n-1}\eta^{-})
+ H_\mu( \eta\mid T^{-n}\xi^{- }\vee T^{-1}\eta^{-}).$$
As $n\to \infty$, $ H_\mu( \xi\mid T\inv \xi^{- }\vee T^{n-1}\eta^{-})\searrow \transent_\mu(T, \xi\mid  \eta^\Z).$
Moreover, $H_\mu( \eta\mid T^{-n}\xi^{- }\vee T^{-1}\eta^{-})\le h_\mu(T, \eta)$.   By \cite[(7.6)]{MR0217258},   if $H_\mu(\xi\vee \eta  \mid T\inv \eta^-)<\infty,$ then as $n\to \infty$, \[H_{\mu}( \eta\mid T^{-n}\xi^{- }\vee T^{-1}\eta^{-})\nearrow   h_\mu(T, \eta).\qedhere\] 
\end{proof}

\subsubsection{Computing entropies}\label{comps}
We observe that it suffices to take the suprema in all definitions of entropy in \cref{fiberentdef} to be over finite measurable partitions.  Indeed, given a measurable partition $\xi$, there is a sequence of finite partitions $\xi_n$ with $\xi_n\nearrow \xi$.  
If $H_\mu(\xi\mid \sigF)<\infty$ then $H_\mu(\xi\mid \xi_n \vee\sigF)\to 0$ and thus $h_\mu(T,\xi\mid \sigF) = \lim h_\mu(T,\xi_n\mid \sigF)$ follows from \cref{prop:standard properties}\eqref{standard2}.
Similarly, $$H_\mu(\xi_n \mid T\inv \xi_n^-\vee \sigF)\ge H_\mu(\xi_n \mid T\inv \xi^-\vee \sigF)\nearrow H_\mu(\xi \mid T\inv \xi^-\vee \sigF)$$
and 
\begin{align*}H_\mu(\xi_n\vee \Xi(\sigF) \mid T\inv( \xi_n^-\vee \Xi(\sigF)))&\ge H_\mu(\xi_n\vee \Xi(\sigF) \mid T\inv (\xi^-\vee \Xi(\sigF)))\\&\nearrow H_\mu(\xi \vee \Xi(\sigF)\mid T\inv (\xi^-\vee \Xi(\sigF))).\end{align*}

In particular, given a sequence  of finite partitions $\xi_n$ with $\xi_n\nearrow \ptpart $ or, more generally,  with $\xi_n\vee \Xi(\sigF)\nearrow \ptpart $, by using  \cref{prop:standard properties}\eqref{standard2} or \cref{claim:mylifeboilsdowntoinfinityminusinfinity}, we may approximate all three entropies relative to an arbitrary finite partition by the entropy relative to some $\xi_n$.  In particular, 
\begin{enumerate}
\item $h_\mu(T\mid \sigF) = \lim _{n\to \infty} h_\mu(T,\xi_n\mid \sigF)$;
\item $\transent_\mu(T\mid \sigF) = \lim _{n\to \infty} \transent_\mu(T,\xi_n\mid \sigF)$;
\item $\maxent_\mu(T\mid \sigF) = \lim _{n\to \infty} h_\mu(T,\xi_n\vee\Xi( \sigF))$.
\end{enumerate}

\subsubsection{SMB-type lemma}
We record the following Shannon-McMillan-Breiman-type formula for metric entropy that is helpful in computing $\maxent_\mu(T\mid \sigF)$.  The proof   essentially follows   \cite[Lem.\  9.3.1]{MR819557}; we provide most of the details.    
\begin{lemma}\label{lemma:pastfinitepartabstract} 
Suppose that $T\colon (X,\mu)\to (X,\mu)$ is an ergodic, measurable automorphism. 
Let $\zeta$ be a measurable sub-invariant partition of $(X,\mu)$ with $h_\mu(T, \zeta)<\infty$ and let  $ \partP$ be a measurable partition of $(  X,  \mu)$ with $H_\mu( \partP\mid  \zeta)<\infty$. 
Then for $ \mu$-a.e.\ $x$,
$$\lim_{n \to \infty} - \frac 1 n \log  { \mu }^{  \zeta}_{x} \left( (    \zeta \vee   \partP) ^{(-n,0]}(x) \right)  = h_{ \mu }(   T, \partP \vee \zeta ).$$
\end{lemma}
\begin{proof}
Given measurable partitions $\xi$ and $\eta$, write  $I(\xi\mid \eta)(x) = -\log\mu^\eta_x(\xi(x))$.  Then 
\begin{align*}
- \log  { \mu }^{  \zeta}_{x} &\left( (    \zeta \vee   \partP) ^{(-n,0]}(x) \right) \\
 &=
 I(\partP\mid \zeta)(x) + 
 \sum_{k=1}^{n-1} I(\partP\vee \zeta\mid T^{-1}\zeta \vee T\inv \partP^{[0,k)})(T^{-k}(x)).
\end{align*}

For $k\in\{0,1,2,\dots\}$, let 
 \begin{align*}
 g_k(x) &:= I\bigl(\zeta\vee \partP\ \mid  T^{-1}\zeta \vee T\inv \partP^{[0,k)} \bigr)(x) \\
 g_*(x) &:= \sup\left\{I\bigl(\zeta\vee \partP\ \mid  T^{-1}\zeta \vee T\inv \partP^{[0,k)} \bigr)(x) : k\in \{0,1,2,\dots\}\right\}\\
 g_\infty(x)& :=I\bigl(\zeta\vee \partP\ \mid  T^{-1}\zeta \vee T\inv \partP^{-} \bigr)(x). 
  \end{align*}
By the Chung-Neveu maximal lemma (see \cite[Lem.\  2.1]{MR0262464}), we have $g_*\in L^1(\mu)$ and thus by martingale convergence,  we have  a.e.\ and  $L^1$ convergence $g_k \to g_\infty$.  
By ergodicity of $T$ and by  applying  the ergodic theorem for sequences of functions (see \cite{MR2028}), we have 
\[
\lim_{n\to \infty}   - \frac{1}{n}\log  { \mu }^{  \zeta}_{x} \left( (    \zeta \vee   \partP) ^{(-n,0]}(x) \right) =\int  g_\infty \, d \mu = h_\mu(T,\zeta\vee \partP). \qedhere \]
\end{proof}

Let  $T\colon (X,\mu)\to (X,\mu)$ be an ergodic, measurable automorphism and let  $\sigF$ be a sub-invariant sub-$\sigma$-algebra with $h_\mu(T,\Xi(\sigF))<\infty$. 
We will say a measurable partition $\partP$  is a  (two-sided) \emph{relative generator} (relative to $\sigF$) if we have $H_{\mu}(\partP\mid \sigF)<\infty$ and $\partP^\Z\vee \sigF=\ptpart$.  
If $\partP$ is a 
relative generator then 
$$ h_{ \mu }(   T,\partP \vee \Xi(\sigF))=\maxent_{ \mu }(   T \mid \sigF).$$
Indeed, if $\partQ$ is any measurable partition with $H_\mu(\partQ\mid \sigF)<\infty$ and if $\partP^\Z\vee \sigF=\ptpart$, then by 
  \cref{Prop:GenARpart} we have
\begin{align*}h_\mu(T, \partQ\vee \Xi(\sigF))&\le h_\mu(T, \partQ\vee \partP \vee \Xi(\sigF))\\&\le 
h_\mu(T, \partP\vee \Xi(\sigF)) + \transent_\mu(T, \partQ\mid ( \partP\vee \Xi(\sigF))^\Z)\\&\le 
h_\mu(T, \partP\vee \Xi(\sigF)) + \transent_\mu(T, \partQ\mid  \partP^\Z \vee \Xi(\sigF))\\
&= h_\mu(T, \partP\vee \Xi(\sigF)). 
\end{align*}

\begin{corollary}\label{cor:SMBgen}
Let  $T\colon (X,\mu)\to (X,\mu)$ be an ergodic, measurable automorphism and let  $\sigF$ be a sub-invariant sub-$\sigma$-algebra with $h_\mu(T,\Xi(\sigF))<\infty$.  If $\partP$ is a 
relative generator then for $\mu$-a.e.~$x$,
$$\lim_{n \to \infty} - \frac 1 n \log  { \mu }^{  \sigF}_{x} \left( (    \sigF \vee   \partP) ^{(-n,0]}(x) \right)  =\maxent_{ \mu }(   T \mid \sigF).$$
\end{corollary}

\subsubsection{Relation to i.i.d.~random walks}\label{computefiberent}
Let $M$ be a compact metric space, let
$\nu_0$ be a probability measure on $C^{0}(M)$,  and let $\mu_0$ be a $\nu_0$-stationary probability measure.  
Let $\what \Sigma=\supp(\nu_0)^\N$ and  let $\what F\colon \widehat \Sigma\times M\to  \widehat \Sigma\times M$ denote  the skew extension preserving the probability measure $\hat \mu = \nu_0^\N\times \mu_0.$
Let $\what \sigF$ denote the sigma algebra on $\what X:= \widehat \Sigma\times M$ given by preimages of the map $\what X\to \widehat \Sigma$. 

Given a measurable partition $\xi$ of $M$, let $\wtd \xi$ denote the partition by atoms of the form $\what \Sigma \times C$ for $C\in \xi$.  If $\xi_n$ is a sequence of finite partitions of $M$ with $\xi_n\nearrow \ptpart_M$ we have $\wtd \xi_n \vee \Xi(\what \sigF)\nearrow \ptpart_{\what \Sigma\times M}.$
In particular, from the discussion in \cref{comps}, we have
\begin{enumerate}
\item $h_{\mu_0}(\nu_0)=\lim h_{\mu_0}(\nu_0,\xi_n)=\lim h_{\hat \mu}(\what F,\wtd\xi_n\mid \what \sigF)=h_{\hat \mu}(\what F\mid \what \sigF),$ and
\item $\transent_{\mu_0}(\nu_0)=\lim \transent_{\mu_0}(\nu_0,\xi_n)=\lim \transent_{\hat \mu}(\what F,\wtd\xi_n\mid \what \sigF)=\transent_{\hat \mu}(\what F\mid \what \sigF).$
\end{enumerate}

\subsection{(Non-)naturality of entropies under iteration}
We record properties of various notions of entropy defined in \cref{fiberentdef} under iteration.  

First, we have the usual naturality of  standard conditional entropy   under iteration.  
\begin{proposition}[{\cite[Fact 2.4.19]{MR2809170}}]
For $n\ge 0$ we have $h_\mu(T^n \mid \sigF)= nh_\mu(T\mid \sigF)$.  

Moreover,  if $T$ is an automorphism and $\sigF$ is \siname\ then for $n\in \Z,$ we have $h_\mu(T^n \mid \sigF)= |n|h_\mu(T\mid \sigF)$.
\end{proposition}

Unlike standard conditional entropy, transverse conditional entropy is only sublinear under iteration.  See \cref{sec:examples} for examples where transverse conditional entropy is strictly sublinear under iteration.   However,  transverse conditional entropy  behaves  naturally  under a mild injectivity hypothesis or for \siname\ $\sigF$.
\begin{proposition}\label{naturalityoftransverse}
For $n\ge 0$, we have $$\transent_\mu(T^n \mid \sigF)\le  n \transent_\mu(T\mid \sigF).$$ 

Moreover, 
\begin{itemize}
\item if $\sigF$ is \siname, or 
\item if the restriction $\restrict T{\sigF(x)}\colon \sigF(x)\to \sigF(T(x))$ of $T$ to $\mu$-a.e.\ atom $\sigF(x)$ of $\sigF$ is injective a.s.,
\end{itemize} then $$\transent_\mu(T^n \mid \sigF)=n \transent_\mu(T\mid \sigF).$$

\end{proposition}
\begin{proof}
For $n\ge 1$ and any measurable partition $\xi$, we have that 
\begin{align}
\transent_\mu(T^n, \xi^{[0,n)} \mid \sigF) 
&=H_\mu( \xi^{[0,n)}\mid T^{-n}\xi^-\vee \sigF)\notag \\
&=\sum_{j=0}^{n-1} H_\mu( T^{-j}\xi \mid \xi^{[j+1,n)}\vee T^{-n}\xi^-\vee \sigF)\notag \\
&=\sum_{j=0}^{n-1} H_\mu( T^{-j}\xi \mid  T^{-j-1}\xi^-\vee \sigF). \label{keytoest}
\end{align}
We bound \eqref{keytoest} above by 
$$\sum_{j=0}^{n-1} H_\mu( T^{-j}\xi \mid  T^{-j-1}\xi^-\vee T^{-j}\sigF)=\sum_{j=0}^{n-1} H_\mu( \xi \mid  T^{-1}\xi^-\vee \sigF)= n \transent_\mu(T,\xi\mid \sigF).$$
The sub-linearity follows.  
Moreover,  if $\sigF$ is \siname\ then $T^{-j}\sigF=\sigF$ and for all $j\ge 1$.  In this case,  the bound for \eqref{keytoest} is equal to \eqref{keytoest} and thus 
$\transent_\mu(T^n, \xi^{[0,n)} \mid \sigF) =n\transent_\mu(T, \xi \mid \sigF) .$ (See also \cref{equalifinv} below.)

Finally, we may alternatively bound \eqref{keytoest} below by 
\begin{align*}
\sum_{j=0}^{n-1} H_\mu( T^{-(n-1)}\xi \mid  T^{-n}\xi^-\vee T^{-(n-j-1)}\sigF)
&\ge \sum_{j=0}^{n-1} H_\mu( T^{-(n-1)}\xi \mid   T^{-n}\xi^-\vee   \sigF)\\
&= n \transent_\mu(T,T^{-(n-1)}\xi \mid \sigF).  
\end{align*}
Take a sequence of finite measurable partitions $\xi_j$ with  $\xi_j\nearrow \ptpart$.   If  the restriction $\restrict T{\sigF(x)}\colon \sigF(x)\to \sigF(T(x))$ of $T$ to $\mu$-a.e.\ atom $\sigF(x)$ of $\sigF$ is injective a.s., then for every $n$, $$T^{-(n-1)}\xi_j \vee \sigF \to \ptpart$$
and so \[\transent_\mu(T,T^{-(n-1)}\xi_j \mid \sigF)\to \transent_\mu(T\mid \sigF).  \qedhere\]
\end{proof}

Finally, we observe that the refinement entropy relative to  $\sigF$ is natural under iteration.  
\begin{proposition}
For $n\ge 0$, we have $$\maxent_\mu(T^n \mid \sigF)=  n \maxFent.$$ 
\end{proposition}
\begin{proof}
Given any measurable partition $\xi$ with $\Xi(\sigF)\prec \xi$, we have $\Xi(\sigF)\prec \xi^{[0,n)}$.  Then  
\begin{align*}
 h_\mu(T^n, \xi^{[0,n)} )
&=H_\mu( \xi^{[0,n)}\mid T^{-n}\xi^-) 
=\sum_{j=0}^{n-1} H_\mu( T^{-j}\xi \mid \xi^{[j+1,n)}\vee T^{-n}\xi^-)\notag \\
&=\sum_{j=0}^{n-1} H_\mu( T^{-j}\xi \mid  T^{-j-1}\xi^-)= n h_\mu(T, \xi ). 
\end{align*}
By taking a sequence of finite measurable partitions with $\xi_j\nearrow \ptpart$, we also have that  $(\Xi(\sigF)\vee \xi_j)^{[0,n)}\nearrow \ptpart$  and the result follows.  
\end{proof}

\subsection{Entropies under natural extension}
Let $T\colon (X,\mu)\to (X,\mu)$ be measure preserving and let $\sigF$ be a sub-invariant sub-$\sigma$-algebra. 
\begin{proposition}\label{prop:howtoNatExt}
Let $\wtd T\colon( \wtd X,\td \mu)\to ( \wtd X,\td \mu)$ be the ergodic theoretic natural extension of $T$ and let  $\tau\colon  (\wtd X,\td \mu)\to (X,\mu)$ denote the  factor map.
Then 
\begin{enumerate}
\item\label{item:natural-extension-refinement-entropy} $\maxent_{\td \mu}(\wtd T\mid \tau\inv \sigF) = \maxent_\mu(T\mid \sigF)$.  
\item\label{item:natural-extension-relative-entropy} $ h_{\td \mu}(\wtd T\mid  \tau\inv \sigF) =  h_\mu(T\mid \sigF)$.  
\end{enumerate}
\end{proposition}
We note the corresponding result for the transverse conditional entropy  $\transent$ does not hold (see \cref{prop:GenAR,hbarhhatsame} below).
\begin{proof}[Proof of \cref{prop:howtoNatExt}]
Write $ \wtd \sigF := \tau\inv \sigF$.
Let $\xi$ and $\eta$ be finite partitions of $\wtd X$ and $X$, respectively.  Let $\zeta= \Xi(\sigF)$ and let $\wtd \eta=\tau\inv  \eta$ and $\wtd \zeta=\tau\inv \zeta$.

For conclusion \eqref{item:natural-extension-refinement-entropy}, we have 
$$h_\mu(T, \eta \vee \zeta)= h_{\td \mu}(\wtd T, \wtd \eta \vee \wtd\zeta)\le 
  \maxent_{\td \mu}(\wtd T\mid  \wtd \sigF).$$
Taking suprema over all   $\eta$, we have $\maxent_\mu(T\mid \sigF)\le \maxent_{\td \mu}(\wtd T\mid  \wtd \sigF).$
 In particular, if $\maxent_\mu(T\mid \sigF)=\infty$, equality holds.
 
We thus assume $\maxent_\mu(T\mid \sigF)<\infty$.  Applying \cref{Prop:GenARpart}, we have
\begin{align*}
h_{\td \mu}(\wtd T, \xi\vee \wtd \zeta) &\le h_{\td \mu}(\wtd T, \xi\vee\wtd \eta \vee  \wtd \zeta)\\
&\le \transent_{\td \mu}(\wtd T, \xi\mid  (\wtd \eta \vee  \wtd \zeta)^\Z)+ h_{\td \mu}(\wtd T, \wtd \eta\vee \wtd \zeta)\\&
= \transent_{\td \mu}(\wtd T, \xi\mid  \wtd \eta^\Z\vee \wtd\sigF^+ )+  h_{\td \mu}(\wtd T, \wtd \eta\vee \wtd \zeta)
\end{align*}
Let $\ptpart$ denote the point partition on $ X$.  
Since the above holds for any finite partition $\eta$ of $X$, by taking a sequence $\eta_n$ of finite measurable partitions of $X$ with $\eta_n \nearrow \ptpart$, we have 
\begin{align*}h_{\td \mu}(\wtd T, \xi\vee \wtd \zeta) &\le \transent_{\td \mu}(\wtd T, \xi\mid  (\tau\inv \ptpart)^+)+ \sup_n h_{\td \mu}(\wtd T, \wtd \eta_n\vee \wtd \zeta)\\
&= \transent_{\td \mu}(\wtd T, \xi\mid  (\tau\inv \ptpart)^+)+ \sup_n h_\mu(T,  \eta_n\vee  \zeta)\\
&\le 0+ \maxent_\mu(T\mid \sigF)
\end{align*}
since $(\tau\inv \ptpart)^+$ is the point partition on $\wtd X$.  
Taking the supremum over all $\xi$, we have $$\maxent_{\td \mu}(\wtd T\mid  \wtd \sigF)\le\maxent_\mu(T\mid \sigF).$$

For \eqref{item:natural-extension-relative-entropy}, 
 we similarly have 
$$h_\mu(T, \eta \mid \sigF)= h_{\td \mu}(\wtd T, \wtd \eta \mid  \wtd \sigF)\le 
  h_{\td \mu}(\wtd T\mid  \wtd \sigF)$$
and so we have $  h_\mu(T\mid \sigF)\le   h_{\td \mu}(\wtd T\mid  \wtd \sigF).$
Applying  \cref{lem:natext},  we obtain 
\begin{align*}
h_{\td \mu}(\wtd T,\xi  \mid  \wtd \sigF) &\le
h_{\td \mu}(\wtd T,\xi \vee \wtd \eta \mid  \wtd \sigF) \\
&= h_{\td \mu}(\wtd T,\wtd \eta  \mid  \wtd \sigF)  + h_{\td \mu}(\wtd T,\xi \mid  \wtd \eta^\Z \vee  \wtd \sigF^+).
\end{align*}
Again taking $\eta_n \nearrow \ptpart$ we have 
$$h_{\td \mu}(\wtd T,\xi  \mid  \wtd \sigF) \le 
 \sup_n h_{\td \mu}(\wtd T, \wtd \eta_n  \mid  \wtd \sigF) + h_{\td \mu}(\wtd T,\xi \mid ( \tau \inv \ptpart)^\Z \vee  \wtd \sigF^+) = h_\mu(  T \mid  \sigF)$$
and so $h_{\td \mu}( \wtd T \mid   \wtd \sigF)\le h_\mu(  T \mid  \sigF)$.  
\end{proof}

\subsection{Relations among  entropies} 
For a  measure-preserving transformation  $T\colon (X,\mu)\to (X,\mu)$ and a sub-invariant sub-$\sigma$-algebra $\sigF$, we have the following elementary inequalities relating   the entropies defined in \cref{fiberentdef}.

\subsubsection{Equivariant family of conditional measures and basic inequalities}  

We begin with the following definition.  
\begin{definition}Let $\sigF$ be a sub-invariant sub-$\sigma$-algebra.  
We say that $\sigF$ admits an \emph{equivariant family of conditional measures} if for $\mu$-a.e.\ $x$,
$$T_*\mu_x^\sigF = \mu_{T(x)}^\sigF.$$
\end{definition}

We remark that if $\sigF$ admits an {equivariant family of conditional measures} then for any set $C$ and almost every $x$,
$$\mu_x ^\sigF (T\inv C) = \mu^\sigF_{T(x)} (C).$$
In particular,  if $\sigF$ admits an {equivariant family of conditional measures} then, for any 
measurable partition $\xi$ of $(X, \mu)$, we have 
\begin{equation}\label{claim:equivar}
H_\mu(T\inv \xi \mid T\inv  \sigF ) =H_\mu(\xi \mid \sigF ) = H_\mu(T\inv \xi \mid \sigF)
\end{equation}
where the first equality is from property \eqref{keyentropy2} of \cref{entfactsbasic}.

\begin{proposition}\label{prop:basicineq}
Let $\xi$ be a measurable partition of  $(X,\mu)$ with $H_\mu(\xi\mid \sigF)<\infty$.   Then 
\begin{equation}\label{basicineq}  h_\mu(T,\xi\mid \sigF)\le  \transent_\mu(T,\xi\mid \sigF)\le h_\mu(T, \xi\vee \Xi(\sigF)).\end{equation}
In particular,
$$  h_\mu(T\mid \sigF)\le  \transent_\mu(T\mid \sigF)\le  \maxent_\mu(T\mid \sigF).$$

Moreover, if
 $\sigF$ admits an equivariant family of conditional measures then   $$  h_\mu(T,\xi\mid \sigF)=  \transent_\mu(T,\xi\mid \sigF) $$ 
 for all  measurable partitions $\xi$  with $H_\mu(\xi\mid \sigF)<\infty$ whence $h_\mu(T\mid \sigF)=  \transent_\mu(T\mid \sigF) .$
\end{proposition}
\begin{proof}
Fix   $\xi$ with $H_\mu(\xi\mid \sigF)<\infty$.  For the right-hand bound in \eqref{basicineq}
we have 
\begin{align*}
H_\mu(\xi\vee &\Xi(\sigF)\mid T^{-1} \xi^-\vee T^{-1}\Xi(\sigF))\\&= 
 H_\mu(\Xi(\sigF)\mid   T\inv \xi^-\vee T\inv \Xi(\sigF))
+H_\mu(\xi\mid T^{-1} \xi^-\vee \Xi(\sigF))
 \\
&\ge  0+
 \transent_\mu(T,\xi \mid \sigF).\end{align*}

 For the left-hand bound in \eqref{basicineq}, since $T\inv \sigF\subset \sigF$ we have 
 \begin{align}
h_\mu(T,\xi\mid \sigF)
&=
\lim_{n\to \infty} \frac 1 n H\left(\xi^{[0,n)} \mid \sigF\right)\notag\\
&=
\lim_{n\to \infty} \frac 1 n\left(
 H_\mu(\xi \mid  T\inv \xi^{[0,n-1)}\vee \sigF)+
  H_\mu(T\inv \xi^{[0,n-1)} \mid \sigF)\notag
 \right) 
\\
&= 
\lim_{n\to \infty} \frac 1 n\sum_{j=0}^{n-1} 
  H_\mu(T^{-j}\xi \mid  T^{-j-1} \xi^{[0,n-j-1)}\vee  \sigF)\notag
\\
&\le 
\lim_{n\to \infty} \frac 1 n\sum_{j=0}^{n-1} 
  H_\mu(T^{-j}\xi \mid  T^{-j-1} \xi^{[0,n-j-1)}\vee T^{-j} \sigF)\label{thisone}
\\&= 
\lim_{n\to \infty} \frac 1 n\sum_{j=0}^{n-1} 
 H_\mu(\xi \mid  T^{-1} \xi^{[0,n-j-1)}\vee \sigF)\notag
\\
&=  H_\mu(\xi \mid  T^{-1} \xi^-\vee \sigF)=:\transent_\mu(T,\xi\mid \sigF)\label{final hit}
\end{align}
where \eqref{final hit} is a standard estimate (cf.\ \cite[pf.\ of Lem.\ II.1.9]{MR884892}).
Moreover, if $\sigF$ admits an equivariant family of conditional measures, then the inequality in \eqref{thisone} becomes an equality by \cref{claim:equivar}. 
\end{proof}
Note that if $\sigF$ is \siname\ for $T\colon X\to X$, then $\sigF$ admits an equivariant family of conditional measures.  We then obtain the following.  

\begin{corollary}[{\cite[Fact 2.3.4]{MR2809170}}]\label{equalifinv}
Suppose  that $\sigF$ is \siname\ for $T\colon X\to X$.  
Then for any $\xi$ with $H_\mu(\xi\mid \sigF)<\infty$ we have  $$\transent_\mu(T,\xi\mid \sigF)= h_\mu(T,\xi\mid \sigF).$$
In particular, if $\sigF$ is \siname\ for $T\colon X\to X$ then 
$$\maxFent  =  \transent_\mu(T\mid \sigF)= h_\mu(T\mid \sigF).$$
\end{corollary}
For the final conclusion in \cref{equalifinv}, note for any partition $\xi$ that 
\begin{align*}
h_\mu( T, \xi\vee \Xi(\sigF)) 
&= H_\mu( \xi\vee \Xi(\sigF) \mid T\inv(\xi^-\vee \Xi(\sigF))) 
\\&=H_\mu( \xi\vee \Xi(\sigF) \mid T\inv(\xi^-)\vee \sigF)  \\& =H_\mu( \xi  \mid T\inv(\xi^-)\vee \sigF)  \\&= \transent_\mu(T,\xi\mid \sigF). 
\end{align*}

\subsubsection{Decomposition and factor theorem}
Combining the above results, 
we obtain the following relating all three notions of entropy.  
\begin{proposition}\label{prop:GenAR}Suppose that 
$\sigF$ is a sub-invariant sub-$\sigma$-algebra with   $h_\mu(T, \Xi(\sigF))<\infty$.  Then
\begin{align}\maxFent&=\transent_\mu(T\mid \sigF) +  H_\mu( \sigF\mid T\inv \ptpart)\label{eq:factAR}\\
&= h_\mu(T\mid \sigF)+ h_\mu(T, \Xi(\sigF))\label{eq:standAR}.\end{align}
\end{proposition}
\begin{proof}
We first consider \eqref{eq:factAR}.  For any measurable partition $\xi$, we have
\begin{align*}
h_\mu(T,\xi\vee\Xi(\sigF)) 
&=H_\mu(\xi\vee \sigF \mid T\inv \xi^-\vee T\inv\sigF)\\
&=H_\mu( \sigF\mid T\inv \xi^-\vee T\inv\sigF)+ H_\mu( \xi \mid T\inv \xi^-\vee \sigF).\\
&=H_\mu( \sigF\mid T\inv \xi^-\vee T\inv\sigF)+ \transent_\mu(T, \xi \mid  \sigF).
\end{align*}
Taking a sequence $\xi_n$ of finite measurable partitions with $\xi_n\nearrow \ptpart$, we have 
\begin{enumerate}
\item $\transent_\mu(T,\xi_n\mid\sigF) \nearrow \transent_\mu(T\mid\sigF) $ 
\item $H_\mu( \sigF\mid T\inv \xi_n^-\vee T\inv\sigF)\searrow H_\mu( \sigF\mid T\inv \ptpart)$.
\end{enumerate}
We thus have
$$\maxent_\mu(T\mid \sigF)=  \transent_\mu(T\mid \sigF)+H_\mu( \sigF\mid T\inv \ptpart).$$

To derive \eqref{eq:standAR}, note by \cref{prop:howtoNatExt} all three terms are preserved by passing to the natural extension.   In particular,   we may assume $T$ is an automorphism.  
Applying the special case \cref{eq:GAR} of  \cref{Prop:GenARpart}  and \cref{natext,equalifinv}, for a measurable partition $\xi$ with $H_\mu(\xi\mid \sigF)<\infty$ we have
\begin{align*}
h_\mu(T,\xi\vee\Xi(\sigF)) &= \transent_\mu(T,\xi\mid \sigF^+) +  h_\mu(T,\Xi(\sigF)) \\
&=  h_\mu(T,\xi\mid \sigF^+) +  h_\mu(T,\Xi(\sigF)) \\
&=  h_\mu(T,\xi\mid \sigF) +  h_\mu(T,\Xi(\sigF)).
\end{align*}
The result follows taking the supremum over all such $\xi$.  
\end{proof}

When  $T$ is an automorphism, we have  $T\inv \ptpart=\ptpart$ and thus obtain the following, which follows directly from  \cref{prop:GenAR} when $h_\mu(T, \Xi(\sigF))<\infty$.  
If $h_\mu(T, \Xi(\sigF))=\infty$, the result also follows  since $$\transent_\mu(T, T(\Xi(\sigF))\mid \sigF)= h_\mu(T, T(\Xi(\sigF))) =h_\mu(T, \Xi(\sigF)).$$
\begin{corollary}\label{hbarhhatsame}
Suppose that $T$ is an automorphism of $(X,\mu)$.  Then for any sub-invariant sub-$\sigma$-algebra $\sigF$, $$\maxFent  =  \transent_\mu(T\mid \sigF).$$
\end{corollary}

In the sequel, we will be interested especially in the case of non-invertible $T$ with the property that the restriction $\restrict {T}{\sigF(x)}\colon \sigF(x) \to \sigF(T(x))$ of $T$ to almost-every atom of $\sigF$ is injective.  In this case $T\inv \ptpart\vee \Xi(\sigF)=\ptpart$ in which case the defect between $ \maxFent$ and $\transent_\mu(T\mid \sigF) $ in \eqref{eq:factAR} of \cref{prop:GenAR} becomes  \begin{equation} H_\mu( \sigF\mid T\inv \ptpart)= H_\mu( \ptpart\mid T\inv \ptpart).\label{eq:thiseq}\end{equation}

\begin{remark}
We recall that when $\what T\colon (\what X, \hat \mu)\to (\what X, \hat \mu)$ is the  factor map induced by $\sigF$, the classical Abramov--Rohlin formula (see \cite{MR0140660,MR1179170,MR0476995}) states (following our notational  convention in \eqref{eq:factorent}) that 
$$h_\mu(T) = h_\mu(T\mid \sigF) + h_{\hat \mu}(\what T)
 = h_\mu(T\mid \sigF) + h_{ \mu}( T,\sigF).$$
We view  \cref{eq:standAR}  of \cref{prop:GenAR} as an analogue of the classical Abramov-Rohlin formula for the quantity $\maxFent$.  Rather than the entropy of the factor system that appears in the Abramov-Rohlin formula, the term replacing the factor entropy is the quantity  $h_\mu(T,\Xi(\sigF))$,  the entropy with respect to the partition into atoms of $\sigF$.
\end{remark}

\subsection{Kullback–Leibler divergence}
To define Furstenberg entropy, we first define and discuss properties of the Kullback–Leibler divergence.
We temporarily forget any dynamical setting and consider an abstract Polish space $X$ equipped with two  probability measures  $\nu$ and $\eta$.

\subsubsection{Lebesgue--Radon--Nikodym decomposition}
The  Lebesgue--Radon--Nikodym decomposition  of $\nu$ relative to $\eta$ is the decomposition   $$\nu= \rho \eta + \nu_s$$ where $ \rho :=\frac{d\nu}{d\eta}\in L^1( \eta)$  is the Radon--Nikodym derivative of the component of $\nu$ absolutely continuous with respect to $\eta$   and $\nu_s\perp \eta$ is the component of $\nu$ singular with respect to $\eta$.  
 We often use the following fact: When $\nu\ll \eta$ (equivalently when $\nu_s(X)=0$) then $1/\rho  \in L^1(\nu)$ and $\eta = \bigl(\frac{1}{ \rho }\bigr)\nu  + \restrict{\eta}{\rho\equiv 0}.$
 In particular, if $\nu\ll \eta$,  the equality 
$\frac{d\eta}{d\nu} = \frac{1}{\rho}  =\left(\frac{d\nu}{d\eta} \right)^{-1}$ holds $\nu$-a.e.

\subsubsection{Kullback--Leibler divergence}
Note the positive part of the function $\log \bigl({\frac{d\eta}{d\nu}}\bigr) $ is clearly $L^1(\nu)$.  
We then define the \emph{Kullback--Leibler divergence} between $\nu$ and $\eta$ to be the quantity 
\begin{equation}\label{eq:KLD1}\KL(\nu \parallel \eta)= -\int \log \left({\frac{d\eta}{d\nu}}\right) \, {d\nu} .\end{equation}
By monotonicity and concavity of $t\mapsto \log(t)$, Jensen's inequality implies
\begin{align*}
- \KL(\nu \parallel \eta) &=\int \log \left({\frac{d\eta}{d\nu}}\right) \, {d\nu}
\le  \log \left( \int{\frac{d\eta}{d\nu}} \, {d\nu}\right) 
\le \log \left(\eta(X)\right) \le \log(1)
\end{align*}
and so \begin{equation}\label{eq:KLpositive}\KL(\nu \parallel \eta)\ge 0.\end{equation}
Moreover,  strict concavity of $t\mapsto \log(t)$ implies \begin{equation}\label{eq:KLzero} \text{$\KL(\nu \parallel \eta)=0$ if and only if $\frac{d\eta}{d\nu}\equiv 1$ $\nu$-a.s. if and only if $\nu=\eta$.  }\end{equation}

Note that if $\nu= \tfrac{d\nu}{d\eta} \eta + \nu_s$ with $\nu_s(X)>0$, then $\frac{d\eta}{d\nu} = 0$ on a set of positive $\nu$-measure in $X$.  In this case,  $-\log\left({\frac{d\eta}{d\nu}}\right)=+\infty$ on a  $\nu$-positive measure set and so $\KL(\nu \parallel \eta)=+\infty$. 
 On the other hand, if $\nu_s(X) =0$, then $$-\log \left({\frac{d\eta}{d\nu}}\right) =  \log\left({\frac{d\nu}{d\eta}}\right)$$ holds $\nu$-a.e. In particular, definition \eqref{eq:KLD1}  is equivalent to the more standard definition: 
\begin{equation}\label{eq:KLD2}\KL(\nu \parallel \eta)= \begin{cases}
\displaystyle\int \log \left({\frac{d\nu}{d\eta}}\right) \, {d\nu}  & \nu\ll \eta\\
+\infty & \text{else}.
\end{cases}\end{equation}

To justify the integral in \eqref{eq:KLD2} and later $L^1$ considerations, note in the case that $ \nu\ll \eta$ that the negative part of $\log \bigl({\frac{d\nu}{d\eta}}\bigr) $ is always $L^1(\nu)$. Indeed if $\rho=  {\frac{d\nu}{d\eta}} $,  the negative part of $\log \rho$ has integral bounded by
$$\int _{0<\rho<1} -\log \rho \, d\nu= \int _{0<\rho<1} -\rho \log \rho \, d\eta \le \frac 1 e.$$
We alternate between definitions \eqref{eq:KLD1} and \eqref{eq:KLD2} freely.

\def\1{\mathbbm 1}

\subsubsection{Convexity properties of Kullback--Leibler divergence} 
Using convexity of the function $\phi\colon t\mapsto t\log t$, we obtain the following facts.
\begin{lemma}\label{lem:superbconvexity}
Let $X,Z$ be Polish spaces, let $\nu,\eta$ be probability measures on $X$, and let $\lambda$ be a probability measure on $Z$. Let $S\colon X\to Z$ be a function that is measurable with respect to both $\nu$ and $\eta$.
Let $z\mapsto \nu_z$ be a measurable map from $Z$ to the space of probability measures on $X$ such that $\nu = \int \nu_z \, d\lambda(z)$.  
Then the following hold:
\begin{enumerate}[label*={(\alph*)}]
\item\label{convex1} $\KL(S_*\nu \parallel  S_*\eta)\le \KL(\nu \parallel \eta)$; moreover if $S$ is injective when restricted to a set $X'\subset X$ with $\nu(X')=\eta(X')=1$, then equality holds.
\item\label{convex2} $\displaystyle\KL(\nu \parallel \eta)\le \int \KL(\nu_z \parallel \eta)\, d \lambda(z)$.
\end{enumerate}
\end{lemma}
\begin{proof}
Inequalities in the lemma hold trivially whenever the righthand sides are infinite.  Also, equality in  \ref{convex1} holds if $S\colon X\to Z$ is  injective on a set of full measure with respect to both $\nu$ and $\eta$ and if $\nu$ has a nontrivial singular component with respect to $\eta$ (in which case both sides are $+\infty$).  
Thus, we may assume  the following: In \ref{convex1}, we  suppose $\nu\ll\eta$ which then implies $S_*\nu\ll S_*\eta$;  in  \ref{convex2}, we  suppose $\nu_z\ll \eta$ for $\lambda$-a.e.~$z$ which then implies $\nu\ll\eta$.  

Let $$\text{$\rho = \frac{d\nu}{d\eta}$, \quad  $\rho_* = \frac{dS_*\nu}{dS_*\eta}$,\quad  and \quad $\rho_z = \frac{d\nu_z}{d\eta}$.}$$  
Let $\ptpart_Z$ denote the point partition on $Z$.  
We have $$\rho_*(S(x)) = \int \rho(y) \, d\eta_x ^{S\inv \ptpart_Z}(y)$$ $\eta$-a.s. 
Indeed,  given a measurable subset $C\subset Z$, the function $\1_C\circ S\colon X\to \R$ is $(S\inv \ptpart_Z)$-measurable and so
\begin{align*}
 \int \1_C(S(x))  & \left(\int \rho(y)  \, d \eta_x ^{S\inv \ptpart_Z} (y)\right) d   \eta (x) \\
 &= \int \1_C(S(y))  \rho(y) \,d\eta(y)= S_*\nu(C) =\int \1_C(z)   \rho_*(z)\, d   S_*\eta (z) .
\\&=\int \1_C(S(x))    \rho_*(S(x))\, d   \eta (x).
 \end{align*}
Moreover, in the  case $S$ is injective on common full measure set, we have $\rho_*(S(x)) =\rho(x)$.
 Similarly, we have  $\rho(x) = \int \rho_z(x) \ d\lambda(z)$.
Indeed for measurable $C\subset X$, we expand $\nu(C) $ as $$ \int  \nu_z(C) \, d \lambda(z) =
\int  \int _C \rho_z(x) \, d\eta(x) \, d \lambda(z) =
\int_C  \int\rho_z(x) \, d\lambda(z)  \, d \eta(x) .$$

Convexity of $\phi\colon t\mapsto t\log t$  then yields 
\begin{align}
\KL(S_*\nu &\parallel S_*\eta)
= \int \log (\rho_*(z)) \, dS_*\nu(z)
=  \int \rho_*(z) \log (\rho_*(z)) \, dS_*\eta(z)
\nonumber \\ &
=  \int \phi (\rho_*(z)) \, dS_*\eta(z)
=  \int \phi (\rho_*(S(x))) \, d\eta(x)
\nonumber \\ &
=  \int \phi \left( \int \rho (y) \, d\eta_x ^{S\inv \ptpart_Z} (y)\right) \, d\eta(x)
\nonumber \\ &
\le  \int  \left ( \int \phi(\rho (y)) \, d\eta_x ^{S\inv \ptpart_Z}(y) \right) \, d\eta(x)
\label{convex:equal} \\ &
=  \int \phi(\rho (y)) \, d\eta(y)
= \KL(\nu \parallel \eta)\nonumber
\end{align} 
with equality in \eqref{convex:equal} if $\rho_*(S(x)) = \rho(x)$ $\eta$-a.s.  Similarly, 
\begin{align*}
\KL(\nu  \parallel \eta )
	&=\int \phi  ( \rho ) \, d \eta 
	=\int \phi \left (\int \rho_z \, d\lambda(z)\right ) \, d \eta\\ &
\le \int \int \phi \left (\rho_z\right) \, d\lambda(z) \, d \eta\\ &
= \int \left( \int \phi \left (\rho_z \right) \, d \eta\right)  \ d\lambda(z) 
=\int \KL (\nu_z  \parallel \eta  )\, d \lambda(z). \qedhere
\end{align*}

\end{proof}

\subsection{Furstenberg entropy}
We return to the setting of a measure-pre\-serving transformation  $T\colon (X,\mu)\to (X,\mu)$ of a standard probability space $(X,\mu)$ and  a sub-invariant sub-$\sigma$-algebra $\sigF$.

Because of the tools we later employ in the proof of our main results, we will primarily be interested in the case that $h_\mu(T, \Xi(\sigF)) <\infty$ or, more generally, the case that 
the conditional measure $\mu^{T\inv \sigF}_x$ is supported on countably many $\sigF$-atoms for a.e.\ $x$.  
This,  in particular,  implies for $\mu$-a.e.\ $x$ that  \begin{equation}\mu_x^\sigF= \frac 1  {\mu_{x}^{T\inv  \sigF}(\sigF(x))} \restrict 
{{\mu_{x}^{ T\inv  \sigF}}}{\sigF(x)}.\end{equation}
We however present a unified definition that works even when  $\restrict{\sigF}{(T\inv \sigF)(x)}$ is uncountable.

\subsubsection{Furstenberg entropy}
We define the  \emph{Furstenberg entropy} of $T$ relative  to  $\sigF$ with respect  $\mu$ to be
$$\furstent_\mu(T\mid \sigF)=\int \KL(T_*\mu_x^\sigF\parallel \mu_{T(x)}^\sigF) \, d \mu(x).$$
By \eqref{eq:KLpositive}, we have $0\le \furstent_\mu(T\mid \sigF)\le \infty$.  

Expanding the Kullback--Leibler divergence, we have
\begin{align*}\furstent_\mu(T\mid \sigF)&
= -\int \int
\log \biggl(\frac{d\mu_{T(x)}^\sigF}{dT_*\mu_x^\sigF}(z)\biggr)  
 \,d(T_*\mu_x^\sigF)(z) \,d \mu(x)\\ 
 & = -\int \int
\log \biggl(\frac{d\mu_{T(x)}^\sigF}{dT_*\mu_x^\sigF}(T(y))\biggr)  
 \,d\mu_x^\sigF(y) \,d \mu(x)\\
  & = -\int 
\log \biggl(\frac{d\mu_{T(y)}^\sigF}{dT_*\mu_y^\sigF}(T(y))\biggr)  \,d \mu(y).
\end{align*}
(In the final expression above and in related expressions below, it is understood that we choose a jointly measurable  family of representatives   $ \rho_z(\cdot)= \frac{d\mu_{T(z)}^\sigF}{dT_*\mu_z ^\sigF}(\cdot)$  with $z\mapsto \rho_z$ constant on atoms of $\sigF$.)

\subsubsection{Related definitions}
We will typically be interested in settings where  $\furstent_\mu(T\mid \sigF)<\infty$. 
Note that the finiteness of Furstenberg entropy implies $T_*\mu_x^\sigF\ll \mu_{T(x)}^\sigF$ for $\mu$-a.e.\ $x$.

To unify our definition with others (cf.\ \cite[p.~1]{MR4561150}) in the literature, we have the following  
  identities that hold whenever $\furstent_\mu(T\mid \sigF)<\infty$.   
\begin{claim} Suppose  $T_*\mu_x^\sigF\ll \mu_{T(x)}^\sigF$ for $\mu$-a.e.\ $x$.  Then
\begin{align*}\furstent_\mu(T\mid \sigF)
 &=\int  \int 
\log \biggl(\frac{dT_*\mu_x^\sigF}{d\mu_{T(x)}^\sigF}(y)\biggr) \, d T_*\mu_x^\sigF (y)\,  d \mu(x)\\
 &=\int  
\log \biggl(\frac{dT_*\mu_x^\sigF}{d\mu_{T(x)}^\sigF}(T(x))\biggr)  
 \,  d \mu(x)
\\ &=
\int \int 
\log \biggl(\frac{dT_*\mu_x^\sigF}{d\mu_{T(x)}^\sigF} (y) \biggr)  \frac{dT_*\mu_x^\sigF}{d\mu_{T(x)}^\sigF} (y)
 \, d\mu_{T(x)}^\sigF(y) \, d \mu(x).
 \end{align*}
\end{claim}

 \subsubsection{Relation to i.i.d.\ random walks}
Recall all notation  in \cref{computefiberent}; in particular,  $M$ is compact metric space, 
$\nu_0$ is a probability on $C^{0}(M)$,  $\mu_0$ is $\nu_0$-stationary,  $\what F\colon \what \Sigma\times M\to \what \Sigma \times M$ is the skew extension preserving $\hat \mu = \nu_0^\N\times \mu_0$, and $\what \sigF$ is the sub-$\sigma$-algebra of sets saturated by preimages of the map $\what \Sigma\times M\to \what \Sigma$.     Then, the following relationship between the Furstenberg entropy of the stationary measure and the skew extension follows directly from definition:
$$\furstent_{\hat \mu}(\what F\mid \what \sigF)=\furstent_{\mu_0}({\nu_0}).$$

 \subsection{Properties of Furstenberg entropy}
We present various properties of  Furstenberg entropy for a measure-preserving transformation $T\colon (X,\mu)\to (X,\mu)$ and a sub-invariant sub-$\sigma$-algebra $\sigF$.

\subsubsection{Elementary bound} We have a naive upper bound on Furstenberg entropy.   
\begin{claim}\label{claim:furstupper}
We have  
\begin{equation}\furstent_\mu(T\mid \sigF) \le h_\mu(T, \Xi(\sigF))\label{furstbound}\end{equation}
\end{claim}
\begin{proof}
If $ h_\mu(T, \Xi(\sigF))=\infty$, there is nothing to prove.  
If $ h_\mu(T, \Xi(\sigF))<\infty$ then for a.e.~$x$, $\mu^{T\inv \sigF}_x$ is supported on countably many $\sigF$-atoms.  
Let $F_0, F_1, F_2, \dots$ denote the distinct (at most countably many) $\sigF$-atoms with positive $\mu^{T\inv \sigF}_x$-measure  in $(T\inv \sigF)(x)$ indexed so that  $\sigF(x) = F_0$.  
Then for $\mu$-a.e.~$x$, we expand the Lebesgue-Radon--Nikodym decomposition  $\mu_{T(x)} ^\sigF$: 
\begin{align*}\rho_x T_*\mu_x^\sigF + \eta_*&=\mu_{T(x)} ^\sigF= T_*\mu^{T\inv \sigF}_x\\&= 
 \sum_{j\ge 0} T_*  \restrict {\mu^{T\inv \sigF}_x}{ F_j} \\&= {\mu_{x}^{T\inv  \sigF}(\sigF(x))}  T_*\mu_x^\sigF + \sum_{j\ge 1} T_* \restrict {\mu^{T\inv \sigF}_x}{ F_j}
\end{align*}
Because $\eta_*\perp T_*\mu_x^\sigF$, it follows that $\rho_x(z) = \frac{d\mu_{T(x)}^\sigF}{dT_*\mu_x^\sigF}(z) \ge  {\mu_{x}^{T\inv  \sigF}(\sigF(x))} $ for $T_*\mu_x^\sigF $-a.e.\ $z$. 
Thus
\begin{align}\furstent_\mu(T\mid \sigF) & = -\int \int \log \rho_x(T(y))\, d\mu_x^\sigF(y) \, d\mu(x)\notag\\& \le -\int  \log {\mu_{x}^{T\inv  \sigF}(\sigF(x))}  \, d\mu(x)\label{jensenrigid}\\&= H_\mu(\sigF\mid T\inv \sigF)= h_\mu(T,\Xi(\sigF)). \notag \qedhere \end{align}
\end{proof}

\subsubsection{Countable-to-one factors} 
We have the following precise defect in \eqref{furstbound} that holds whenever $h_\mu(T, \Xi(\sigF))<\infty$.
\begin{claim}\label{claim:altfursten} Suppose $\mu^{T\inv \sigF}_x$ is supported on countably many $\sigF$-atoms for almost every $x$.   
Then for $\mu$-a.e.\ $x$,  {and for $\mu_x^\sigF$-a.e.\ $y$,}
$$\frac{dT_*\mu_x^\sigF}{d\mu_{T(x)}^\sigF}(T(y)) = \frac{\mu_y^{T\inv\sigB}(\sigF(x))}{\mu_x^{T\inv\sigF}(\sigF(x))}.$$

In particular, if $h_\mu(T, \Xi(\sigF))<\infty$ then 
\begin{align}\furstent_\mu(T\mid \sigF)
&= \int \log \frac{\mu^{T\inv\sigB}_x(\sigF(x))} {\mu_x^{T\inv\sigF}(\sigF(x))} \, d\mu(x) \notag \\&= H_\mu (\sigF\mid T\inv \sigF) - H_\mu(\sigF\mid T\inv \sigB) \notag \\ 
&= h_\mu (T,\Xi(\sigF)) - H_\mu(\sigF\mid T\inv \sigB). \label{furst3} 
\end{align}

\end{claim}

\begin{proof}
Let $F_0, F_1, F_2, \dots$ denote the distinct (at most countable) $\sigF$-atoms of  $(T\inv \sigF)(x)$ with positive $\mu_x^{T\inv \sigF}$-measure 
indexed so that  $\sigF(x) = F_0$.  
Given a measurable $C\subset \sigF(T(x))$, we have $$\mu_{T(x)}^\sigF(C) ={\mu_x^{T\inv \sigF}}(T\inv C) =  \sum_{j\ge 0} \mu^{T\inv \sigF}_x(T\inv C\cap F_j)$$
and (since $\mu_x ^{T\inv \sigF}(F_0)>0$),
$$T_*\mu_x^\sigF(C) = \mu_x^\sigF (T\inv C) =\dfrac 1 {\mu^{T\inv \sigF}_x(F_0)}\mu^{T\inv \sigF}_x(T\inv C\cap F_0).$$

The function $z\mapsto  \mu^{T\inv \sigB}_{z}(F_0)$ on $(T\inv \sigF)(x)$ is $T\inv \sigB$-measurable and thus there is a ($\sigB$-)measurable function $q\colon \sigF(T(x)) \to [0,1]$ such that $\mu^{T\inv \sigB}_{z}(F_0)=q( T(z))$ for $\mu_x^{T\inv \sigF}$-a.e.~$z$. 
As $C$ is $\sigB $-measurable, $T\inv C$ is $T\inv \sigB$-measurable and thus for any measurable $C\subset \sigF(T(x))$,
\begin{align*}
 {\mu^{T\inv \sigF}_x(F_0)} T_*\mu_x^\sigF(C) &=\mu^{T\inv \sigF}_x(T\inv C\cap F_0)
 \\
 &=
 \int _{T\inv C}
 \mu^{T\inv \sigB}_z(F_0)
 \, d \mu _x ^{T\inv \sigF}( z)\\
 &=
 \int _{T\inv C}
q(T(z))
 \, d \mu _x ^{T\inv \sigF}( z)
\\ &=
 \int _{{  C}}
{ q(w)}
 \, d \mu _{T(x)} ^{ \sigF}({ w})
\end{align*}
using that $T_* \mu _x ^{T\inv \sigF}=  \mu _{T(x)} ^{ \sigF}.$
{ Thus,
for $\mu_{T(x)}^\sigF$-a.e.~$w\in \sigF(T(x))$, $$ \frac{dT_*\mu_x^\sigF}{d\mu_{T(x)}^\sigF}({ w}) = \frac{q(w)}{\mu_x^{T\inv\sigF}(F_0)}.$$
Since $T_* \mu _x ^{ \sigF}\ll \mu _{T(x)} ^{ \sigF}$, for $\mu_x^\sigF$-a.e.\ $z$, we have 
\[\frac{dT_*\mu_x^\sigF}{d\mu_{T(x)}^\sigF}(T(z))=\frac{\mu_{ z}^{T\inv\sigB}(\sigF(x))}{\mu_x^{T\inv\sigF}(\sigF(x))}. \qedhere\]
}
\end{proof}

\subsubsection{(Non-)naturality under iteration}
As observed for transverse conditional entropy,  Furstenberg entropy in full generality is only sublinear under iteration.  
See \cref{sec:examples} for examples where Furstenberg entropy is strictly sublinear under iteration.
However, as in \cref{naturalityoftransverse},  Furstenberg entropy is  natural under iteration  under a mild injectivity hypothesis.
\begin{proposition}\label{subaddfurstent}
For any sub-invariant sub-$\sigma$-algebra $\sigF$ and $n\ge 0$,
\begin{equation}
\label{furstentsubadd}
\furstent_\mu(T^n\mid \sigF) \le n \furstent_\mu(T\mid \sigF).\end{equation}

Moreover, if the restriction $\restrict T{\sigF(x)}\colon \sigF(x)\to \sigF(T(x))$ of $T$ to $\mu$-a.e.\ atom $\sigF(x)$ of $\sigF$ is injective on a set of full $\mu_x^\sigF$-measure, then $$\furstent_\mu(T^n\mid \sigF) = n \furstent_\mu(T\mid \sigF).$$
\end{proposition}

\cref{furstentsubadd} of \cref{subaddfurstent} suggests that for a general transformation  $T$ and a general sub-invariant $\sigma$-algebra $\sigF$,  it might be more natural to study the quantity $$\inf _{n\ge 1} \frac 1 n \furstent_\mu(T^n\mid \sigF)$$
or to incorporate  into the definition of   $\furstent_\mu(T\mid \sigF)$ an appropriate notion of ``folding entropy'' when the injectivity hypothesis fails.
\begin{proof}[Proof of \cref{subaddfurstent}]
Write  $$h_n = \furstent_\mu(T^n\mid \sigF) = \int \KL(T^n_* \mu^\sigF_x\parallel  \mu^\sigF_{T^n(x)}) \, d \mu(x).$$

We first claim $n\mapsto h_n$ is monotonic nondecreasing.
Let $\sigG= T^{-1}\sigF$. 
Note that  $$\text{$T_* \mu^{\sigG}_x= \mu^\sigF_{T(x)}$,\quad  $\int \mu^\sigF_y \, d\mu_x ^{\sigG}(y) = \mu^{\sigG}_x$, \quad and so \quad  
$\int T_*\mu^\sigF_y \, d\mu_x ^{\sigG}(y) = \mu^\sigF_{ T(x)}$}.
$$
 Applying
\cref{lem:superbconvexity}\ref{convex2} and invariance of $\mu$,  for $n\ge 2$ we have 
\begin{align*}
h_{n-1}&=
\int  \KL\left( T^{n-1}_* \mu^\sigF_x  \parallel  \mu^\sigF_{T^{n-1}(x)}\right) \, d \mu(x)
\\&
=\int  \KL\left( T^{n-1}_* \mu^\sigF_x  \parallel  \mu^\sigF_{T^{n-1}(x)}\right) \, d T_*\mu(x)
\\&
=\int  \KL\left( T^{n-1}_* \mu^\sigF_{T(x)}  \parallel  \mu^\sigF_{T^n(x)}\right) \, d \mu(x)
\\&
= \int  \KL\left( \int T_*^{n-1}T_* \mu^\sigF_{y} \, d\mu_x ^{T\inv \sigF}(y)  \parallel  \mu^\sigF_{T^n(x)}\right) \, d \mu(x)
\\&
\le \int  \int\KL\left(  T_*^{n-1}T_* \mu^\sigF_{y}  \parallel  \mu^\sigF_{T^n(y)}\right) \, d\mu_x ^{T\inv \sigF}(y) \, d \mu(x)
\\&
=
 \int  \KL\left( T^n_* \mu^\sigF_x  \parallel  \mu^\sigF_{T^n(x)}\right) \, d \mu(x)=h_n.
\end{align*}

The proposition follows trivially if $h_1=\infty$.  Thus we suppose $h_1<\infty$.  To show the proposition, it suffices to show for every $n\in \{2,3,4,\dots\}$ that $h_n<\infty$ and $h_n \le h_{n-1}+ h_1$ with equality if the restriction $\restrict T{\sigF(x)}\colon \sigF(x)\to \sigF(T(x))$ of $T$ to $\mu$-a.e.\ atom $\sigF(x)$ of $\sigF$ is injective a.s.
Suppose we know $h_{n-1}<\infty$.  Then, since $h_1<\infty$, we have joint absolute continuity
$${T^{n}_*\mu_x^\sigF}\ll {T_*\mu_{T^{n-1}(x)}^\sigF}\ll {\mu_{T^{n}(x)}^\sigF}$$
for a.e.~$x$.
In particular, 
$$\frac{dT^{n}_*\mu_x^\sigF} {d\mu_{T^{n}(x)}^\sigF}= 
\frac {dT^{n}_*\mu_x^\sigF} {dT_*\mu_{T^{n-1}(x)}^\sigF} \cdot \frac {dT_*\mu_{T^{n-1}(x)}^\sigF} {d\mu_{T^{n}(x)}^\sigF}$$
and so 
\begin{align*}
h_n&=
 \int  \KL\left( T^n_* \mu^\sigF_x  \parallel  \mu^\sigF_{T^n(x)}\right) \, d \mu(x)\\
&=\int \left(\int \log \left(\frac {dT^{n}_*\mu_x^\sigF} {dT_*\mu_{T^{n-1}(x)}^\sigF}(y) \cdot \frac {dT_*\mu_{T^{n-1}(x)}^\sigF} {d\mu_{T^{n}(x)}^\sigF} (y)\right)  \, d T^n_*\mu_x^\sigF(y)\right) \, d\mu(x)\\
&=
\int \KL\left(
{T^{n}_*\mu_x^\sigF} \parallel {T_*\mu_{T^{n-1}(x)}^\sigF}
\right)  d \mu(x)\\ &\quad \quad \quad+
\int \left(\int \log \left(  \frac {dT_*\mu_{T^{n-1}(x)}^\sigF} {d\mu_{T^{n}(x)}^\sigF} (y)\right)  \, d T^n_*\mu_x^\sigF (y)\right) \, d\mu(x).
\end{align*}
(Additivity in the final line is justified because the expansions and bounds below show the integrands are $L^1$.)

Applying \cref{lem:superbconvexity}\ref{convex1} to the first integral, we  have 
\begin{align}
\int \KL&\left(
{T^{n}_*\mu_x^\sigF} \parallel {T_*\mu_{T^{n-1}(x)}^\sigF}
\right)  d \mu(x) \nonumber 
=
\int \KL\left(
{T_*T^{n-1}_*\mu_x^\sigF} \parallel {T_*\mu_{T^{n-1}(x)}^\sigF}
\right)  d \mu(x)
\nonumber \\&\le
\int \KL\left(
{T^{n-1}_*\mu_x^\sigF} \parallel {\mu_{T^{n-1}(x)}^\sigF}
\right)  d \mu(x)\label{furstsubad}
 =
h_{n-1}. 
\end{align}
Moreover, since ${T^{n-1}_*\mu_x^\sigF} \ll {\mu_{T^{n-1}(x)}^\sigF}$ (by the hypothesis $h_{n-1}<\infty$), assuming the restriction  $\restrict T{\sigF(x)}\colon \sigF(x)\to \sigF(T(x))$ of $T$
is injective on a set of full  $\mu_x^\sigF$-measure for $\mu$-a.e.~$x$, 
the inequality in \eqref{furstsubad} becomes an equality.  
We rearrange the second integral as follows:
\begin{align}
\int &\left(\int \log \left(  \frac {dT_*\mu_{T^{n-1}(x)}^\sigF} {d\mu_{T^{n}(x)}^\sigF} (y)\right)  \, d T^n_*\mu_x^\sigF (y)\right) \, d\mu(x)\nonumber\\
&=
\int \int  \left(\int\log \left(  \frac {dT_*\mu_{T^{n-1}(x)}^\sigF} {d\mu_{T^{n}(x)}^\sigF} (y)\right)  \, d T^n_*\mu_z^\sigF (y)\right)  d \mu_{x}^{T^{-(n-1)}\sigF}(z)\, d\mu(x)\nonumber\\
&=
\int \left(\int   \log \left(  \frac {dT_*\mu_{T^{n-1}(x)}^\sigF} {d\mu_{T^{n}(x)}^\sigF} (y)\right)  \, d T_* (\mu_{T^{n-1}(x)}^{\sigF})(y)\right)\, d\mu(x) \nonumber\\
&=\int \KL\left({T_*\mu_{T^{n-1}(x)}^\sigF} \parallel {\mu_{T^{n}(x)}^\sigF} \right)\, d\mu(x)\nonumber\\
&=
\int \KL\left({T_*\mu_{x}^\sigF} \parallel {\mu_{T(x)}^\sigF} \right)\, dT^{n-1}_*\mu(x)=h_1.\nonumber
\end{align}
The proposition then follows.
\end{proof}

\subsubsection{Extreme values of Furstenberg entropy}
Like most other notions of entropy,   Furstenberg entropy has natural geometric interpretations that become especially transparent when it attains its extreme values:  $\furstent_\mu(T\mid \sigF) $ measures the average amount of (de)correlation between $T_*\mu_y^\sigF$ as $y$ varies across $\sigF$-atoms of $T\inv \sigF$.  This is made precise in the next two claims.

Recall from \cref{claim:furstupper} that $\furstent_\mu(T\mid \sigF) \le h_\mu(T, \Xi(\sigF))$.  When $h_\mu(T, \Xi(\sigF))<\infty,$ attaining this upper bound implies complete decorrelation among the measures $T_*\mu_y^\sigF$ as $y$ varies across $\sigF$-atoms of $T\inv \sigF$.

\begin{claim}
Suppose $h_\mu(T, \Xi(\sigF))<\infty$ and $\furstent_\mu(T\mid \sigF)= h_\mu(T, \Xi(\sigF))$.  Then for $\mu$-a.e.\ $x$ and $\mu_x^{T\inv \sigF}$-a.e.\ $z\in (T\inv \sigF)(x)\sm \sigF(x)$, 
$$T_*\mu^\sigF_x \perp T_*\mu^\sigF_z.$$
\end{claim}
\begin{proof}Indeed,  in the proof of \cref{claim:furstupper}, by monotonicity of the integral, equality in \cref{jensenrigid}
implies 
$$\frac{dT_*\mu_x^\sigF}{d\mu_{T(x)}^\sigF} (y) =\frac{1}{ \mu_x^{T\inv \sigF}(\sigF(x))}$$ for 
$T_*\mu_x^\sigF$-a.e.\ $y$ and $\mu$-a.e.\ $x$.   
Let $F_0, F_1, F_2, \dots$ be the distinct (at most countably many) $\sigF$-atoms in $(T\inv \sigF)(x)$ with positive $\mu_x^{T\inv \sigF}$-measure 
indexed with $\sigF(x) = F_0$.  
Then uniqueness of the Lebesgue-Radon--Nikodym decomposition implies \begin{align*}\mu_{T(x)} ^\sigF 
&= T_*\mu^{T\inv \sigF}_x  
= {\mu_{x}^{T\inv  \sigF}(\sigF(x))}  T_*\mu_x^\sigF + \sum_{j\ge 1} T_* \restrict {\mu^{T\inv \sigF}_x}{ F_j}\\
&=  {\mu_{x}^{T\inv  \sigF}(\sigF(x))} T_*\mu_x^\sigF + \eta
\end{align*}
where $\eta\perp T_*\mu_x^\sigF$. 
\end{proof}

At the other extreme, when $\furstent_\mu(T\mid \sigF)=0$, the measures $T_*\mu_y^\sigF$ are completely correlated as $y$ varies across $\sigF$-atoms of $T\inv \sigF$.

\begin{claim}\label{claim:zero-furstenberg-entropy}
If $\furstent_\mu(T\mid \sigF)=0$ then $T_*\mu^\sigF_x= \mu^\sigF_{T(x)}$ for $\mu$-a.e.\ $x$.  In particular, $\furstent_\mu(T\mid \sigF)=0$ if and only if for $\mu$-a.e.\ $x$ and $\mu^{T\inv \sigF}_x$-a.e.\ $y$, 
$$T_*\mu_x^\sigF= T_*\mu_y^\sigF.$$
\end{claim}

\begin{proof}
This follows directly from definition and \cref{eq:KLzero}.
\end{proof}

This immediately implies the following.  
\begin{corollary}
 $\sigF$ admits an {equivariant family of conditional measures} if and only if $\furstent_\mu(T\mid \sigF)=0$.
\end{corollary}

\subsubsection{Relationship with other entropies}
Finally, we have the following which relates the Furstenberg entropy with the standard and transverse conditional entropies. 
Although the following almost certainly holds in more generality, we only formulate the proposition under the additional assumption that $h_\mu(T, \Xi(\sigF))<\infty$.  In this case, the following is a direct consequence of \eqref{furst3} of  \cref{claim:altfursten} and  \cref{prop:GenAR}.

\begin{proposition}\label{entropygap}
Let $\sigF$ be a sub-invariant sub-$\sigma$-algebra with $h_\mu(T, \Xi(\sigF))<\infty$.  Then 
$$\furstent_\mu(T\mid \sigF)+ h_\mu(T\mid \sigF)=  \transent_\mu(T\mid \sigF).$$
\end{proposition}

\subsection{Examples}\label{sec:examples}
We conclude this section by presenting  examples  that showcase various definitions and properties of entropies for i.i.d.\ random walks with non-injective generators.  
\subsubsection{Random expanding maps}
For $r\ge 2$, let $E^r(S^1)$ denote the set of $C^r$ expanding maps $f\colon S^1\to S^1$.  Let $\nu_0$ be a finitely supported\footnote{We take $\nu_0$ finitely supported for simplicity.  One could take non-atomic measures $\nu$ on $E^r$ satisfying appropriate moment conditions.}  probability measure on $E^r$ and set  $\supp( \nu_0)= \{f_1, \dots, f_\ell\}$.  Let $\lambda:= \min\{|f'_i(x)|: 1\le i \le \ell, x\in S^1\}$.  By assumption $\lambda>1$. 

Let $m$ denote the Lebesgue measure on $S^1$ and let $d(\cdot, \cdot)$ denote the distance on $S^1$.    Given a function $\rho\colon S^1\to (0,\infty)$, let $$\|\rho\|_{\log-L}:= \sup_{x\neq y}\frac{|\log(\rho(x))-\log(\rho(y))|}{d(x,y)}$$ denote the $\log$-Lipschitz semi-norm.  
A standard invariant-cone argument
(see \cite[(3.1) and Lem.~3.1]{MR1402416})
shows 
that for all  $L>0$ sufficiently large, 
for any 
function $\rho\colon S^1\to (0,\infty)$ with $\|\rho\|_{\log-L}\le L$, and for any $1\le j\le \ell$, if $\mu:= \rho m$ denotes the measure with density $\rho$ and if $(f_j)_*\mu = \hat\rho m$, then $$\|\hat \rho\|_{\log-L}\le L.$$  
In particular, if $\int \rho \, dm =1$ (so that $\rho m$ is a probability) and if $\|\rho\|_{\log-L}<\infty$, then there exists $C>1$ such that for any $n\ge 1$ and any sequence of indices $i_0,\dots, i_{n-1}\in \{1,\dots,\ell\}$ of length $n$,
\begin{equation}\label{eq:bdd}\frac 1 C \le \frac {d(f_{i_{n-1}}\circ\cdots\circ f_{i_0})_*(\rho m)}{d\rho m} \le C.\end{equation}

We may ensure $\nu_0$ was chosen so that there is no absolutely continuous probability measure invariant under every $f_j \in \supp (\nu_0)$.  On the other hand, since $\nu_0$ is finitely supported, there exists a unique $\nu_0$-stationary  probability measure  $\mu_0$ equivalent to the Lebesgue measure on $S^1$ with a $\log$-Lipschitz density.  (See e.g.\ \cite{MR1359097,MR2399927}.)

We follow the notation of the associated skew-extension and fiber $\sigma$-algebra in \cref{skewext}.  
Since $\mu_0$ is not a common invariant measure for $\nu_0$-a.e.\ map $f$, we have
$$\furstent_{\mu_0}({\nu_0})= \furstent_{\hat \mu}(\what F\mid \what \sigF)>0.$$
On the other hand, \cref{eq:bdd} ensures there is $\hat C$ such  that for all $n\ge 1$,
$$\furstent_{\mu_0}({\nu_0^{\otimes n}})= \furstent_{\hat \mu}(\what F^n \mid \what \sigF)\le \hat C.$$
In particular, $$\lim_{n\to \infty} \frac 1 n \furstent_{\hat \mu}(\what F^n\mid \what \sigF) =0.$$

On the other hand, we have (cf.~\cite{MR1944407})
$$\log(\lambda)\le h_{\mu_0}({\nu_0})= h_{\hat \mu}(\what F \mid \what \sigF).$$
In particular, the standard entropy 
$ h_{\hat \mu}(\what F^n\mid \what \sigF)=nh_{\hat \mu}(\what F \mid \what \sigF)\ge n \log(\lambda)$ grows linearly in $n$.  Moreover,  by \cref{entropygap}, the transverse entropy also grows linearly in $n$ since 
$$ \lim_{n\to \infty } \frac 1 n \transent_{\hat \mu}(\what F^n\mid \what\sigF)= h_{\hat \mu}(\what F\mid \what \sigF)\ge \log(\lambda).$$ 

\subsubsection{Non-injective affine contractions}
Let $\mathcal C\subset [0,1]$ denote the standard middle third Cantor set.  Let $\mu_0$ denote the Cantor-Lebesgue 
measure on $\mathcal C$. (That is, the measure that gives mass $ {2^{-n}}$ to the level-$n$ basic Cantor intervals; i.e.~intersections of the form $[j3^{-n},(j+1)3^{-n}]\cap \mathcal C$ with infinite cardinality.)

Let $I$ be an interval properly containing $[0,1]$.  E.g.~let $I=[-1,2]$.  Fix $r\ge 0$.    For $1\le j\le 4$, let $T_j\colon I\to I$ be a $C^r$ map such that:
\begin{enumerate}
\item $T_1$ maps $[0,\frac 1 3]$ onto $[0,\frac 1 9]$ by an increasing  affine contraction and maps
$[\frac 2 3,1]$ onto $[0,\frac 1 3]$ by an increasing  affine isometry;
\item $T_2$ maps $[0,\frac 1 3]$ onto $[\frac 2 3,1]$ by  an increasing affine isometry and maps
$[\frac 2 3,1]$ onto $[\frac 2 9,\frac 1 3]$ by  an increasing affine contraction;
\item $T_3$ maps $[0,\frac 1 3]$ onto $[\frac 2 3,1]$ by  an increasing affine isometry and maps
$[\frac 2 3,1]$ onto $[\frac 8 9, 1 ]$ by  an increasing affine contraction;
\item $T_4$ maps $[0,\frac 1 3]$ onto $[\frac 2 3,\frac 7 9]$ by  an increasing affine contraction and maps
$[\frac 2 3,1]$ onto $[0,\frac 1 3]$ by  an increasing affine isometry.
\end{enumerate}
Outside $[0,\frac 1 3]\cup [\frac 2 3,1]$, extend $T_j$ by a $C^r$-map, for $r\ge 1$.  
Let $\nu_0:= \frac 1 4 (\delta_{T_1}+\delta_{T_2}+\delta_{T_3}+\delta_{T_4}).$ 
Then $\mu_0$ is $\nu_0$-stationary (and ergodic) though $\mu_0$ is not invariant under any $T_i$.  

As above, 
we follow the notation of the associated skew-extension and fiber $\sigma$-algebra in \cref{skewext}.  
Since the  dynamics of each $T_i$ restricted to the fundamental intervals on $\mathcal C$ (the support of $\mu_0$) is either a contraction or an isometry, the standard entropy of the random system, which coincides with the standard conditional entropy, vanishes: $$h_{\mu_0}({\nu_0})= h_{\hat \mu}(\what F\mid \what \sigF)=0.$$

To calculate the Furstenberg entropy, we first compute:
\begin{align*}
\frac {d\mu_0}{d(T_1)_*\mu_0}(x)&= 
\begin{cases}
\frac{1}{3} &x\in [0, \frac 1 9]\cap \mathcal C\\
1&x\in [\frac 2 9,  \frac 1 3]\cap \mathcal C\\
\infty &x\in [ \frac 2 3, 1]\cap \mathcal C
\end{cases}
&\frac {d\mu_0}{d(T_2)_*\mu_0}(x)&= 
\begin{cases}
\infty &x\in [0, \frac 1 9]\cap \mathcal C\\
\frac{1}{2} &x\in [ \frac 2 9,\frac 1 3]\cap \mathcal C\\
1&x\in [\frac 2 3,   1 ]\cap \mathcal C\\
\end{cases}\\
\frac {d\mu_0}{d(T_3)_*\mu_0}(x)&= 
\begin{cases}
\infty &x\in [ 0,\frac 1 3]\cap \mathcal C\\
1&x\in [\frac 2 3,  \frac 7 9]\cap \mathcal C\\
\frac{1}{3} &x\in [ \frac 8 9,1]\cap \mathcal C\\
\end{cases}
&\frac {d\mu_0}{d(T_4)_*\mu_0}(x)&= 
\begin{cases}
1 &x\in [0,\frac 1 3]\cap \mathcal C\\
\frac 1 2 &x\in [\frac 2 3,  \frac 7 9 ]\cap \mathcal C\\
\infty &x\in [\frac 8 9 , 1]\cap \mathcal C
\end{cases}
\end{align*}
and so 
\begin{align*}
\furstent_{\mu_0}({\nu_0})&= \furstent_{\hat \mu}(\what F\mid \what\sigF)\\
	&= -\tfrac 1 4 \Bigl(
	\tfrac 3 4 \log \left(\tfrac 1 3\right)+\tfrac 1 2 \log \left(\tfrac 1 2\right)+\tfrac 3 4 \log \left(\tfrac 1 3\right)+\tfrac 1 2 \log \left(\tfrac 1 2\right)
	\Bigr)\\
	&= \tfrac 3 8 \log 3 + \tfrac 1 4 \log 2.
\end{align*}
We omit calculations,  but may similarly compute 
\begin{align*}\furstent_{\mu_0}({\nu_0^{\otimes 2}})&= \furstent_{\hat \mu}(\what F^2\mid \what \sigF)\\&=
\tfrac{3}{4}  \log (2)+\tfrac{3}{8}  \log (3) +\tfrac{5}{64}  \log (5)
\\ &<2 \furstent_{\mu_0}({\nu_0}).
\end{align*}
Again, Furstenberg entropy is strictly sub-linear under iteration.

Since standard entropy vanishes, by \cref{entropygap}, the transverse entropy coincides with the Furstenberg entropy.  However, we may directly compute  transverse entropy as follows: Letting $\xi$ denote the partition of $\mathcal C$ into 4 intervals of the form $[j/9, (j+1)/9]\cap \mathcal C$, we claim that $\transent_{\mu_0}({\nu_0},\xi)=\transent_{\mu_0}({\nu_0})$. The proof follows from the following computation:  First, since each $T_i$ maps atoms of $\xi$ into atoms of $\xi$, we have $\xi^{-} _\omega= \xi$ for every $\omega$.
We then calculate 
\begin{align*}
H_{\mu_0}(\xi\mid T_1\inv \xi)&=H_{\mu_0}(\xi\mid T_3\inv \xi)=-\left (\tfrac 3 4 \log \left(\tfrac 1 3\right )  + \tfrac 1 4 \log (1)\right)\\
H_{\mu_0}(\xi\mid T_2\inv \xi)&=H_{\mu_0}(\xi\mid T_4\inv \xi)=-\left (\tfrac 1 2 \log \left(\tfrac 1 2\right )  + \tfrac 1 2 \log (1)\right).
\end{align*}
It follows that 
$$\transent_{\mu_0}({\nu_0},\xi)= \transent_{\hat \mu}(\what F,\xi\mid \what \sigF)= \tfrac 3 8 \log 3 + \tfrac 1 4 \log 2.$$
Since $h_{\mu_0}({\nu_0})=0$, by \cref{entropygap} we have 
$$\transent_{\hat \mu}(\what F,\xi \mid \what \sigF)=\transent_{\mu_0}({\nu_0},\xi)= \transent_{\mu_0}({\nu_0})=\furstent_{\mu_0}({\nu_0})$$ and so $$\transent_{\hat \mu}(\what F \mid \what \sigF)=\transent_{\hat \mu}(\what F,\xi \mid \what \sigF)=\furstent_{\mu_0}({\nu_0}).$$

Finally, we verify \cref{claim:altfursten} in this example by calculating
$H_{\hat \mu}(\what \sigF\mid \what F\inv \sigB).$
Every point $(\omega, x)$ in $\what \Sigma\times \mathcal C$ has exactly 3 preimages (in $\what \Sigma\times \mathcal C$) under the skew-extension $\what F$ with the following weights on the associated conditional measure:
\begin{itemize}
\item For $x\in[0,\frac 1 9]$, we have
${\hat \mu}^{\what F\inv\sigB}_{\what F\inv(\omega,x)}:= 
\tfrac 1 2 \delta_{(T_1\omega, 3x)} + \tfrac 1 4 \delta_{ (T_1\omega, x+ \tfrac 2 3)} + \tfrac 1 4 \delta_{ (T_4\omega, x+ \tfrac 2 3)}.$
\item For $x\in[\frac 2 9,\frac 1 3]$,
we have ${\hat \mu}^{\what F\inv\sigB}_{\what F\inv (\omega,x)}:= 
\tfrac 1 4 \delta_{(T_1\omega, x+\tfrac 2 3)} + \tfrac 1 2 \delta_{ (T_2\omega, 3 x)} + 
\tfrac 1 4 \delta_{ (T_4\omega, x+ \tfrac 2 3)}.$

\item For $x\in[\frac 2 3,\frac 7 9]$, we have 
${\hat \mu}^{\what F\inv\sigB}_{\what F\inv(\omega,x)}:= 
\tfrac 1 4\delta_{(T_2\omega, x-\tfrac 2 3)}+\tfrac 1 4\delta_{ (T_3\omega, x-\tfrac 2 3)}+\tfrac 1 2\delta_{ (T_4\omega, 3x-2)}.$

\item For $x\in[\frac 8 9, 1 ]$, we have 
${\hat \mu}^{\what F\inv\sigB}_{\what F\inv (\omega,x)}:= 
\tfrac 1 4\delta_{(T_2\omega, x-\tfrac 2 3)}+\tfrac 1 4\delta_{ (T_3\omega,  x-\frac 2 3)}+\tfrac 1 2\delta_{ (T_3\omega, 3x - 2 )}.$
\end{itemize}
Moreover, for each non-trivial (i.e.~with infinite cardinality) interval of the form $J=[j/9,(j+1)/9]\cap \mathcal C$ and $\omega\in \what \Sigma$, the preimage   $\what F\inv (\{\omega\}\times J)$ in $\what F\inv (\{\omega\}\times \mathcal C)$ has relative measure $\frac 1 4$.  
It follows from definition that
\begin{align*}H_{\hat \mu} (\what \sigF\mid \what F\inv\sigB)
&=- \tfrac 1 4 \Bigl(
(\tfrac 3 4\log (\tfrac 3 4) + \tfrac 1 4\log (\tfrac 1 4))+
(\tfrac 1 2\log (\tfrac 1 2) + \tfrac 1 2\log (\tfrac 1 4))\\&\quad \quad \quad+
(\tfrac 3 4\log (\tfrac 3 4) + \tfrac 1 4\log (\tfrac 1 4))+(\tfrac 1 2\log (\tfrac 1 2) + \tfrac 1 2\log (\tfrac 1 4))
\Bigr)\\
&=-\tfrac 3 8 \log 3 +\tfrac 7 8 \log 4.
\end{align*}
Since $h_{\hat \mu}(\what F, \Xi(\what \sigF))=\log 4$, we verify \cref{furst3} of \cref{claim:altfursten} for this example:
\begin{align*}\furstent_{\hat \mu}(\what F\mid \what\sigF)&=\log 4 - \left(-\tfrac 3 8 \log 3 +\tfrac 7 8 \log 4\right)\\&= \tfrac 3 8 \log 3 +\tfrac 1 8 \log 4= \tfrac 3 8 \log 3 +\tfrac 1 4 \log 2.\end{align*}
Since the number of preimages under $\what F$ is $3$, we have $H_{\hat \mu} (\what \sigF\mid \what F^{-n}\sigB)\le n\log 3$ for all $n \ge 1$.  In particular, unlike the above example, 
$$\furstent_{\hat \mu}(\what F^n\mid \what\sigF)\ge n \left(\log 4-\log 3\right)$$
grows linearly in $n$.

\subsubsection{A two-point process} 
Let $M=\{x_1,x_2\}$.  Let $T_1,T_2\colon M\to M$ be given by $T_1(x_1)= T_1(x_2)= x_2$ and $T_2(x_1)= x_2$ and $T_2(x_2)= x_1$.
With $\nu_0= \frac 1 2 (\delta_{T_1}+ \delta_{T_2})$ and $\mu_0 = \frac 1 3 \delta_{x_1}+ \frac 2 3 \delta_{x_2}$, we have that $\mu_0$ is $\nu_0$-stationary.  
Direct computation  verifies strict sub-linearity and \cref{furst3} of \cref{claim:altfursten} in this example: 
\begin{enumerate}
\item $\furstent_{\mu_0}(\nu_0)=\frac{1}{6} \log \left(\frac{27}{4}\right)$.
\item $H_{\hat \mu} (\what \sigF\mid \what F^{-1}\sigB)=\frac{1}{6} \log \left(\frac{256}{27}\right)$.
\item $\furstent_{\mu_0}(\nu_0)=\log(2)-H_{\hat \mu} (\what \sigF\mid \what F^{-1}\sigB)$.
\item $\furstent_{\mu_0}(\nu_0^{\otimes 2})=\frac{1}{4} \log \left(\frac{27}{4}\right)<2\furstent_{\mu_0}(\nu_0)$.
\item $H_{\hat \mu} (\what \sigF\mid \what F^{-2}\sigB)= \frac{1}{4} \log \left(\frac{1024}{27}\right)$.
\item $\furstent_{\mu_0}(\nu_0^{\otimes 2})=\log(4)-H_{\hat \mu} (\what \sigF\mid \what F^{-2}\sigB)$.
\end{enumerate}

\section{Formulation of results for abstract skew extensions}\label{formulation and results}
\subsection{Model of random dynamics and fibered systems with holonomies} 
We present a more general framework for our results than what we previously presented in \cref{sec:iid}.  This, in particular, allows for random walks that are identically distributed but non-independent and for applications to fibered dynamics admitting holonomies.  

\subsubsection{Measurable base dynamics}
We consider a standard probability space $(\widehat \Omega,\what \sigB, \hat \nu)$.  Let $\hat \theta \colon \widehat \Omega\to \widehat \Omega$ be  an ergodic, measurable, $\hat \nu$-preserving transformation.  

Though not strictly necessary, for our applications of interest, $\hat \theta$ will be assumed non-invertible.  
Let $\tau\colon (\Omega,\sigB, \nu) \to (\what \Omega, \what \sigB, \hat \nu)$ be the ergodic-theoretic  natural extension equipped with the (unique) invertible, ergodic,  $\nu$-preserving transformation $\theta\colon \Omega\to \Omega$ satisfying $\hat \theta \circ \tau = \tau \circ \theta.$  
Let  $\sigA\subset \sigB$ be the sub-$\sigma$-algebra $\sigA= \tau\inv \what \sigB$.  Then $\sigA$ is sub-invariant for $\theta$; by uniqueness of the natural extension, we have $\sigB= \sigA^+$.

Let $M$ be a compact smooth manifold.  Fix $\beta>0$ and let $\widehat \Omega\to \diff^{1+\beta}(M)$,  $$  \omega\mapsto f_{  \omega},$$ be a ($\what \sigB$-)measurable\footnote{We equip $\Diff^{1+\beta}(M)$ with the standard Borel structure
inherited from $\Diff^1(M)$ with its $C^1$ topology.  We implicitly impose a Polish model of every standard probability space appearing.}  function on $\widehat \Omega$.  We also write $\Omega\to\diff^{1+\beta}(M)$,   $\omega\mapsto f_\omega$, for the $\sigA$-measurable function given by $$f_\omega := f_{\tau (\omega)}.$$

Given \footnote{Throughout, we will use $\omega$ to represent elements of both  $\Omega$ and $\widehat \Omega$.} $\omega\in \Omega$ and $n\in \Z$, we write $$f_{\omega}^{(n)}= \begin{cases}
f_{\theta^{n-1}(\omega)}\circ \cdots \circ f_{\omega} & n\ge 1\\
\id & n=0\\
\bigl(f_{\theta \inv \omega}\circ \cdots \circ f_{\theta^n(\omega)}\bigr)\inv &
n<0\\
\end{cases}$$ 
with similar conventions for $f^{(n)}_{  \omega}$ for $  \omega\in \widehat \Omega$ and $n\ge 0$.  

\subsubsection{Fibered dynamics}
Let $\what X= \widehat \Omega\times M$ and equip $\what X$ with the map $\what F\colon \what X\to \what X$, 
$$\what F(  \omega, x)= \bigl(\hat \theta(  \omega), f_{ \omega}(x)\bigr).$$
Similarly, let $X=    \Omega\times M$ be equipped with the map $  F\colon   X\to   X$, 
$  F(  \omega, x)= \bigl(  \theta(  \omega), f_{ \omega}(x)\bigr).$
By abuse of notation, extend $\tau\colon \Omega \to \widehat \Omega$ to $\tau\colon X\to  \what X$ by declaring $$\tau (\omega, x) = (\tau (\omega), x).$$
Then $\tau$ induces $\what F$ as a factor of $F$.  
Also write $\pi_1\colon X\to \Omega$, $\pi_2\colon X\to M$, $\hat \pi_1\colon \what X\to \widehat \Omega$, and  $\hat \pi_2\colon \what X\to M$  for the natural factor maps.  
We equip   $\what X$ (resp.\ $X$) with the $\sigma$-algebra induced by $\what \sigB$ (resp.\ $\sigB$) and the Borel $\sigma$-algebra $\sigB_M$ on $M$. 

We have the following standard facts which follow from the construction of (ergodic-theoretic) natural extensions.  
\begin{claim}\label{unique natext}
Let $\hat \mu$ be an $\what F$-invariant probability measure with $(\hat \pi_1)_*\hat \mu = \hat \nu$.  Then there exists a unique $F$-invariant probability measure $\mu$ on $X$ such that  $F\colon (X,\mu) \to (X,\mu)$ is the ergodic theoretic natural extension of $\what F\colon (\what  X,\hat\mu) \to (\what  X,\hat \mu)$.  In particular, 
\begin{enumerate}
\item $(\pi_1)_*\mu= \nu$, 
\item $\tau_*\mu = \hat \mu$.
\item $\mu$ is ergodic if and only if $\hat \mu$ is ergodic.
\end{enumerate}
\end{claim}

Conversely, let $\mu$ be any ergodic, $F$-invariant probability measure on $X$ projecting to $\nu$ on $\Omega$.  Taking $\hat \mu=  \tau_*\mu$, we have that $\hat \mu$ is ergodic, $\what F$-invariant, and $( \hat \pi_1)_* \hat \mu = \hat \nu$.  Thus $\mu$ is the unique $F$-invariant probability measure guaranteed by \cref{unique natext} for this $\hat \mu$.

\subsubsection{Further $\sigma$-algebras of interest}
We define the following additional sub-$\sigma$-algebras:
\begin{enumerate}
\item $\sigF= \pi_1\inv \sigB$ is  the sub-$\sigma$-algebra of sets saturated by fibers of $\pi_1$.  Each atom of $\sigF$ is of the form $\{\omega \}\times M$. 
 \item Similarly, $\widehat \sigF= \hat \pi_1\inv \what \sigB$.  
\item $\sigC= \pi_1\inv \sigA = \tau\inv \what \sigF$ is the sub-$\sigma$-algebra on $X$ of sets formed by  saturating sets in $\sigA$ by fibers of $\pi_1$.   Note that  $\sigF= \sigC^+$.  
\item Let $\wtd \sigA$ denote the product $\sigma$-algebra $\sigA\otimes \sigB_M$.  
\end{enumerate}

\subsubsection{Fiber and sample measures}
We simplify notation and write $\hat \mu_\omega= \hat \mu_\omega^{\what \sigF}$ and $ \mu_\omega= \mu_\omega^{\sigF}$.  We refer to the  collection $\{\hat \mu_\omega\}$ as the \emph{fiber measures} of $\hat \mu$ and refer to the  collection $\{\mu_\omega\}$ as the \emph{fiber measures} or \emph{sample measures} of $\mu$.  The dimension properties   of the sample measures $\mu_\omega$, and their  relation to standard conditional entropy, have been studied in \cite{MR2218998,MR2032494}.  We instead study dimension properties of the fiber measures $\{\hat \mu_\omega\}$ (along stable manifolds) and their relation to transverse conditional entropy (see \cref{entformtrans}).

\subsection{Standing hypotheses: finite entropy and integrability}\label{sec:standing ass skew}
For  our main applications below, we will make the following finite entropy assumption: 
\begin{equation*} h _{\hat \nu}(\hat \theta, \ptpart):= H_{\hat \nu}(\what \sigB\mid \hat \theta\inv  \what \sigB)<\infty \tag{H}\label{H}\end{equation*}
and the integrability assumption:
\begin{equation*}
\omega\mapsto \log ^+ \|f_\omega \|_{C^{1+\beta}}, \ \  \omega\mapsto \log ^+ \|f_\omega\inv \|_{C^{1+\beta}}  \in L^1(\hat \nu).
\tag{L1}\label{L1}
\end{equation*}

\subsection{Smooth dynamics of the skew extension}
For this subsection, we need only assume the integrability assumptions \eqref{L1} of \cref{sec:standing ass skew}.
We let $\hat \mu$ be an ergodic, $\what F$-invariant probability measure on $\what X$ with $(\hat \pi_1) _* \hat \mu = \hat \nu$.  
We may apply Oseledec's theorem to the fiberwise derivative cocycle $(  \omega,x)\mapsto D_xf_{ \omega}^{(n)}$ over the action of $\what F$  on $(\what X, \hat \mu)$ and obtain Lyapunov exponents 
$$\lambda_1<\lambda_2<\dots <\lambda_p$$
listed without multiplicity.  
We let $s= \max \bigl(\{0\} \cup \{i: \lambda_i<0\} \bigr )$.  Throughout, we will always assume  $0<s\le p$ unless stated otherwise.  
In particular, the Lyapunov exponents satisfy either 
$$\lambda_1<\lambda_2<\dots<\lambda_s<0\le \lambda_{s+1}<\dots $$
or $$\lambda_1<\lambda_2<\dots<\lambda_s=\lambda_p<0.$$
Also set  $\lambda_{p+1} =+\infty$.  

Let $V^0_{\omega}(x)= \{0\}$ and for $1\le i\le p$,  $$V^i_{\omega}(x) := \left\{v\in T_xM: \limsup_{n\to \infty}  \frac 1 n \log \|D_x f_{\omega }^{(n)}v\|\le \lambda_i
\right\}$$
denote the associated forward  Oseledec subspace.  
The \emph{multiplicity} of $\lambda_i$ is the almost surely constant value of \begin{equation}\label{multiplicity}m_i:=\dim \left(V^i_{\omega}(x)/V^{i-1}_{\omega}(x)\right).\end{equation}

We also write $\lambda_{\btm}(\hat \mu) = \lambda_1$ and $\lambda_{\top}(\hat \mu)= \lambda_p$.  

Given $(\omega,x)\in \what X$ and $1\le i\le s$, let 
$$W^i_{\omega}(x):=\left\{y\in   M:  \limsup_{n\to \infty}  \frac 1 n \log d_M( f_{ \omega}^{(n)}(x),f_{ \omega}^{(n)}(y)) \le \lambda_i\right\}$$
be the \emph{$i$th fast  stable set} in $M$ through $x$ determined by $\omega$. 
We observe that fast stable sets are defined using only forwards dynamics; in particular, given  $(\omega,x)\in X$, we also write  $W^i_{\omega}(x):=W^i_{\tau(\omega)}(x).$

Applying a fibered version of Pesin's stable manifold theorem (cf.\ \cite[Thm.\  III.3.2]{MR1369243} for the case of $C^{1+\mathrm{Lip}}$
 fibered systems or discussion and references in \cite{MR4599404} for the setting of $C^{1+\text{H\"older}}$-fibered systems), we obtain the following.  
\begin{lemma} Assume  the integrability assumption \eqref{L1} of \cref{sec:standing ass skew}.
Then for $1\le  i \le s$ and for 
$\hat \mu$-a.e.\ $(\omega,x)\in \what X$, the set $W^i_{\omega}(x)$ is an injectively immersed, $C^{1+\beta}$ submanifold tangent to $V^i_{\omega}(x)$ at $x$.
\end{lemma}

We also write $$W^i{(\omega,x)}:= \{\omega\}\times W^i_{\omega}(x)$$
for the corresponding \emph{$i$th fast fiberwise stable set} in $\what X$ or $X$ through $(\omega,x)$; we view $W^i{(\omega,x)} $ 
as an injectively immersed submanifold of $\{\omega\}\times M$.  
By definition, the sets $W^i{(\omega,x)}$ form atoms of a partition of $\what X$ or $ X$.  Typically, these partitions are not measurable.  

We  have the following key definitions: 
\begin{definition}\label{def:subord}
Let $1\le i\le s$.  A measurable partition $  \xi$  of $(\what X,\hat \mu)$ is \emph{subordinate} to the partition into fiber $W^i$-manifolds if for $\hat \mu$-a.e.\ $(\omega, x)\in \what X$ 
the following hold:
\begin{enumerate}
\item $\xi(\omega,x)\subset W^i{(\omega,x)}$;  
\item each atom $\xi(\omega,x)$ contains an open (w.r.t.\ the immersed topology) neighborhood of $(\omega,x)$ in $W^i{(\omega,x)}$;
\item each atom $\xi(\omega,x)$ is precompact (w.r.t.\ the immersed topology) in $W^i{(\omega,x)}$.
\end{enumerate}
We say a measurable partition $ \xi$ of $(\what X, \hat \mu)$  (resp.\ $\xi$ of $(  X,   \mu)$) is \emph{Markov} or \emph{sub-invariant}
 if 
\begin{enumerate}[resume]
\item  $\what F\inv \xi \prec \xi$ \ (resp.\ $ F\inv \xi \prec \xi$)
\end{enumerate}
\end{definition}

Adapting \cite{MR693976}, the integrability assumption  \eqref{L1} of \cref{sec:standing ass skew} and standard arguments in the literature ensure the existence of  sub-invariant, measurable partitions $ \hat\xi ^i$ of $\what X=\widehat \Omega \times M$ subordinate  to the partition into fiber  $W^i$-manifolds.  See \cref{sss:1} below for an outline of their construction.

\subsection{Pointwise dimensions and exact dimensionality}
We continue to assume  the integrability assumption \eqref{L1} of \cref{sec:standing ass skew}.
  Given $1\le i\le s$, let $\hat \xi^i$ be a measurable partition of $( \what  X,  \hat \mu)$ subordinate  to the partition into fiber  $W^i$-manifolds.
  Given $(\omega, x)\in \what X$, let $\{\hat \mu^{\hat \xi^i }_{(\omega,x)}\}$ be a family of conditional measures relative  to the partition $\hat \xi^i $.  
   Let $B(x,r)$ denote the ball in $M$ centered at $x$ of radius $r>0$.  By minor abuse of notation, we identify $M$ and $\{\omega\}\times M$ and thus view $\hat \mu^{\hat \xi^i }_{(\omega,x)}$ as a measure on both $M$ and also on $\{\omega\}\times M$.  
   
   \subsubsection{Pointwise dimensions of factor measures}\label{factormeasures}
Let $\hat \mu$ be an ergodic, $\what F$-invariant probability measure on $\what X$ with $(\hat \pi_1)_*\hat \mu= \hat \nu$.   
We define the  $i$th \emph{upper} and \emph{lower pointwise dimensions} of $\hat \mu$  at $(\omega,x)$ to be $$\overline {\dim}^i(\hat \mu, (\omega,x)  ) := \limsup_{r\to 0}\dfrac{\log \left(\hat  \mu^{\hat \xi^i }_{(\omega,x)} (B(x,r)) \right)}{\log r}$$  $$ \underline {\dim}^i(\hat \mu, {(\omega,x)}  ) := \liminf_{r\to 0}\dfrac{\log \left(\hat \mu^{\hat \xi^i }_{(\omega,x)} (B(x,r)) \right)}{\log r}$$
and  the  $i$th \emph{pointwise dimension} of $\hat\mu$ at  ${(\omega,x)}$ to be $$ {\dim}^i(\hat\mu, {(\omega,x)}  ) := \lim_{r\to 0}\frac{\log \left(\hat\mu^{\hat \xi^i }_{(\omega,x)} (B(x,r)) \right)}{\log r}$$  whenever the limit exists.  

One verifies that  the functions
$\overline{\dim}^i(\hat \mu, {(\omega,x)}  )$ and $\underline{\dim}^i(\hat \mu, {(\omega,x)}  )$ are measurable and independent of the choice of $\hat \xi^i$.  
Finiteness of entropy implies the following.
\begin{claim}\label{invdim}Suppose $\furstent_{\hat \mu}(\what F\mid \what \sigF)<\infty.$  Then the functions 
${(\omega,x)}\mapsto \overline{\dim}^i(\hat \mu, {(\omega,x)}  )$ and ${(\omega,x)}\mapsto \underline{\dim}^i(\hat \mu, {(\omega,x)}  )$ are $\what F$-invariant.
\end{claim}
Indeed, in this case we have $(f_{\omega })_*\hat \mu_{\omega }\ll \hat \mu_{\hat \theta({\omega })}$.  Since the maps $f_\omega$ are bi-Lipschitz %and since conditional measures defined relative to $\hat \xi^i $ satisfy equivariance up to restriction and renormalization, 
the local dimension is preserved.   
In particular, ergodicity of $\hat \mu$ implies the upper and lower pointwise dimensions are a.s.\ independent of $(\omega,x)\in \what X$.  Let $\overline{\dim}^i(\hat \mu)$ and $\underline{\dim}^i(\hat \mu  )$ denote their a.s.\ constant values. 

\subsubsection{Exact dimensionality for factor measures}\label{ptwisefactormeasures}
The first main result of this paper is the following exact dimensionality assertion.  
\begin{proposition}[Exact dimensionality of factor measures along stable manifolds]\label{prop:dimensions are defined}  
Assume   the integrability assumption  \eqref{L1} and finite entropy assumption \eqref{H} of \cref{sec:standing ass skew}. Let $\hat \mu$ be an ergodic, $\what F$-invariant probability measure on $\what X$ with $(\hat \pi_1)_*\hat \mu= \hat \nu$.  Then for each $1\le  i \le s$,
$$\overline{\dim}^i(\hat \mu)= \underline{\dim}^i(\hat \mu).$$
\end{proposition}

Write ${\dim}^i(\hat \mu )$ for the common value guaranteed by Proposition \ref{prop:dimensions are defined}.
Set $\dim^0(\hat \mu  )=0$.  
For $1\le i\le s$ define the \emph{$i$th transverse dimension of $\hat \mu$} to be 
$$\gamma^i(\hat \mu ) = \dim^i(\hat \mu )- \dim^{i-1}(\hat \mu).$$

\subsubsection{Generalized setup and exact dimensionality}

As above, let  $\{\hat \xi^i \}$ be a measurable partition of $(\what X, \hat \mu)$ subordinate  to the partition into fiber   $W^i$-manifolds.  Write $\xi^i= \tau\inv (\hat \xi^i)$ for the associated  measurable partition of $(X, \mu)$.  We note that $\xi^i$ is fully $\wtd \sigA$-saturated and thus is not subordinate  to the partition into fiber  $W^i$-manifolds in $(X, \mu)$.  

Let $\zeta$ be a measurable partition of $X$ with the following properties:
\begin{enumerate}
\item $\zeta$ is sub-invariant (that is, $F\inv \zeta\prec \zeta$), and 
\item $\Xi(\sigC)\prec \zeta$.  
\end{enumerate}
We similarly study the dimension properties of the measures $\{ \mu^{\xi^i\vee \zeta }_{(\omega,x)}\}$.

As above, let $B(x, r)$ denote the Riemannian ball of radius $r>0$ centered at $x$ in $M$.  
Given $(\omega,x)\in X$, let $$\wtd B((\omega,x), r):= \tau\inv (\{\tau (\omega)\} \times B(x, r))
= \tau\inv (\{\tau (\omega)\} ) \times B(x, r).$$
That is, $\wtd B((\omega,x), r)$ is the $\wtd \sigA$-saturation of the ball $\{\omega\} \times B(x, r)$ in the fiber $\{\omega\}\times M$.

\begin{definition}We  define the \emph{$i$th upper} and \emph{lower pointwise dimensions} of the measure $\mu$ relative to $\zeta$ to be
\begin{enumerate}
\item $\displaystyle \overline {\dim}^i( \mu, (\omega,x) \mid \zeta ) := \limsup_{r\to 0}\dfrac{\log \left(  \mu^{ \xi^i\vee \zeta }_{(\omega,x)} \wtd B((\omega,x), r)\right)}{\log r}$
\item $ \displaystyle\underline {\dim}^i( \mu, {(\omega,x)} \mid \zeta ) := \liminf_{r\to 0}\dfrac{\log \left( \mu^{
 \xi^i \vee \zeta}_{(\omega,x)} \wtd B((\omega,x), r) \right)}{\log r}.$
\end{enumerate}
We also define the \emph{$i$th pointwise dimension} of the measure $\mu$ relative to $\zeta$  to be 
\begin{enumerate}[resume]
\item $\displaystyle {\dim}^i( \mu, {(\omega,x)} \mid \zeta ) := \lim_{r\to 0}\dfrac{\log \left( \mu^{
 \xi^i \vee \zeta}_{(\omega,x)} \wtd B((\omega,x), r) \right)}{\log r}$, assuming the limit exists. 
\end{enumerate}
\end{definition}

As in \cref{invdim}, the following follows immediately from definition.  
\begin{claim}\label{smelltest}Suppose $h_{\mu}(F, \zeta) <\infty.$  Then 
${(\omega,x)}\mapsto \overline{\dim}^i(  \mu, {(\omega,x)} \mid \zeta )$ and ${(\omega,x)}\mapsto  \underline{\dim}^i(  \mu, {(\omega,x)} \mid \zeta )$ are $F$-invariant.
\end{claim}
Indeed, from the discussion in \cref{finiteentass} below, we have $h_\mu(F, \zeta\vee \xi^i)<\infty$.  

As we always assume $\mu$ to be $F$-ergodic, we write $\overline{\dim}^i(  \mu \mid \zeta )$ and $\underline{\dim}^i(  \mu \mid \zeta )$ for the almost surely constant values of the upper and lower pointwise dimension from \cref{smelltest}.    

Our primary applications for this paper correspond to the  following two cases of this setup:
\begin{enumerate}
\item $\zeta= \Xi(\sigC)$ and $\mu$ is the natural extension of an ergodic, $\what F$-invariant probability measure $\hat \mu$ as guaranteed by \cref{unique natext}.  
In this case, $\overline{\dim}^i(  \mu \mid \zeta )$ and $\underline{\dim}^i(  \mu \mid \zeta )$ recover  dimensions of the factor measure $\hat \mu$ discussed in \cref{factormeasures}.\item $\zeta= \Xi(\sigF)$ and $\mu$ is any ergodic, $F$-invariant probability measure on $X$ with $( \pi_1)_*\mu = \nu$. In this case $\zeta\vee \xi^i$ is a measurable partition of $(X,   \mu)$ subordinate  to the partition into fiber   $W^i$-manifolds.
 Then $\overline{\dim}^i(  \mu \mid \zeta )$ and $\underline{\dim}^i(  \mu \mid \zeta )$ recover the dimensions of the sample measures $\mu_\omega$ along the fiber  $i$th stable manifold.  
For $C^2$ skew extensions, these were shown to exist and their properties studied in \cite{MR2218998,MR2032494}.  
\end{enumerate}

We have the following generalization of \cref{prop:dimensions are defined}  

\begin{proposition}[Generalized exact dimensionality]\label{prop:dimensions are defined general.}  
Assume   the integrability assumption  \eqref{L1}  of \cref{sec:standing ass skew}. 
Let  $\mu$ be an ergodic, $F$-invariant probability measure on $X$ with $( \pi_1)_* \mu= \nu$ and 
let $\zeta$ be a sub-invariant measurable partition of $X$ with $\Xi(\sigC)\prec \zeta$ and $h_\mu(F, \zeta)<\infty$.

 Then for every $1\le i\le s$, $$\overline{\dim}^i(\mu \mid \zeta)= \underline{\dim}^i( \mu\mid \zeta).$$
\end{proposition}
As before, we similarly write  $\dim^0( \mu\mid \zeta )=0$ and   $$\gamma^i( \mu \mid \zeta) = \dim^i( \mu\mid \zeta )- \dim^{i-1}( \mu\mid \zeta).$$
\begin{claim}\label{bdd on dimn 2} In the setting of  \cref{prop:dimensions are defined general.}, we have 
$\gamma^i( \mu \mid \zeta) \le m_i.$
\end{claim}

\begin{remark}\label{finiteentass}We comment on the hypothesis that $h_\mu(F, \zeta)<\infty$ in \cref{prop:dimensions are defined general.} and its relationship to assumption \eqref{H} of \cref{sec:standing ass skew} in \cref{prop:dimensions are defined}

By \cref{prop:GenAR}, if $h_\mu(F, \zeta)<\infty$ we have 
$$\maxent_\mu(F\mid  \sigma( \zeta))= h_\mu(F\mid \sigma( \zeta))+ h_\mu(F, \zeta) .$$  
From a fibered Margulis-Ruelle inequality (see \cite{MR1314494}), the integrability assumption  \eqref{L1} of \cref{sec:standing ass skew} implies that 
$h_\mu(F\mid \sigma( \zeta))=h_\mu(F\mid \sigma( \zeta)^+)\le h_\mu(F\mid \sigF)<\infty$.  
Thus,  under assumption   \eqref{L1} of \cref{sec:standing ass skew}, we have 
that $h_\mu(F, \zeta) <\infty$ implies $\maxent_\mu(F\mid \sigma( \zeta))<\infty$.
  In particular, assuming  \eqref{L1} of \cref{sec:standing ass skew}, we have $\maxent_\mu(F\mid \sigma (\zeta))<\infty$ if and only if $h_\mu(F, \zeta)<\infty$.

  Since $\sigC= (\tau\circ \pi_1)\inv( \what \sigB)$, assumption \eqref{H} of \cref{sec:standing ass skew} implies $h_\mu(F,\Xi(\sigC))=h _{\hat \nu}(\hat \theta, \ptpart)<\infty$ and thus (combined with our standing assumption   \eqref{L1} of \cref{sec:standing ass skew}) 
  we have  $\maxent_\mu(F\mid \sigC)<\infty$.  
  
\end{remark}

\begin{remark}\label{rem:factorsample}
When $\hat \mu$ is an ergodic, $\what F$-invariant probability measure on $\what X$ with $(\hat \pi_1)_*\hat \mu= \hat \nu$ and when $ \mu$ is the natural extension  measure guaranteed by \cref{unique natext}, we will instead use the following notation: 
\begin{enumerate}
\item $ \dim^i(\hat \mu)=\dim^i(\mu\mid \Xi(\sigC))$ and 
\item $\dim^i(\mu)= \dim^i(\mu\mid \Xi(\sigF))$.  
\end{enumerate}
In the latter case, we also write  $$\gamma^i( \mu ) = \dim^i( \mu )- \dim^{i-1}( \mu).$$
\end{remark}

When $\zeta = \Xi(\sigF)$, we recover the exact dimensionality of the sample measures along the filtration by fiber stable manifolds.  This is formulated in  \cite[Thm.\  2.1]{MR2218998}, as a modification of the proof for i.i.d.\ random dynamics in \cite[Prop.\  2.1]{MR2032494}.  We note that as formulated in  \cite{MR2218998,MR2032494}, the maps $f_\omega$ are assumed to be $C^2$; because of an implicit recoding in \cite{MR2218998}  appealing to the separability of $C^2(M)$ and use of Lusin's theorem (see discussion following  \cite[Thm.\ 2.1, p.\ 4]{MR2218998}), 
it seems the exact dimensionality of the measures $\{ \mu^{  \xi^i \vee \Xi(\sigF)}_{(\omega,x)}\}$ is not formally present in the literature  when the maps $f_\omega$ are only assumed to be $C^{1+\text{H\"older}}.$

\subsection{Relationships among dimensions}
Assume   the integrability assumption  \eqref{L1} and finite entropy assumption \eqref{H} of \cref{sec:standing ass skew}. 
Let $\hat \mu$ be an ergodic, $\what F$-invariant probability measure on $\what X$ with $(\hat \pi_1)_*\hat \mu= \hat \nu$ and  let $\mu$ be the unique ergodic $F$-invariant probability measure on $X$ guaranteed by \cref{unique natext}.

Following the notation in \cref{rem:factorsample}, we have the following relationships between the above definitions.  

\begin{proposition}\label{claim:bound on transverse dimension}
For $1\le i\le s$, 
\begin{enumerate}
\item $\dim^i(\mu)\le \dim^i(\hat \mu)$, and 
\item $\gamma^i( \mu )\le \gamma^i(\hat \mu ) \le m_i.$ 
\end{enumerate}
\end{proposition}

\subsection{Entropy formulae and invariance principle} 
Following the seminal work of Ledrappier and Young \cite{MR819557,MR819556}, relations between the pointwise dimensions of $\mu$, (standard) conditional entropy, and Lyapunov exponents are well established in the existing literature.  We especially emphasize \cite{MR953818,MR968818,MR2218998,MR1369243,MR2032494} for results on  standard conditional entropy for random dynamics by  $C^2$-diffeomorphisms.  We recover the results of \cite{MR2218998,MR2032494} as well as new entropy formulas for notions of entropy we introduced above.

\subsubsection{Entropy formula for standard conditional entropy $h_\mu(F\mid  \sigF)$}
\begin{theorem}[{cf. \cite[Thm.\ 2.1]{MR2218998} and  \cite[Thm.\ 2.3]{MR2032494}}]
\label{entformtang}
Assume   the integrability assumption  \eqref{L1} of \cref{sec:standing ass skew}. 
Let $\mu$ be an ergodic, $F$-invariant probability measure with $( \pi_1)_*\mu= \nu$.  Then 
$$h_\mu(F\mid  \sigF)= h_{\hat \mu}(\what F\mid  \what \sigF) =h_\mu(F\mid\sigC)= -\sum_{j=1}^s \gamma^j(\mu) \lambda_j.$$
(Above, $\hat \mu = \tau_*\mu$.) 
\end{theorem}
As the above is not formally shown for actions by $C^{1+\text{H\"older}}$ diffeomorphisms in the literature, we will provide a proof of \cref{entformtang} unified with the proofs of our main new results.  Indeed, \cref{entformtang}  is simply a special case of \cref{thm:realthm}, below, with $\zeta = \Xi(\sigF)$.

\subsubsection{Entropy formula for transverse conditional entropy $\transent_{\hat \mu}(\what F\mid  \what \sigF)$}
Similarly, a special case of   \cref{thm:realthm}, below, with $\zeta = \Xi(\sigC)$ recovers the following entropy formula for transverse conditional entropy.
\begin{theorem}\label{entformtrans}
Assume   the integrability assumption   \eqref{L1} and finite entropy assumption \eqref{H} of \cref{sec:standing ass skew}.
Let $\hat \mu$ be an ergodic, $\what F$-invariant probability measure on $\what X$ with $(\hat \pi_1)_*\hat \mu= \hat \nu$.  
Then 
$$\transent_{\hat \mu}(\what F\mid \what  \sigF)=  -\sum_{j=1}^s \gamma^j(\hat \mu)\lambda_j.$$
\end{theorem}

\subsubsection{Entropy formula for Furstenberg entropy $\furstent_{\hat \mu}(\what F\mid \what\sigF)$}
From \cref{entropygap,entformtrans,entformtang}, we immediately obtain the following.
\begin{theorem}\label{entformskewfurst}
Assume   the integrability assumption   \eqref{L1} and finite entropy assumption \eqref{H} of \cref{sec:standing ass skew}.
Let $\hat \mu$ be an ergodic, $\what F$-invariant probability measure on $\what X$ with $(\hat \pi_1)_*\hat \mu= \hat \nu$ and  let $\mu$ be the unique ergodic $F$-invariant probability measure $\mu$ on $X$ guaranteed by \cref{unique natext}.  Then 
\begin{align*}
\furstent_{\hat \mu}(\what F\mid \what\sigF)&= \transent_{\hat \mu}(\what F\mid \what \sigF)-  h_{\hat   \mu}(\what F\mid \what  \sigF)\\
&= \transent_{\hat \mu}(\what F\mid \what \sigF)-  h_{  \mu}(  F\mid    \sigF)\\
&=-\sum_{i=1}^s \lambda_i (\gamma^i(\hat \mu) - \gamma^i(\mu)).
\end{align*}
\end{theorem}
We only state and prove \cref{entformskewfurst} under our standing assumption that  $s>0$ and so  there is at least one negative exponent.  In principle, the same proof recovers the classical (sufficient but not necessary criterion) {\it invariance principle} in our setting asserting that $\furstent_{\hat \mu}(\what F\mid \what\sigF)=0$ whenever $\lambda_{\btm}(\hat \mu)\ge 0$. 

From a fibered Margulis-Ruelle inequality (see \cite{MR1314494}), if all exponents of $\hat \mu$ are negative, then the standard conditional entropy vanishes, $h_{\hat   \mu}(\what F\mid \what  \sigF)=0$.  It follows that $\gamma^i(\mu)=0$ for every $1\le i\le s$ and we obtain the following entropy formula for measures $\hat \mu$ that are contracting on average in every fiber.  
\begin{corollary}\label{cor:contracting}
Assume   the integrability assumption   \eqref{L1} and finite entropy assumption \eqref{H} of \cref{sec:standing ass skew}.
Let $\hat \mu$ be an ergodic, $\what F$-invariant probability measure on $\what X$ with $(\hat \pi_1)_*\hat \mu= \hat \nu$.  
Assume $\lambda_{\top}(\hat \mu)<0$.
Then 
$$\furstent_{\hat \mu}(\what F\mid \what  \sigF)=  -\sum_{j=1}^s \gamma^j(\hat \mu)\lambda_j.$$

\end{corollary}

\subsubsection{Rigidity and invariance principle}
We have the following extension of the invariance principle  which follows directly from \cref{entformskewfurst}, \cref{claim:bound on transverse dimension}, \cref{entropygap}, \cref{prop:basicineq}, and \cref{claim:zero-furstenberg-entropy}.
\begin{theorem}\label{invprinc:general}
Assume   the integrability assumption   \eqref{L1} and finite entropy assumption \eqref{H} of \cref{sec:standing ass skew}.
Let $\hat \mu$ be an ergodic, $\what F$-invariant probability measure on $\what X$ with $(\hat \pi_1)_*\hat \mu= \hat \nu$ and  let $\mu$ be the unique ergodic $F$-invariant probability measure $\mu$ on $X$ guaranteed by \cref{unique natext}. 
The following are equivalent:
\begin{enumerate}
\item $\furstent_{\hat \mu}(\what F\mid \what\sigF)= 0$;
\item $h_{\hat \mu}(\what F\mid  \what \sigF) = \transent_{\hat \mu}(\what F\mid  \what \sigF)$;
\item for each $1\le i\le s$, $\dim^i(\mu) = \dim^i(\hat \mu)$;
\item for each $1\le i\le s$, $\gamma^i(\mu) = \gamma^i(\hat \mu)$;
\item the map $\omega \mapsto \mu^\sigF_\omega$ is $\sigA$-measurable on $\Omega$;
\item for $\nu$-a.e.\ $\omega$,  $\mu^\sigF_\omega = \hat \mu^{\what \sigF}_{\tau( \omega)}$ (viewed as probability measures on $M$);
\item for $\hat \nu$-a.e. $\omega\in \what \Omega$, $(f_\omega)_* \hat \mu_\omega = \hat \mu_{\hat \theta(\omega)}$;
\item $H_\mu(\wtd \sigA\mid F\inv \wtd\sigA)= H_\nu( \sigA\mid \theta \inv \sigA).$
\end{enumerate}
\end{theorem}

\subsection{Main technical theorem} 

The main technical result of our paper is the following adaptation of the  results of Ledrappier and Young from \cite{MR819556,MR819557}.
\begin{theorem}\label{thm:realthm}
Assume the integrability assumption \eqref{L1} of \cref{sec:standing ass skew}. 
Let $\mu$ be an ergodic, $F$-invariant probability measure on $  X$ with $ ( \pi_1)_*  \mu=   \nu$ and let $\zeta$ be a sub-invariant measurable partition with $\Xi(\sigC)\prec \zeta$.  

Suppose that $ h_\mu(F, \zeta)<\infty$.  Then 
 $$\maxent_\mu(F\mid \sigma(\zeta))= h_{\mu}(F, \xi^0\vee \zeta) -\sum_{j=1}^s \gamma^j( \mu\mid \zeta) \lambda_j.$$
\end{theorem}

Specializing to our main setting of interest where $\zeta= \Xi(\sigC)$, we have $\xi^0\vee \zeta= \xi^0$ and since $\xi^0= \tau\inv \ptpart$ is the preimage of the point partition on $\what X$,  we have $$h_{\mu}(F, \xi^0\vee \zeta)= H_{\hat \mu}(\ptpart\mid\what  F\inv \ptpart)$$ and deduce  the following.  
\begin{corollary}\label{thm:realthm2}
Assume   the integrability assumption  \eqref{L1} and finite entropy assumption \eqref{H} of \cref{sec:standing ass skew}. 
Let $\hat \mu$ be an ergodic, $\what F$-invariant probability measure on $\what X$ with $(\hat \pi_1)_*\hat \mu= \hat \nu$ and  let $\mu$ be the unique ergodic $F$-invariant probability measure $\mu$ on $X$ guaranteed by \cref{unique natext}.
Then $$\maxFentC= H_{\hat \mu}(\ptpart\mid\what  F\inv \ptpart)-\sum_{j=1}^s \gamma^j(\hat \mu)\lambda_j.$$
\end{corollary}

We remark that  the statement (and idea of proof) of \cref{thm:realthm2} is similar to the main result of \cite{MR2530169} for $C^2$ endomorphisms if we replace the quantity $\maxFentC$ by the total metric entropy $h_\mu(F)$ of a ($C^2$, non-degenerate) endomorphism $F\colon M\to M$ (i.e., a self-covering map.)   Indeed, one may view the term $ H_{\hat \mu}(\ptpart\mid\what  F\inv \ptpart)$ in  \cref{thm:realthm2} as a replacement for the ``folding entropy'' for endomorphisms; in this case, \cref{thm:realthm2} is rather similar to
\cite[Thm.\  1.2]{MR2530169}.  
We note however that the results in \cite{MR2530169} are only formulated for $C^2$ endomorphisms and also rely on  \cite[Prop.\ 4.1.2]{MR2448134} which is only formally proved in the case of hyperbolic measures (i.e., measures without zero Lyapunov exponents) with indications for modifications for general measures.  We provide complete details for our \cref{thm:realthm}, establishing \cref{prop:entropyisunstable} below which provides a replacement of \cite[Prop.\ 4.1.2]{MR2448134}.

\subsection{\texorpdfstring{Proofs of \cref{entformtrans,entformtang}}{Proofs of Theorems 4.13 and 4.14}}
\cref{entformtrans,entformtang}   follow directly from \cref{thm:realthm}.  
Indeed, in the setting of \cref{entformtang}, we set $\zeta=  \Xi(\sigF)$ which is strictly invariant and so $h_{\mu}(F, \xi^0\vee \zeta)= 0$.
Then \cref{entformtang} follows directly from  \cref{thm:realthm} after appropriate notation changes.

In the setting of 
\cref{entformtrans}, 
by \cref{prop:howtoNatExt} we have $\maxFentC= \maxent_{\hat \mu} (\what F\mid \what \sigF)$.  
We apply  \eqref{eq:factAR} of \cref{prop:GenAR} and, using that the restriction of $\what F$ to every atom of $\what \sigF$ is a diffeomorphism (and thus injective), we also apply \cref{eq:thiseq}  to conclude that
\begin{align*}
\transent_{\hat  \mu}(\what F\mid  \what \sigF) &= 
\maxent_{ \hat  \mu}(  \what F\mid    \what \sigF)- H_{\hat \mu}(\what \sigF\mid\what  F\inv \ptpart) \\&=
\maxent_{   \mu}(  F\mid    \sigC)- H_{\hat \mu}(\ptpart\mid\what  F\inv \ptpart) 
\end{align*}
and so \cref{entformtrans} then follows directly  from \cref{thm:realthm2}.

The remainder of the paper is devoted to the proof of  \cref{thm:realthm}.  

\subsection{Exact dimensionality of projected fiber measures}
Let $\hat \mu$ be an ergodic, $\what F$-invariant probability measure on the one-sided skew extension $\what X= \what \Omega\times M$.  
When all fiberwise Lyapunov exponents of $\hat \mu$ are negative, \cref{ptwisefactormeasures} implies the fiberwise conditional measures  $\hat \mu_\omega$ are exact dimensional for a.e.~$\omega$.  To obtain exact dimensionality of the measures $\hat \mu_\omega$, roughly, we view the fiberwise conditional measures $\{\hat \mu_\omega:\omega\in \what \Sigma\}$ for the one-sided skew shift as images under a  quotient $\tau \colon \Sigma\times M\to \what \Sigma\times M$  along a system of stable holonomies.  

Let $\mu_0 =( \hat \pi_2)_* \hat \mu$ denote the projection of the measure $\hat \mu$ to $M$.  When this projection can be viewed as a projection along a system of unstable holonomies, we similarly show $\mu_0$ is exact dimensional with the same a.s.~dimension as $\hat \mu_\omega$.  
A  \emph{one-sided generator} for  the abstract measurable dynamical system $\hat \theta\colon (\what \Omega, \hat \nu)\to (\what \Omega, \hat \nu)$ is a finite entropy partition $\partP$ of $(\what \Omega, \hat \nu)$ with 
$ \partP^{-} = \ptpart$.  
Given a one-sided generator $\partP$, we have $h_{\hat \nu}(\hat \theta)= h_{\hat \nu}(\hat \theta,\partP)\le H_{\hat \nu}(\partP).$
In the case of the one-sided shift $\what \Sigma$ (on a finite alphabet), the partition into $0$-cylinders is a one-sided generator.  
We formalize the notion of the projection being taken relative to a system of unstable holonomies by requiring the maps $\omega\mapsto f_\omega$ being constant on atoms of such a partition.  

\begin{theorem}[Exact dimensionality of projection to fiber]\label{thm:exactdimproj}
Let $\hat \mu$ be an ergodic, $\what F$-invariant probability measure on the one-sided skew extension $\what X= \what \Omega\times M$ and let $\hat \nu=(\hat \pi_1)_* \hat \mu $.  
Suppose the measurable system $\hat \theta\colon (\what \Omega, \hat \nu)\to (\what \Omega, \hat \nu)$ has finite entropy $h_{\hat \nu}(\hat \theta)<\infty$ and admits a one-sided generator $\hat \partQ$.  Suppose the function $\omega\mapsto f_\omega$ is $\sigma(\hat \partQ)$-measurable.  Suppose also the  integrability assumption  \eqref{L1} of \cref{sec:standing ass skew} and suppose that $\lambda_{\top}(\hat \mu)<0$.  

Then there exists $0\le d\le \dim M$ such that 
\begin{enumerate}[label*={(\alph*)}]
\item \label{projdima} for $\hat \nu$-a.e.~$\omega$, the measure $\hat \mu_\omega$ is exact dimensional with dimension $d$, and
\item \label{projdimb} the measure $\mu_0 =(\hat \pi_2)_*\hat \mu$ is exact dimensional with dimension equal to $d$.  
\end{enumerate}
\end{theorem}
The proof of \cref{thm:exactdimproj} appears in \cref{sss:proj msr}.  Roughly, the result follows from the  ``almost product structure'' of hyperbolic measures in smooth ergodic theory.  The ``almost product structure''  was established in generality in \cite{MR1709302} extending the case of surfaces in \cite{MR0684248}.  The assumption that $\hat \theta\colon (\what \Omega, \hat \nu)\to (\what \Omega, \hat \nu)$ admits a one-sided generator ensures the measure $\mu_0$ is the quotient under the ``expanding directions'' of the dynamical system.  
A similar  argument appears in section 6 (see especially (C2)) of Feng's paper \cite{MR4557759} where exact dimensionality of the measure $\mu_0$ is established simultaneously with the exact dimensionality of Feng's analogue of our measures $\hat \mu_\omega$.  
Although the  proof of \cref{thm:exactdimproj} is essentially included in the proof of our main results (compare \cref{sss:ineqI} and \cref{sss:proj msr}), we  choose to separate the proof of  \cref{thm:exactdimproj} from our main abstract theorems  to emphasize the additional structure needed  in \cref{thm:exactdimproj}.

\section{Preparatory constructions} \label{sss:preparatory}
We keep all notation introduced in \cref{formulation and results}.  In particular, we consider the fibered dynamics  $\what F\colon \what X\to \what X$ and an ergodic, $\what F$-invariant probability measure $\hat \mu$ on $\what X$ fibering onto $\hat \nu$.  Then $F\colon X\to X$ represents the ergodic theoretic natural extension preserving an ergodic, $ F$-invariant probability measure $ \mu$ on $ X$ fibering onto $ \nu$.

We revisit the majority of the preparatory lemmas and constructions developed to prove the main results of \cite{MR819556,MR819557}.  
\subsection{\texorpdfstring{Sub-invariant partitions subordinate to $W^i$-manifolds}{Sub-invariant partitions subordinate to Wⁱ-manifolds}}\label{good parts}\label{sss:1}

Recall we assume the Lyapunov exponents (for the cocycle $D_xf^{(n)}_\omega$ over $F\colon X\to X$) indexed without multiplicity satisfy 
 $$\text{$\lambda_1<\lambda_2<\dots<\lambda_s<0\le  \lambda_{s+1}$ \quad or \quad  
$\lambda_1<\lambda_2<\dots<\lambda_s=\lambda_p<0.$}$$

Write  $\hat \xi^0=\ptpart$ for the point partition on $\what X$.  
 Following the constructions in \cite[Prop.\  3.1]{MR693976},  \cite[Lem.\  3.1.1]{MR819556}, and  \cite[Lem.\  9.1.1]{MR819557}, for  $1\le i\le s$,  there exists a  measurable partition $\hat  \xi^i$ of $(\what X, \hat \mu)$  with the following properties:   
\begin{enumerate}[{label*={(\arabic*)}, ref= {(\arabic*)} }]
	\item \label{kkk1} $ \hat \xi^i$  is subordinate  to the partition into fiber  $W^i$-manifolds (see \cref{def:subord});  
	\item  $ \what F\inv \hat \xi^i\prec \hat\xi^i$;
	\item  \label{kkk3} for $0\le i\le s-1$, we have $\hat \xi^{i+1}\prec \hat \xi^{i}$.  
\end{enumerate}
See \cite[Prop.\  IV.2.1]{MR1369243} where such partitions are constructed for one-sided i.i.d.\ random dynamics.  We will outline their construction in \cref{sec:sss} below.  

Recall $\tau\colon X\to \what X$ is the  factor map of natural extensions.  Throughout, we will always write $$\xi^i:= \tau\inv \hat \xi^i$$ for the associated measurable partition of $(X,\mu)$.      We note that $\xi^i$ is saturated by $\wtd \sigA$-atoms and so is not subordinate to fiber $W^i$-manifolds.  
However, we still have that $\xi^i$ satisfy the following:
\begin{enumerate}[,label*={(\arabic*)}, ref= {(\arabic*)} ]
\item   $   F\inv   \xi^i\prec  \xi^i$;
	\item  \label{kkk6}$\bigvee_{n=0}^\infty F^{n}(  \xi^i)=\ptpart$ is the point partition on $X$;
	\item   for $0\le i\le s-1$, we have $  \xi^{i+1}\prec   \xi^{i}$.  
\end{enumerate}
Note that property \ref{kkk6} of the partition $\xi^i$ follows from recurrence and using that each atom $\hat \xi^i(x)$  is precompact in $W^i{(\omega,x)}$ for a.e.\ $(\omega, x)\in \what X$.

\subsection{\texorpdfstring{Entropy properties of partitions in \cref{good parts}}{Entropy properties of partitions in §5.1}}
Let $\hat \xi^i$ be a partition as in \cref{good parts}.  Recall we write  $\xi^i:= \tau\inv \hat \xi^i$. We first observe the following which follows from definition.  

\begin{claim}For every $0\le  i \le s$, we have 
$h_{\hat \mu}(\what F, \hat \xi^i) = h_{ \mu}( F,  \xi^i). $
\end{claim}

We now let  $\zeta$ be a sub-invariant measurable partition of $(X,\mu)$ with $\Xi(\sigC)\prec \zeta$ and $h_\mu(F, \zeta)<\infty$.  This, in particular, implies (see \cref{finiteentass}) that $\maxent_\mu(F\mid  \sigma( \zeta))<\infty$.  
We will study the contribution to the entropy $\maxent_\mu(F\mid  \sigma( \zeta))$ witnessed by the partitions of the form $\xi^i:= \tau\inv \hat \xi^i$. 

As in \cite[Lem.\  3.1.2]{MR819556}, 
all future entropy calculations will be independent of various choices made (in \cref{sss:outling} below) in the construction of the partitions in \cref{good parts}.
\begin{lemma}[{cf.\  \cite[Lem.\  3.1.2]{MR819556}}] \label{lem:sameentropy}
For each $1\le i\le s$,  let $  \xi_1^i$ and $   \xi_2^i$ be two  partitions of $(  X,  \mu)$ as in \cref{sss:1}.   
Then $$h_{  \mu}(  F,   \xi_1^i\vee \zeta)= h_{  \mu}(  F,   \xi_2^i\vee \zeta).$$
\end{lemma}

Similarly, \cref{Prop:GenARpart}
implies the partitions constructed in \cref{sss:1} ``locally maximize'' entropy.  
\begin{lemma}\label{cor:entropysame}
Let $\xi$ be any measurable partition of $(X,\mu)$ satisfying
$   F\inv   \xi\prec  \xi$ and $\bigvee_{n=0}^\infty F^{n}(  \xi)=\ptpart$.  
 Let  $\partP$ be   measurable partition with  
 $H_\mu(\partP\mid \xi\vee \zeta)<\infty$.  Then 
$$h_\mu(F, \xi\vee \zeta\vee \partP )= h_\mu(F, \xi\vee \zeta).$$
In particular, 
$$h_\mu(F, \xi\vee \zeta)= \maxent_\mu(F \mid \sigma(\xi\vee \zeta)).$$
\end{lemma}

\cref{lemma:pastfinitepartabstract} and \cref{cor:entropysame} then give the following.
\begin{lemma}\label{lemma:pastfinitepart}
Let $ \calP$ be a measurable partition of $(  X,  \mu)$ with $H_\mu( \calP\mid  \xi^i\vee \zeta)<\infty$. 
Then for $ \mu$-a.e.\ $(\omega ,x)$,
$$\lim_{n \to \infty} - \frac 1 n \log  { \mu }^{  \xi^i\vee \zeta}_{(\omega, x)} \left( (   \xi^i\vee \zeta \vee   \calP) ^{(-n,0]}(\omega, x) \right)  = h_{ \mu }(   F, \xi^i \vee \zeta ).$$
\end{lemma}

The following proposition is the analogue of the main technical result (Corollary 5.3, corollary of Proposition 5.1) of \cite{MR819556}.  
Roughly, the proposition says that the refinement entropy is  captured by the partition 
$\zeta\vee\xi^s$ where $\xi^s=\tau\inv\hat \xi^s$ for  a sub-invariant partition  $\hat \xi^s$ subordinated to the full stable manifolds whose properties are outlined in \cref{good parts}.     
  The proof will  follow directly from   the inequalities in \cref{sec:outline}  below. 

\begin{proposition}\label{prop:entropyisunstable}
Assume   the integrability assumption  \eqref{L1} of \cref{sec:standing ass skew}. 
Let $  \mu$ be an ergodic, $  F$-invariant probability measure on $  X$ with $ ( \pi_1)_*  \mu=   \nu$ and let $\zeta$ be a sub-invariant measurable partition of $X$ with $\Xi(\sigC)\prec \zeta$. 

Suppose that $  h_\mu(F, \zeta)<\infty$.  Then 
		  $$ \maxent_\mu(F\mid \sigma(\zeta)) =  h_\mu( F,  \xi^s\vee \zeta).$$

		 In particular, under the integrability assumption  \eqref{L1} we have 
 		 $$ \maxent_\mu(F\mid \sigF) =h_\mu(F\mid \sigF) = h_\mu(  F,\xi^s\vee \Xi(\sigF))$$
		 and  under the integrability assumption  \eqref{L1} and finite entropy assumption \eqref{H} of \cref{sec:standing ass skew} we have
 		 $$ \maxent_\mu(F\mid \sigC) 
		 =  h_\mu( F, \xi^s\vee \Xi(\sigC))
		 = h_{\hat \mu}(\what F, \hat \xi^s\vee \Xi({\what {\sigF}})).$$ 
  \end{proposition}

\subsection{\texorpdfstring{Outline of construction of partitions in \cref{good parts}}{Outline of construction of partitions in §5.1}}
\label{sss:outling}\label{sec:sss}
To construct local (and thus global) stable manifolds $W^i_\omega(x)$ parametrized by the one-sided skew-shift $(\omega,x)\in \what X$, one invokes one-sided Lyapunov charts and the graph transform method for one-sided dynamics.  See  for instance \cite[]{MR968818} for a formulation of one-sided Lyapunov charts.   However,  later arguments we employ require   the use of two-sided Lyapunov charts.  We thus only summarize the outcomes of the constructions of local stable-manifolds using one-sided charts.

Given $\omega\in \widehat \Omega$ and $y,z\in M$, write $d^s_\omega(y,z)$ for the induced Riemannian distance from $y$ to $z$ along $W^s_\omega(y)$; here we define 
$d^s_\omega(y,z)=\infty$ if $W^s_\omega(y)$ is not an injectively immersed submanifold or if $z\notin  W^s_\omega(y)$.
Moreover, we note that for $\hat \mu$-a.e.\ $(\omega,x)\in \what X$, we have that $W^{i}_{\omega}(x)$ is a properly embedded (with respect to immersed topologies) submanifold in  $W^s_\omega(x)$ and so the restriction of $d^s_\omega(\cdot ,\cdot )$ to $W^i_\omega(x)$ induces the immersed topology on $W^i_\omega(x)$.

Given $\hat \mu$-a.e.\ $(\omega, x)\in \what X$, there is an embedded $C^{1+\beta}$ disc $W^{s}_{\loc,{\omega}}(x)$ in $M$ such that the following hold:
\begin{enumerate}
\item $W^{s}_{\loc,{\omega}}(x)$ is an open (w.r.t.\ immersed topology) subset of $W^{s}_{{\omega}}(x)$.
\item For all such $(\omega,x)$, every $z\in W^{s}_{\loc,{\omega}}(x)$, and every $1\le  i \le s$, 
there is a $C^1$ embedded open disk $\bar W^{i}_{\loc,{(\omega,x)}}(z)$  characterized by 
  $$ \bar W^{i}_{\loc,{(\omega,x)}}(z) := \left\{y\in W^{s}_{\loc,{\omega}}(x) :
\limsup_{n\to \infty} \frac 1 n \log d^s_{\hat \theta^n{\omega}}(f_\omega^{(n)}(y), f_\omega^{(n)}(z))\le \lambda_i  \right\}.$$
\item The collection $\left \{\bar W^{i}_{\loc,{(\omega,x)}}(z) : z\in W^{s}_{\loc,{\omega}}(x) \right \}$ forms a continuous foliation of $W^{s}_{\loc,{\omega}}(x)$  with trivial holonomy.  
\item Both the  family of $C^1$ discs $(\omega,x)\mapsto W^{s}_{\loc,{\omega}}(x)$ and their associated foliations  $\left\{\bar W^{i}_{\loc,{(\omega,x)}}(z) : z\in W^{s}_{\loc,{\omega}}(x) \right\}$ depend measurably on $(\omega,x)\in \what X$.  
\end{enumerate}

Let $\Lambda_1\subset \what X$ be a compact set such that 
 \begin{enumerate}
 \item $W^{s}_{\loc,{\omega}}(x)$ is a  $C^1$ embedded disc for every $(\omega, x)\in \Lambda_1$, 
 \item if $(\omega, x)\in \Lambda_1$ and $(\omega, y)\in \Lambda_1$, then $W^{s}_{\loc,{\omega}}(y)\cap W^{s}_{\loc,{\omega}}(x)=\emptyset$, and 
 \item  the family  $(\omega,x)\mapsto W^{s}_{\loc,{\omega}}(x)$ of $C^1$ embedded discs varies continuously on $\Lambda_1$;
\end{enumerate}
moreover, by setting $\Lambda_2= \bigcup_{(\omega,x)\in \Lambda_1} W^{s}_{\loc,{\omega}}(x)$, we assume
 \begin{enumerate}[resume]
\item $\hat \mu(\Lambda_2)>0$, and  
\item for each $1\le  i \le s$,  the family of $C^1$ embedded discs $(\omega,z)\mapsto \bar W^{i}_{\loc,{(\omega,x)}}(z)$ varies continuously on $\Lambda_2$.
\end{enumerate}

   Recall $\hat \pi_2\colon \what X\to M$ denotes the natural projection.  Then $(\hat \pi_2)_*(\restrict{\hat \mu}{\Lambda_2})\ll (\hat \pi_2)_*{\hat \mu}.$ Let $x_0\in M$ be in the support of $(\hat \pi_2)_*{\hat \mu}$ and such that the Radon--Nikodym  derivative (as computed via Lebesgue-Besicovitch differentiation in $M$) $\frac{d(\hat \pi_2)_*(\restrict{\hat \mu}{\Lambda_2})}{ d(\hat \pi_2)_*{\hat \mu}}$ is positive at $x_0$; then $(\hat \pi_2)_*(\restrict{\hat \mu}{\Lambda_2})(B(x_0,\rho))>0$ for all $\rho>0$.

Given $\rho>0$, let   \begin{equation}\label{bighatballs}
  S_\rho= \what B_{\Lambda_2}(x_0, \rho):= \Lambda_2\cap \hat \pi_2\inv( B(x_0, \rho)).
   \end{equation}
Given ${(\omega,y)}\in S_\rho$, there is a unique   $(\omega, x)\in \Lambda_1$ such that $y\in W^{s}_{\loc,{\omega}}(x)$;
let $$V^{i}_{\omega,\rho}(y)=\bar W^{i}_{\loc,{(\omega,x)}}(y)\cap B(x_0, \rho)$$
and let 
$$V^{i}_{(\omega,y),\rho}= \{\omega\}\times V^{i}_{\omega,\rho}(y).$$

Taking $\rho_0>0$ sufficiently small, we may assume for $0<\rho<\rho_0$ and all $(\omega,y)\in S_\rho$, by taking $(\omega,x)\in \Lambda_1$ with $y\in W^{s}_{\loc,{\omega}}(x)$, for 
every $1\le i\le s$ that the following hold:  
\begin{enumerate}[label*={(\alph*)}]
	\item $V^{i}_{(\omega,y),\rho}$ is a (possibly non-connected) open neighborhood of $y$ with compact closure  in  $\{\omega\}\times \bar W^{i}_{\loc, (\omega,x)}(y);$ 
%	\item $V^{i}_{(\omega,y),\rho}$ is a connected open disk containing $y$ with compact closure  in  $\{\omega\}\times \bar W^{i}_{\loc, (\omega,x)}(y);$ 
	\item there is $\kappa>0$ such that if $z\in\bar  W^i_{\loc,(\omega,x)}(y)$, if  $(\omega,z)\in S_\rho$, and if  $d^s_\omega(y,z)<\kappa$ 
	 then $V^{i}_{(\omega,y),\rho}=V^{i}_{(\omega,z),\rho}$.
\end{enumerate}
%%%Note that for every $(\omega,x)\in S_\rho$ and every $y\in  V^s_{\omega,\rho}(x)$, 
%%%the set $V^{i}_{(\omega,y),\rho}$ is well defined and 
%%%the collection $\{V^{i}_{(\omega,y),\rho}: y\in V^{s}_{\omega,\rho}(x)\}$ forms a continuous foliation by $C^1$ leaves of the open disk $\{\omega\}\times V^{s}_{\omega,\rho}(x)$.

For  $(\omega',y)\in S_\rho$, set $D^i(\omega',y ) = V^{i}_{(\omega', y),\rho}= \{\omega'\}\times V^{i}_{\omega',\rho}(y)$. 
Then $$\bar  \xi^i(\omega',y) = \begin{cases}  D^i(\omega',y) & (\omega',y)\in S_\rho,  \\ \what X \sm S_\rho & (\omega',y) \notin S_\rho\end{cases}$$ 
defines a partition of a full measure subset of $(\what X, \hat \mu)$.  
Take $$\hat \xi^i:= (\bar \xi^i)^-\vee \Xi
(\what \sigF).$$

The partition  $\hat \xi^i$ clearly is sub-invariant.
By ergodicity, for any $0<\rho<\rho_0$ and a.e.\ $(\omega',y)$, we have that $\what F^{n}(\omega',y) \in S_\rho$ for some $n\ge 0$.  
It follows that  $\hat \xi^i(\omega',y) \subset W^i(\omega',y)$ for a.e.\ $(\omega',y)$.  
Since $D^{i-1}(\omega',y)\subset D^{i}(\omega',y)$ for every $(\omega',y)\in S_\rho$, we have $\bar  \xi^i\prec \bar  \xi^{i-1}$ and thus $\hat  \xi^i\prec \hat  \xi^{i-1}$.  
Since $\what F$ restricts to a diffeomorphism on every fiber, $\hat \xi^i(\omega,x)\subset W^i{(\omega,x)}$ and is precompact in $W^i{(\omega,x)}$.
By the same arguments as in \cite[Section 3]{MR693976}, for Lebesgue-almost every choice of $\rho<\rho_0$, $\hat \xi^i(\omega,x)$ contains an open (w.r.t.\ the immersed topology) neighborhood of $(\omega,x)$ in $W^i{(\omega,x)}$.  
Properties \ref{kkk1}--\ref{kkk3}  in \cref{sss:1} then hold for each $1\le i\le s$ for such a choice of $\rho<\rho_0$.

 Note that 
$\hat \xi^i$ has the following characterization: 
 \begin{claim} \label{in same tilde chart}There is $\kappa>0$ such that for all $0<\rho<\rho_0$ the following are equivalent for almost every $(\omega, x)\in \what X$ and $y\in W^i_\omega(x)$:
 \begin{enumerate}
 \item  for any  $m\ge 0$,
 $\what F^{m}(\omega,y)\in S_\rho$ if and only if $\what  F^{m}(\omega, x) \in S_\rho$;  moreover, for all such $m$,
$d_{\hat \theta^m(\omega)}^s(f_\omega^{(m)}(x), f_\omega^{(m)}(y))<\kappa$.
\item  $\hat \xi^i(\omega,x) = \hat\xi^i(\omega,y)$.
\end{enumerate}
\end{claim}

This implies the following.  
\begin{corollary}\label{cor:refinedynamically}
For $n\ge 0$, we have $\hat \xi^i\vee \what F^{-n}\hat \xi^{i-1}= \hat \xi^{i-1}$.  
\end{corollary}
\begin{proof}
If $(\omega,z)\in (\hat \xi^i\vee \what F^{-n}\hat \xi^{i-1})(\omega,y)\subset (\hat \xi^i)(\omega,y)$ then 
for all $m\ge 0$, we have 
$\what F^m(  \omega, z)\in S_\rho$ if and only if 
$\what F^m(  \omega, y)\in S_\rho$;  
moreover for all such $m$,  we have that
$d_{\hat \theta^m(  \omega)}^s(f_{  \omega}^{(m)}(z), f_{  \omega}^{(m)}(y))<\kappa$.  
Since $(\what F^{-n} \hat \xi^{i-1}\vee\hat \xi^{i})(\omega,y)\subset 
 (\what F^{-n} \hat \xi^{i-1}\vee \Xi(\what \sigF))(\omega,y)
\subset W^{i-1}(\omega,y)$, \cref{in same tilde chart} implies  $(\omega,z)\in \hat \xi^{i-1}(\omega,y)$, showing the non-trivial inclusion of sets.  
\end{proof}

\subsection{Two-sided Lyapunov charts and related constructions}\label{sec:twosidecharts}
We review the properties of two-sided Lyapunov charts and the standard  ``globalized'' dynamics constructed relative to such charts.  Relative to the choice of globalized dynamics, we define and review properties of ``fake'' stable and center-stable manifolds, following constructions from \cite{MR2630044}.  
\subsubsection{Notation}
Given $\mu$-a.e.\ $(\omega,x)\in X$, let 
$T_xM= \bigoplus_{i=1}^p E^i_\omega(x)$ denote the associated Oseledec splitting for the two-sided derivative cocycle on $(X,\mu)$.

Let $d=\dim M$ and write $\R^d= \bigoplus _{i = 1}^p \R_{i}$ where $\dim \R_{i} = m_i$,  the multiplicity of the $i$th Lyapunov exponent defined in \cref{multiplicity}.
Equip each $\R_{i}$ with the Euclidean inner product.  Given $v\in \R^d$, write $v= \sum v_i$ with $v_i\in \R_i$ and $\|v\|= \max\|v_i\|$.
Write $V_{\le i} = \bigoplus _{1\le j\le i} \R_{j}$ and $V_{> i}  
=\bigoplus _{ j> i} \R_{j}$.  

Write $\R^d(r)$ for the ball centered at $0$ of radius $r$ in $\R^d$.  
We fix $\eps>0$ sufficiently small (depending on the Lyapunov exponents of $\mu$) and set $\lambda_0= \max \{-\lambda_1, \lambda_p\}$.  

\subsubsection{Two-sided Lyapunov charts and globalization of local dynamics}\label{sec:charts}
The following summarizes properties of two-sided Lyapunov charts (see discussion in \cite[\S 4]{MR4504676} for modifications from the standard charts in the literature) as well as the standard trick of globalizing local dynamics in charts obtained by interpolating   via a bump function  between the local dynamics and the derivative (see, e.g., Prop.~7 and discussion in ~\cite[\S2]{MR730270}.)

\begin{proposition}\label{Lyapcharts} There exists a constant  $\hat \delta<1$ such that 
for every $0<\epsilon<\eps$, there exist  an  $F$-invariant set $\Lambda_0\subset X$ with $\mu(\Lambda_0)=1$, a measurable function $\ell\colon \Lambda_0\to [1,\infty)$,  and, for $(\omega,x)\in  \Lambda_0$ a $C^\infty $-map $\Phi_{(\omega,x)}\colon \R^d(1)\to M$ and a $C^{1+\beta}$ diffeomorphism $G_{(\omega,x)}\colon \R^d\to \R^d$ such that the following properties hold:    
\begin{enumerate}\item (Properties of charts)
\begin{enumerate}
\item $\Phi_{(\omega,x)}\colon \R^d(1)\to M$ is a diffeomorphism onto its image with $\Phi_{(\omega,x)}(0) = x$; 
\item $\|\Phi_{(\omega,x)}\|_{C^1}\le 1$ and $\|\Phi_{(\omega,x)}\inv \|_{C^1}\le  \ell(\omega,x) $;  
\item $D_0\Phi_{(\omega,x)}\R_i = E^i_{\omega}( x)$;
\item $\ell(F^n(\omega,x))\le   e^{|n|\epsilon}\ell(\omega,x)$ for $n\in \Z$.  
\end{enumerate}
\item (Properties of globalized dynamics in charts) 
\begin{enumerate}
\item $G_{ (\omega,x)}(0) = 0$.
\item For $v\in \R^d$ with $\|v\|<\hat \delta$, $$G_{(\omega,x)}(v) = \Phi_{F(\omega,x)}\inv \circ f_{\omega}\circ \Phi_{(\omega,x)}(v).$$
and $$G_{ F\inv(\omega,x)}\inv(v)  = \Phi_{F\inv(\omega,x)}\inv \circ f_{ \theta^{-1}(\omega)}\inv \circ \Phi_{(\omega,x)}(v).$$
\item For $v\in \R_i$, $e^{\lambda_i- \epsilon}\|v\|\le \| D_0G_{ (\omega,x)}v\|\le e^{\lambda_i+ \epsilon}\|v\|$.
\item $\|G_{(\omega,x)}- D_0G_{ (\omega,x)}\|_{C^1}\le \epsilon$ and  $\|G_{ F\inv(\omega,x)}\inv - D_0G_{ F\inv(\omega,x)}\inv \|_{C^1}\le \epsilon$.
\item $\Hol(DG_{(\omega,x)})<C$ and $ \Hol(DG_{F\inv(\omega,x)}\inv)\le C$ for some uniform $C$ independent of $\omega$, $x$, and $0<\epsilon\le \eps$.
\end{enumerate}
\item (Measurability) $(\omega,x)\mapsto \Phi_{(\omega,x)}$,  $(\omega,x)\mapsto G_{(\omega,x)}$,  $(\omega,x)\mapsto G_{F\inv(\omega,x)}\inv$  are measurable. 
\end{enumerate}
\end{proposition}
Given $0<\epsilon<\eps$ and corresponding function $\ell\colon \Lambda_0\to [1,\infty)$, we refer to the collection $\{\Phi_{(\omega,x)}\}$ 
as a family of \emph{$(\epsilon, \ell)$-Lyapunov charts.}

\subsubsection{``Fake'' pseudo-stable manifolds} 
Fix a choice of $(\epsilon, \ell)$-Lyapunov charts as in \cref{Lyapcharts}.  Recall that the globalized dynamics of $\{G_{(\omega,x)}\}\colon \R^d\to \R^d$ described in \cref{Lyapcharts} is not canonical and depends on the choices (of bump function and interpolation to derivative) in the construction.  However, dynamical 
objects built from the (non-canonical) globalized dynamics are, at times, useful to understand objects (at points whose future orbit remains in charts) defined by the true dynamics $F$.  

We follow (the proof of)  \cite[Prop.~3.1]{MR2630044} and construct ``fake'' stable and center-stable manifolds for our choice of globalized dynamics.  The existence of such objects follows from standard  (pseudo-)stable manifolds theorems using  either the Perron--Irwin method or Hadamard graph transform method.  (See discussion in  step 1 of the proof of \cite[Prop.~3.1]{MR2630044}, which invokes  graph-transform methods such as those from \cite{MR0501173}. For discussion on regularity of the manifold, see discussion in \cite[\S6.3]{MR0983869}.)

Let $G_{(\omega,x)}^{(0)}$ be the identity map on $\R^d$ and for $n\ge 1$, write $$G_{(\omega,x)}^{(n)}:= G_{F^{n-1}(\omega,x)}\circ \cdots \circ G_{F(\omega,x)}\circ G_{(\omega,x)}.$$
\begin{lemma}[``Fake'' pseudo-stable manifolds]\label{lem:localunstable}
For every $1\le i\le p$,  every  $(\omega,x)\in \Lambda_0$,  and every $v\in \R^d$  there is a  $C^{1+\text{H\"older}}$ function $$h^i_{(\omega,x),v} \colon V_{\le i}  \to V_{>i}$$
such that  the following hold:
\begin{enumerate}
\item There is a constant $c>0$ (independent  of $(\omega,x)$ and $v$) such that  
\begin{itemize}
\item $\Hol(Dh^i_{(\omega,x),v}) <c$  if  $\lambda_i<0$,  and 
\item $\Hol[\beta'](Dh^i_{(\omega,x),v}) <c$  for some $\beta'>0$ if $\lambda_i=0$.
\end{itemize}
\item $\|Dh^i_{(\omega,x),v}\|\le 1/3$ and $h^i_{(\omega,x),v}$ is bounded if $\lambda_{i+1}>0$.
\end{enumerate}
Write $  W^i_{(\omega,x)}(v)$ for the graph of $ h^i_{(\omega,x),v} \colon V_{\le i}  \to V_{>i}.$
\begin{enumerate}[resume]
\item If $u\in   W^i_{(\omega,x)}(v)$ then $  W^i_{(\omega,x)}(v)=   W^i_{(\omega,x)}(u)$.
\item For $u\in   W^i_{(\omega,x)}(v)$, 
	$$e^{\lambda_1 - 2 \epsilon} \|v-u\| \le  \|G_{ (\omega,x)} (v) -G_{ (\omega,x)}(u)\| \le e^{\lambda_i+ 2\epsilon} \|v-u\|.$$
\item $   W^i_{(\omega,x)}(v)$ coincides with the set of all $u\in \R^d$ for which  
$$ \limsup_{n\to \infty}\frac 1 n \log \|G_{  (\omega,x)}^{(n)}(u) - G_{  (\omega,x)}^{(n)}(v)\|\le \lambda_i + 2\epsilon.$$
\item $G_{ (\omega,x)}(W^i_{(\omega,x)}(v)) = W^i_{F(\omega,x)}(G_{ (\omega,x)}(v))$.
\end{enumerate}
Moreover, 
\begin{enumerate}[resume]
\item The family $\{h^i_{(\omega,x),v}\}$ depends continuously on $v\in \R^d$ and measurably on $(\omega,x)\in   X$.  
\item For $1\le j<i\le \min\{s+1, p\}$ we have $   W^j_{(\omega,x)}(v)\subset    W^i_{(\omega,x)}(v)$; moreover 
\label{good holonomies} for $1\le i\le \min\{s+1, p\}$, the foliation  of $   W^i_{(\omega,x)}(v)$ by  submanifolds of the form $$\{   W^{i-1}_{(\omega,x)}( y) : y\in    W^i_{(\omega,x)}(v)\}$$ has uniformly Lipschitz holonomies.
\end{enumerate}
Finally, when $v=0$, 
\begin{enumerate}[resume]
\item $h^i_{(\omega,x),0}(0) = 0$ and $D_0h^i_{(\omega,x),0}(0) = 0$.
\end{enumerate}
\end{lemma}
For proof and explanation of the Lipschitzness of holonomies in conclusion \eqref{good holonomies}, see \cite[Appendix]{MR1709302} when $i\le s$ and  \cite{MR4504676} for the general case (especially when $i=s+1$).

In all further uses, we will only apply \cref{lem:localunstable} for $1\le i\le s$ (when $\lambda_{s+1} >0$) or for $1\le i\le s+1$ (when $\lambda_{s+1} =0$).  In the fully contracting setting where $s=p$, declare  $$W^{s+1}_{(\omega,x)}(v)= W^s_{(\omega,x)}(v)= W^p_{(\omega,x)}(v)=\R^d.$$  
For $1\le i\le s+1$ and for $(\omega,x)\in X$, write 
 \begin{equation}\label{eq:localfake}
  W^i_{(\omega,x), \delta}(v):= \{u\in   W^i_{(\omega,x)}(v): \|u\|< \delta\}  \end{equation}
and 
\begin{equation}
W^i_{ \loc,\delta}(\omega,x):=  \Phi_{ (\omega,x)}(  W^i_{(\omega,x),  \delta}(0)).\end{equation}
Also write $$W^0(\omega,x)=W^0_{ \loc,\delta}(\omega,x)= x.$$  

For $1\le i\le s$, we note that choice of local $i$th fast stable manifold $W^i_{ \loc,\delta}(\omega, x)$ depends on both the future and past itinerary  $\omega$ as it is defined relative to  two-sided charts.  
In particular, if $\tau(\omega) = \tau(\omega')$ then,  while we have $W^i_\omega(x) = W^i_{\omega'}(x)$, in general we would expect
$$W^i_{ \loc,\delta}(\omega,x)\neq W^i_{ \loc,\delta}(\omega',x).$$
To account for our definition of local stable manifolds dependent on both future and past itineraries, for $0\le i\le s$, we also write 
\begin{equation}\wtd W^i_{ \loc,\delta}(\omega,x):= \tau\inv \bigl(\{\tau( \omega)\}\times W^i_{ \loc,\delta}(\omega,x)\bigr)\end{equation}
for the saturation of $\{\omega\}\times W^i_{ \loc,\delta}(\omega, x)$ by atoms of the $\sigma$-algebra $\wtd \sigA$.

\subsection{Center-stable sets and properties of   stable manifolds}

Given $0< \delta\le  1$ and $(\omega,x)\in \Lambda_0$, define the associated \emph{center-stable sets}:  
 $$S^{cs}_{(\omega,x),\delta} = \{y\in \R^d(\delta): \| G_{(\omega,x)}^{(n)}(y) \| <\delta : \text{ for all } n\ge 0\}.$$
We note by definition that \begin{equation} \label{eqiv of cs} G_{(\omega,x)} (S^{cs}_{(\omega,x),\delta})\subset S^{cs}_{F(\omega,x),\delta}.\end{equation}

Recall $\lambda_0= \max \{-\lambda_1, \lambda_p\}=\max\{|\lambda_i|, 1\le i\le p\}$.  Fix $0<\epsilon<\eps$ and \begin{equation}\label{delta0}\delta_0:= \frac{1}{4e^{\lambda_0 +2\epsilon}}\min\{\hat \delta, \kappa\}\end{equation}
where $\hat \delta>0$ is as given by \cref{Lyapcharts} and $\kappa$ is as in the construction of $\hat \xi^i$ in \cref{sec:sss}.  
 
\begin{claim}\label{claim:neatfacts}
For   $0<\delta\le \delta_0$ and every $1\le i\le s$, using the notation of \cref{eq:localfake} the following hold for $(\omega,x)\in \Lambda_0$:
\begin{enumerate}[label= {(\alph*)}, ref= {(\alph*)}]
\item \label{yuaa} $S^{cs}_{(\omega, x),2\delta}= \{ v\in \R^d(2\delta):\|\Phi_{F^n(\omega,x)}\inv \circ f_\omega^{(n)}\circ \Phi_{(\omega,x)}(v)\|< 2\delta 
 \text{ for all $n\ge0 $}\}$; 
\item  $W^i_{\loc,2\delta}(\omega, x)=\Phi_{(\omega,x)}(   W^i_{(\omega,x),  2\delta}(0))$ is an open neighborhood of $x$ in $W^{i}_\omega(x)$;\item \label{yucc} $S^{cs}_{(\omega, x), 2\delta}\subset   W^s_{(\omega,x),  2\delta}(0)$ if $\lambda_{s+1}>0$ and 
$S^{cs}_{(\omega, x), 2\delta}\subset   W^{s+1}_{(\omega,x),  2\delta}(0)$ if $\lambda_{s+1}=0$.  
\end{enumerate}
Moreover, for $z\in S^{cs}_{(\omega,x),\delta}$    and $1\le  i \le s$, 
\begin{enumerate}[resume,label= {(\alph*)}, ref= {(\alph*)}] 
\item  \label{yudd} $  W^i_{(\omega,x),  2\delta} (z)  \subset  S^{cs}_{(\omega, x), 2\delta}$;

\item \label{yuee} $   W^i_{(\omega,x),  2\delta} (z)
= G_{ (\omega,x)}\inv (  W^i_{F(\omega,x),  2\delta} (G_{(\omega,x)}(z)))\cap \R^d (2\delta)  $;
\item  \label{yuff} if $(\omega', \Phi_{(\omega,x)}( z)) \in  \Lambda_0$ for some $\omega'$ with $\tau(\omega) = \tau(\omega')$ then 
 $$\Phi_{(\omega,x)} \inv\left( W^i_{\omega'}( \Phi_{(\omega,x)}( z)) \right) \cap S^{cs}_{(\omega, x), 2\delta} \subset    W^i_{(\omega,x), 2\delta}(z).$$ 
\end{enumerate}
\end{claim}

\begin{proof}
Conclusions \ref{yuaa}--\ref{yucc} follow using that the globalized dynamics and dynamics in charts coincide on $\R^d(2\delta)$ and applying dynamics inside charts.  
Conclusions \ref{yudd} and \ref{yuee} follow using that $G_{(\omega, x)}$ is an $\epsilon$-small $C^1$ perturbation of $D_0 G_{(\omega, x)}$ for a sufficiently small $\epsilon>0$.  A standard computation then shows  for 
$z\in \R^d(\delta)$ with $G_{(\omega,x)} (z)\in \R^d(\delta)$, 
$$ G_{(\omega,x)} (  W^i_{(\omega,x),2\delta}(z))\subset W^i_{F(\omega,x),2\delta}(G_{(\omega,x)}(z)).$$

For conclusion \ref{yuff}, let $z_0 = \Phi_{(\omega,x)}( z)\in M$  
and consider  
 $$y\in W^i_{\omega'}( z_0)  \cap \Phi_{(\omega,x)}(S^{cs}_{(\omega, x), 2\delta}) =W^i_{\omega}( z_0)  \cap \Phi_{(\omega,x)}(S^{cs}_{(\omega, x), 2\delta}) .$$ 
We have  $$\limsup_{n\to \infty}  \frac 1 n \log d_M\left ( f_{ \omega}^{(n)}(z_0),f_{ \omega}^{(n)}(y)\right)=\limsup_{n\to \infty}  \frac 1 n \log d_M\left ( f_{ \omega'}^{(n)}(z_0),f_{ \omega'}^{(n)}(y)\right) \le \lambda_i$$
and, since the Lipschitz distortion by the charts   $\|\Phi_{F^n(\omega,x)}\inv \|_{C^1}$ increases  slow-exponentially, 
$$ \limsup_{n\to \infty}\frac 1 n \log \|G_{  (\omega,x)}^{(n)}(z) - G_{  (\omega,x)}^{(n)}( \Phi_{(\omega,x)}\inv( y))\|\le \lambda_i + \epsilon.$$
Since $z\in S^{cs}_{(\omega, x), \delta} $,  there is some $n\ge 0$ such that 
 $$G_{(\omega,x)}^{(n)} (\Phi_{(\omega,x)} \inv(y))\in W^i_{F^n(\omega,x), 2\delta}(G_{(\omega,x)}^{(n)}(z)).$$ 
Let $n_0$ be the minimal such $n$.  If $n_0>0$, then \ref{yuee} implies $G_{(\omega,x)}^{(n_0-1)} (\Phi_{(\omega,x)} \inv(y))\notin \R^d(2\delta)$, contradicting the choice of $y\in \Phi_{(\omega,x)} (S^{cs}_{(\omega, x), 2\delta})$.  
\end{proof}

\subsection{Partitions adapted to   Lyapunov charts}

Recall our family of two-sided $(\epsilon, \ell)$-Lyapunov charts $\{\Phi_{(\omega,x)}\}$ in \cref{Lyapcharts}.  Fix $0< \delta\le \delta_0$ (where $\delta_0$ is defined in \cref{delta0}). Recall $\pi_2\colon X=\Omega\times M\to M$ is the canonical projection.

\begin{definition}\label{def:adapted}We say a measurable partition $\partP$ of $(X,\mu)$ is \emph{adapted to $(\{\Phi_{(\omega,x)}\}, \delta)$} if, for  almost every ${(\omega,x)}\in \Lambda_0$, 
 $$\pi_2 (\partP^-\vee\Xi(\sigC)){(\omega,x)} \subset \Phi_{(\omega,x)}( S^{cs}_{\delta, (\omega,x)}).$$
 \end{definition}

 We note that $ S^{cs}_{\delta, (\omega,x)}$ is defined for $(\omega,x)\in X$, rather than for $(\omega,x)\in \what X$.  
Write \begin{equation}
\wtd S^{cs}_{\delta, (\omega,x)}:= \tau\inv (\{\tau(\omega)\}\times  \Phi_{(\omega,x)}( S^{cs}_{\delta, (\omega,x)}))\end{equation}
for the $\wtd \sigA$-saturation of the image of $S^{cs}_{\delta, (\omega,x)}$ under $\Phi_{(\omega,x)}$.  By \eqref{eqiv of cs} and since $\wtd \sigA$ is sub-invariant, we have  $$F(\wtd S^{cs}_{\delta, (\omega,x)}) \subset \wtd S^{cs}_{\delta, F(\omega,x)}.$$
Then $\partP$ is adapted to $(\{\Phi_{(\omega,x)}\}, \delta)$ if and only if 
$$(\partP^-\vee\Xi(\sigC)){(\omega,x)}\subset \wtd S^{cs}_{\delta, (\omega,x)}.$$

 \begin{lemma}\label{lem:subpartfinent}
 For every   $0<\delta\le \delta_0$, there exists a   measurable partition $\partP$ of $(X,\mu)$ adapted to $(\{\Phi_{(\omega,x)}\}, \delta)$ with $H_\mu(\partP)<\infty.$
 \end{lemma}
The proof follows {\it mutatis mutandis}  as in   \cite[Lem.\  2.4.2]{MR819556}.  
We provide an outline as  objects appearing in the proof will be needed in the proof of \cref{thm:exactdimproj}.
\begin{proof}[Proof sketch of \cref{lem:subpartfinent}]
Fix $\ell_0$ sufficiently large so that the set $\Lambda:=\{(\omega,x)| \ell(\omega,x) \le \ell_0\}$ has positive $\mu$-measure. 
For $(\omega,x)\in X$, let 
\[\begin{aligned}
T_+(\omega,x)&=\inf \{m\ge1:F^m(\omega,x)\in\Lambda\},\\
T_-(\omega,x)&=\inf\{m\ge1:F^{-m}(\omega,x)\in\Lambda\},\\
T(\omega,x)&=\max\{T_+(\omega,x),T_-(\omega,x)\}.
\end{aligned}\]
Let $$r(\omega,x) = \begin{cases}
\frac {\delta}{\ell_0} e^{-T(\omega,x)(\lambda_0 +3\epsilon)} &(\omega,x)\in \Lambda\\
1 &(\omega,x)\notin \Lambda.
\end{cases} $$
We have $-\log r(\omega,x) \in L^1(\mu)$ and so may apply  Ma\~n\'e's lemma which we summarize in 
  \cref{lemma:MANE} below.
  Let $\calP$ be the Ma\~n\'e's  partition built in \cref{lemma:MANE}.  Then  $H_\mu(\partP)<\infty.$

  From the construction of  $\calP$ in \cref{lemma:MANE}, $\calP$ has the following property:  if $n\ge T_+(\omega,x) +1$, if $(\omega',y)\in \partP^{[0,n)}(\omega,x)$, and if \begin{equation}\label{eq:pesinadapt} f_{\theta^k(\omega)} = f_{\theta^k(\omega') } \text{ for all $0\le k\le n-1$},\end{equation}
  then 
  $$\pi_2\bigl(\partP^{[0,n)}(\omega,x) \cap \sigF(\omega',y)\bigr)\subset \Phi _{(\omega,x)}(\R^d(\delta)).$$
Indeed, we show something stronger:
Let $0=r_0$ and let $1\le  r_1<r_2<\dots <n$ denote the distinct times  $r_j\ge 1$ with $F^{r_j}(\omega,x)\in \Lambda$.  For $j\ge 1$, 
$$\partP(F^{r_j}(\omega,x))\subset \Omega \times \Phi _{F^{r_j}(\omega,x)}\bigl(\R^d(\delta e^{-T(F^{r_j}(\omega,x))(\lambda_0 +3\epsilon)})\bigr).$$
Thus for $j\ge 1$, if $\omega'$ satisfies \eqref{eq:pesinadapt} then from the dynamics inside Lyapunov charts (along the backwards orbit from $F^{r_{j}}(\omega,x)$ to $F^{r_{j-1}}(\omega,x)$), for $0\le k\le  r_j -r_{j-1}$, we have
$$F^{-k}\bigl(\partP (F^{r_j}(\omega,x)  ) \cap \sigF (F^{r_j}(\omega',y) )\bigr)\subset \Omega \times \Phi _{F^{r_{j}-k}(\omega,x)}(\R^d(\delta)).$$

  We remark that \eqref{eq:pesinadapt} holds for all $(\omega',y)\in \sigC(\omega,x)$ and all $n\ge 0$.  In particular, this establishes that $\partP$ is adapted to $(\{\Phi_{(\omega,x)}\}, \delta)$.
\end{proof}

The proof of \cref{lem:subpartfinent} only needs the following adaptation of  a lemma of Ma\~n\'e from \cite{MR627789}.

\begin{claim}[{\cite[Lem.\  2]{MR627789}}]\label{lemma:MANE}
Let $r\colon (X,\mu)\to (0,1]$ be a measurable function such that $$\int -\log (r(\omega,x))\, d \mu(\omega,x)<\infty.$$
Then there exists a measurable partition $\calP$ of $(X,\mu)$ with $H_\mu(\calP)<\infty$ such that for $\mu$-a.e.\ $(\omega,x)$,
$$\diam(\pi_2(\partP(\omega,x)))\le r(\omega, x).$$
\end{claim}
\begin{proof}[Proof sketch of \cref{lemma:MANE}]
For each $n\ge 0$, pick a finite partition  $\partQ_n$ of $M$ with $\diam (\partQ_n(x) )\le e^{-n}$ for every $x\in M$ and $|\partQ_n|\le C e^{d n}$ for some uniform constant $C$.  
Take the countable partition of $X$ by sets of the form  $U_n:= \{(\omega,x)\in X: -(n+1) <\log r(\omega,x)\le -n\}$.  Define $\partP$ as follows: if $(\omega,x)\in U_n$ set $$\partP(\omega,x)= U_n\cap \left(\Omega \times \partQ_{n+1}(x)\right).$$
Since $\log (\mathrm{card}(\restrict{\partP}{U_n}))\le \log C+ d(n+1)$ and since 
$$\sum (n+1) \mu(U_n)\le \int -\log (r(\omega,x))\, d \mu(\omega,x) +1<\infty,$$ then,
exactly as in the proof of \cite[Lem.\  2]{MR627789}, we have $H_\mu(\partP)<\infty$.  
\end{proof}

As an immediate corollary of    \cref{claim:neatfacts}\ref{yuff} we have the following: 
\begin{claim}\label{claim:insides}
 Let $0<\delta\le\delta_0$ 
 and let $\partP$ be adapted to $(\{\Phi_{\omega,x} \}, \delta)$.  
  Then for $1\le i\le s$ and  $  \mu$-almost every $(\omega,x)\in X$, 
 	$$(\partP^-\vee \xi^i)(\omega,x)\subset  \wtd  W^i_{\loc,2\delta}(\omega,x).$$ 
 \end{claim}

 \subsection{Properties of partitions adapted to Lyapunov charts}\label{little props about eta}
Let   $\partP$ be a measurable partition of $(X,\mu)$ adapted to $(\{\Phi_{(\omega,x)}\}, \delta_0)$ with $H_\mu(\calP)<\infty$.  
Let $\wtd S= \tau\inv S_\rho$ where $S_\rho$ was chosen in \cref{bighatballs} in the construction of the partitions $\hat \xi^i$ in \cref{sec:sss}.  Let $E=  \{  X\sm \wtd S, \wtd S\}$.  

Recall we have a sub-invariant measurable partition $\zeta$ with  $\Xi(\sigC)\prec \zeta$.  
 Take $\partQ$ an arbitrary  partition with $H_\mu(\partQ)<\infty$ and take $\eta_*= (\partP\vee \zeta\vee E\vee \partQ)^-$.  For $0\le i\le s$ take $$\eta_i = \eta_*\vee \xi^i.$$  
For notational conventions later on, we also write 
 $\eta_{s+1}:= \eta_*$.  Note that if  no   exponent is zero, it follows exactly as in   \cite{MR743818} (c.f \cite[Prop.\ (4.1)]{MR743818} and  \cite[Prop.\ (4.5)]{MR743818}) that $\eta_*=\eta_{s+1}= \eta_s.$
 In the presence of zero exponents we may have  $ \eta_s \neq \eta_{s+1}.$  From Claim \ref{claim:neatfacts}\ref{yucc} and Claim \ref{claim:insides}, for almost every $(\omega,x)$ we have
 \begin{enumerate}
\item   $\eta_{s+1}{(\omega,x)} \subset \wtd S^{cs}_{(\omega,x),\delta_0}\subset   \wtd W_{\loc,2\delta_0}^{s+1}(\omega,x)$;
\item    $\eta_{i}({\omega,x}) \subset  \wtd  W^i_{\loc,2\delta_0} (\omega,x)$  for $1\le i\le s$.
\end{enumerate}

 As in \cite[Lems.\  3.3.1--3.3.2]{MR819556} and \cite[Lems.\  11.1.2--11.1.3]{MR819557},  we have  the following.  
 
 \begin{lemma}\label{lemma:44} 
For each $0 \le i\le s$ and for $\mu$-almost every ${(\omega,x)}$, 
\begin{enumerate}[label= {(\alph*)}, ref= {(\alph*)}]
\item \label{PP1} $\wtd  W^i_{\loc,2\delta_0} (\omega,x)\cap \eta_{i+1}{(\omega,x)} = \eta_i {(\omega,x)}$.    
\item \label{PP2}  $(F\eta_i) (\omega,x)= \eta_{i} (\omega, x) \cap( F\eta_{i+1})(\omega, x)  $. 
\end{enumerate}
 \end{lemma} 
 
 \begin{proof}
The result is clear when $i=0$.  Suppose $1\le i\le s$.  
For \ref{PP1}, first consider 
$$(\omega',y) \in  \wtd  W^i_{\loc,2\delta_0}(\omega,x)\cap \eta_{i+1}{(\omega,x)}.$$  
Note that $\tau(\omega) =\tau(\omega')$.  
We claim $(\omega',y)\in  \xi^i(\omega,x)$.  
Indeed,  $\eta_*$ refines  $\{\wtd S, X\sm \wtd  S\}$; if  $(\omega',y)\in \eta_*(\omega,x)$, it follows  for all $m\ge 0$ that  $F^m(\omega',y)\in \wtd  S$ if and only if $F^{m}(\omega, x)\in \wtd  S$.  
 Moreover, if $
 (\omega',y) \in \wtd  W^i_{\loc,2\delta_0} (\omega,x)\cap \eta_{i+1}{(\omega,x)}$ then $y\in W^i_{\omega}(x)$; 
  using the dynamics in $W^i_{\omega}(x)$ relative to charts along the orbit of 
   $(\omega,x)$ and using that Lyapunov charts are $1$-Lipschitz (even when restricted to the intrinsic metric in stable manifolds $W^i_{F^n(\omega,x),2\delta_0}(0)$ in charts),   
    we have $$d^s(f^{(m)}_{\omega'}(y), f^{(m)}_{\omega}(x))=d^s(f^{(m)}_{\omega}(y), f^{(m)}_{\omega}(x))\le 2\delta_0 e^{m(\lambda_i+2\epsilon)} \le { \kappa}$$ which, by \cref{in same tilde chart} and our choice of ${ \delta_0}\le { \kappa/4 }$, implies  $\hat \xi ^i(\tau (\omega'),y)= \hat  \xi^i(\tau (\omega),x)$   and so 
    $ \xi ^i(\omega',y)=   \xi^i(\omega,x)$.
    The  reverse inclusion is  Claim \ref{claim:insides}.

For \ref{PP2},   
  \cref{claim:neatfacts}\ref{yuee}, ensures
\begin{align*}  
 F\inv( \wtd  W^i_{\loc,2\delta_0} (\omega,x))\cap  \eta_{i+1}(F\inv(\omega,x)) 
 &\subset  \wtd  W^i_{\loc,2\delta_0} (F\inv( \omega,x)).
\end{align*}
Using part \ref{PP1}  (at $F\inv( \omega,x)$) we  have 
\begin{align*}
 \eta_{i}(\omega,x) \cap (F\eta_{i+1})(\omega,x) 
&= 
 \left( \wtd  W^i_{\loc,2\delta_0} (\omega,x) \cap \eta_{i+1}(\omega,x)\right) \cap (F\eta_{i+1})(\omega,x)  \\
&= 
 \wtd  W^i_{\loc,2\delta_0} (\omega,x) \cap (F\eta_{i+1})(\omega,x)  \\
&= F\left( F\inv \bigl( \wtd  W^i_{\loc,2\delta_0} (\omega,x)\bigl) \cap \eta_{i+1}\bigl( F\inv (\omega,x) \bigr) \right)\\
&\subset  F\left( \wtd  W^i_{\loc,2\delta_0} (F\inv(\omega,x)) \cap \eta_{i+1}\bigl( F\inv (\omega,x) \bigr) \right)\\
&= F\left (\eta_{i}\bigl( F\inv (\omega,x) \bigr)  \right)\\
&=(F\eta_i)(\omega,x).\end{align*}
The reverse inclusion holds trivially.
\end{proof}

\subsection{One further partition}

Recall our family of forwards $(\epsilon, \ell)$-Lyapunov $\{\Phi_{(\omega,x)}\}$ relative to which  
 local fast stable manifolds $W^i_{\loc,\delta}(\omega,x)\subset M$ in $M$ are parametrized by local graphs $W^i _{(\omega, x), \delta}(0)$ of functions in $\R^d$.  

We   have the following elementary fact.  We note that we consider atoms of the partition $ (\partP \vee   \xi^i)^{(-n,0]}$ rather than $(\partP ^{(-n,0]})\vee   \xi^i$ and so the proof is a simpler version of \cite[Lems.\  9.3.3, 9.3.2]{MR819557}  (again using \cref{lemma:MANE}). 

 \begin{lemma}\label{lem:partadaptedtolocalentropy}
There exists a partition $\partP$ of $(  X,  \mu)$ with $H_{  \mu}(\partP)<\infty$ and a measurable $n_0\colon X\to \N \cup\{\infty\}$ which is finite a.s.\ such  that for every $1\le i\le s$, $\mu$-almost every $(\omega,x)\in   X$,   and  every $n\ge n_0(\omega, x)$,
$$ \left ((\partP \vee   \xi^i)^{(-n,0]}(\omega, x)\right) \subset \wtd W^i_{ \loc, 
e^{n(\lambda_i + 2 \epsilon)}}(\omega,x).$$ 
\end{lemma}
\begin{proof}
Fix $0<\delta\le \delta_0<1$.  
 Recall the atoms of $\hat \xi^i$ are precompact in $W^{i}(\omega,x)\subset \{\omega\}\times M$ for $(\omega,x)\in \what X$.  Thus, there is a measurable function $c\colon (X,\mu)\to (0,\infty)$ such that  for $1\le i\le s$, 
$$\left\{ (\omega',y)\in \xi^i(\omega, x) : d_M(x,y)\le c(\omega,x)\right\} \subset \wtd W^i_{\loc,\delta}(\omega,x).$$
Fix a set $\Lambda\subset X$ with $\mu(\Lambda)>0$ on which the measurable function  $\ell\colon X\to (0,\infty)$ is  bounded from above by $1\le \ell_0<\infty$ and the measurable function $c\colon X\to (0,\infty)$ is bounded from below by $0<c_0\le 1$.

As in the proof of \cref{lem:subpartfinent}, for $(\omega,x)\in X$, let
\[\begin{aligned}
T_+(\omega,x)&=\inf\{m\ge1:F^m(\omega,x)\in\Lambda\},\\
T_-(\omega,x)&=\inf\{m\ge1:F^{-m}(\omega,x)\in\Lambda\},\\
T(\omega,x)&=\max\{T_+(\omega,x),T_-(\omega,x)\}.
\end{aligned}\]
Set $q= \min\{c_0,\delta/\ell_0\}$ and define
\[r(\omega,x)=
\begin{cases}
q  e^{-T(\omega,x)(\lambda_0+3\epsilon)}&(\omega,x)\in\Lambda,\\
1 &(\omega,x)\notin\Lambda.
\end{cases}\]
Then $-\log r$ is in $L^1(\mu)$.  Let $\partP$ be the measurable partition obtained in  \cref{lemma:MANE}.  

Let $n_0(\omega,x) = T_-(\omega,x)+1$.  For $n\ge n_0(\omega,x)$, let $k_0$ be the maximal $k\le (n-1)$ with $F^{-k}(\omega,x)\in \Lambda$. 
Since $T_-(F^{-k_0}(\omega,x) )\ge n-k_0$ and since  $\Phi _{F^{-k_0}(\omega,x)}\inv$ is $\ell_0$-Lipschitz,
we have
\begin{align*}\pi_2(\partP(F^{-k_0}(\omega,x)))&\subset \Phi _{F^{-k_0}(\omega,x)}(\R^d(\ell_0q e^{-T(F^{-k_0}(\omega,x))(\lambda_0 +3\epsilon)}))\\ &\subset \Phi _{F^{-k_0}(\omega,x)}(\R^d(\delta  e^{(n-k_0)(\lambda_1 -3\epsilon)})).\end{align*}
Also, since  $q\le c_0$, 
   \[   (\partP\vee \xi^i)(F^{-k_0}(\omega,x))\subset  \wtd W^i_{\loc,\delta}(F^{-k_0}(\omega,x)).   \]
Take $(\omega',y)\in( \partP\vee \xi^i)(F^{-k_0}(\omega,x))$.  
Then 
\begin{align*}\Phi_{F^{-k_0}(\omega,x)}\inv(y)& \in 
W^i_{F^{-k_0}(\omega,x), \delta }(0)\cap \R^d\left(\delta e^{(n-k_0)(\lambda_1 -3\epsilon)})\right)
\\& \subset W^i_{F^{-k_0}(\omega,x), \delta e^{(n-k_0)(\lambda_i +2\epsilon)}}(0).\end{align*}
Let $\omega_{-k_0} = \theta^{-k_0}(\omega)$.  
By the local graph invariance and contracting dynamics along  fake stable manifolds relative to charts, we have
 $$\Phi_{(\omega,x)}\inv(f^{(k_0)}_{\omega_{-k_0}}(y))\in 
W^i_{(\omega,x), \delta e^{n(\lambda_i +2\epsilon)}}(0).$$
and so using that $\Phi_{(\omega,x)}$ is $1$-Lipschitz, $$\left(\omega,f_{\omega_{-k_0}}^{(k_0)}(y)\right)\in  \wtd W^i_{ \loc, 
e^{n(\lambda_i + 2 \epsilon)}}(\omega,x).$$ 
Since $f_{\omega_{-k_0}}^{(j)} = f_{\omega'}^{(j)}$ for every $(\omega',y)\in \xi^i(F^{-k_0}(\omega,x))$ and $j\ge 0$,
we also have $$F^{k_0}(\omega',y)\in \wtd W^i_{ \loc, 
e^{n(\lambda_i + 2 \epsilon)}}(\omega,x)$$ for every $(\omega',y)\in (\partP\vee \xi^i)(F^{-k_0}(\omega,x)).$
Since $k_0\le (n-1)$, we have 
\[(\partP\vee \xi^i)^{(-n,0]}(\omega,x)\subset F^{k_0}(\partP\vee \xi^i)(F^{-k_0}(\omega,x))\subset \wtd W^i_{ \loc, e^{n(\lambda_i + 2 \epsilon)}}(\omega,x). \qedhere\]
\end{proof}

\section{Ledrappier--Young, recapitulated}\label{LYagain}
To establish all  announced results, it remains to establish \cref{thm:realthm,claim:bound on transverse dimension,bdd on dimn 2,thm:exactdimproj,prop:entropyisunstable}.
These results \cref{thm:realthm,claim:bound on transverse dimension,bdd on dimn 2,prop:entropyisunstable} follow  immediately from the (in)equalities in \cref{sec:outline} below. \cref{thm:exactdimproj} follows from essentially the same proof as  inequality \ref{ineqI}; we formulate this in \cref{sss:proj msr} below.  
This section re-runs the main arguments from the Ledrappier--Young papers \cite{MR819556,MR819557}, establishing these inequalities.

\subsection{Key (in)equalities}\label{sec:outline}
Recall the measurable partitions $\hat \xi^i$ of $(\what X, \hat \mu)$ described in \cref{sss:1}.
Recall we  write $\hat \xi^0=\ptpart$ for the point partition on $\what X$ and for $1\le i\le s$, write $ \xi^i$ for the measurable partition of  $(  X,   \mu)$ given by $\xi^i := \tau\inv \hat \xi^i$ where $\tau\colon X\to \what X$ is the factor map.  

Fix $\zeta $ a sub-invariant measurable partition with $ \Xi(\sigC)\prec \zeta$.  
For every  $0\le  i \le s$, let \begin{equation}\label{eqhi}h_i:= h_\mu(F, \xi^i\vee \zeta).\end{equation}
Recall (by \cref{lem:sameentropy}) that the value of $h_i$ is independent of any choices made in constructing $\xi^i$.  
Also note by definition that $h_0 = h_{\mu}(F, \xi^0\vee \zeta) $.

Under the assumption that $ h_\mu(F, \zeta)<\infty$ we have $ \maxent_\mu(F\mid \sigma(\zeta))<\infty$ (and thus $h_i<\infty$ for every $i$) as discussed in \cref{finiteentass}.
For $1\le  i \le s$, we show  the following (in)equalities hold:
\begin{enumerate}[{label*={(\Roman*)}, ref= {(\Roman*)} }]
\item \label{ineqI} $  h_i -   h_{i-1}\ge -\lambda_i  \big( \overline{\dim}^i(   \mu\mid \zeta) - \overline{\dim}^{i-1}(    \mu\mid \zeta)\big) $
\item  \label{ineqII}$  h_i -   h_{i-1}\le- \lambda_i  \big( \underline{\dim}^i(    \mu\mid \zeta) - \underline{\dim}^{i-1}(    \mu\mid \zeta)\big) \le -\lambda_i  m_i $ 
\item  \label{ineqIIII}$   h_s =     \maxent_\mu(F\mid \sigma(\zeta))$.  
\end{enumerate}
Moreover, if $\zeta'$ is sub-invariant and if $\zeta\prec \zeta'$ then 
\begin{enumerate}[{resume,label*={(\Roman*)}, ref= {(\Roman*)} }]
\item  \label{ineqV} $\overline{\dim}^i(    \mu\mid \zeta') - \overline{\dim}^{i-1}(    \mu\mid \zeta')\le \underline{\dim}^i(    \mu\mid \zeta) - \underline{\dim}^{i-1}(    \mu\mid \zeta).$
\end{enumerate}

\cref{thm:realthm}, \cref{bdd on dimn 2}, \cref{claim:bound on transverse dimension}, and  \cref{prop:entropyisunstable} then follow immediately from \ref{ineqI}--\ref{ineqV}.

\subsection{\texorpdfstring{Proof of inequality \ref{ineqI} and \cref{thm:exactdimproj}}{Proof of inequality (I) and Theorem 4.20}}
Recall the map $\tau\colon (X,\mu) \to (\what X, \hat \mu)$, $\tau(\omega, x) = (\tau(\omega), x).$
We write \begin{equation}\wtd B((\omega , x), r):= \tau\inv (\{\tau(\omega) \}\times B(x,r))\end{equation} for the $\sigA$-saturation of the ball of radius $r>0$ centered at $x$ in $\{\omega \}\times M\subset \what X$.  
We also  write $\wtd B^i((\omega , x), r) = \wtd B((\omega , x), r)\cap 
 \wtd  W^i_{\loc,2\delta}(\omega,x).$

\begin{lemma}\label{lem:ooo}
For each $1\le  i \le s$ and every $\epsilon>0$ sufficiently small, 
there exists a partition $\partP$ of $  X$ with $H_\mu(\partP)<\infty$ and 
a measurable function $n_0\colon X\to \N$ such that for $  \mu$-a.e.\  $(\omega,x)$, the following hold for all $n\ge n_0(\omega,x)$:  
\begin{enumerate}[label= {(\alph*)}, ref= {(\alph*)}]
	\item\label{hh1} $\displaystyle \dfrac{\log     \mu^{     \xi ^{i-1}\vee\zeta}_{(\omega,x)}    \wtd B^{i-1}\left((\omega,x) , e^{n(\lambda_i+2\epsilon)} \right) } {n(\lambda_i + 2\epsilon)}\le \overline{\dim}^{i-1}(   \mu \mid \zeta) + \epsilon$;
	\item\label{hh2} { $\displaystyle - \frac 1 n \log   \mu^{  \xi ^{i-1} \vee\zeta}_{(\omega,x)} ( \xi^i\vee  \partP \vee\zeta)  ^{(-n,0]}(\omega,x) \ge h_{i-1} - \epsilon$}; 
	\item \label{hh3}$\displaystyle(  \xi^i\vee  \partP ) ^{(-n,0]}(\omega,x)  \subset \wtd W^i_{ \loc, e^{n(\lambda_i +2\epsilon)}} (\omega,x)\subset \wtd B^i \left((\omega,x), e^{n(\lambda_i +2\epsilon)}\right)$;
	\item \label{hh4} $\displaystyle - \frac 1 n \log\mu^{ \xi ^{i} \vee\zeta}_{(\omega,x)} ( \xi ^{i}\vee\zeta \vee \partP)  ^{(-n,0]}(\omega,x) \le h_i  + \epsilon$;  
	\item \label{hh5}$  \wtd B^{i-1}\left((\omega,x) , e^{n(\lambda_i+2\epsilon)} \right) \subset    \xi^{i-1}( \omega,x)$.
		\end{enumerate}Moreover, for infinitely many $n\ge  n_0{(\omega,x)}$,  
\begin{enumerate}[label= {(\alph*)}, ref= {(\alph*)}, resume]
	\item  \label{hh6} $\displaystyle \dfrac{ \log  \mu^{  \xi ^{i}\vee\zeta}_{(\omega,x)}  \wtd B^i\left({(\omega,x)} , {2} e^{n(\lambda_i+2\epsilon)} \right) } {n(\lambda_i +2\epsilon)}\ge \overline{\dim}^{i}(  \mu \mid \zeta) - \epsilon$.
	\end{enumerate} 
\end{lemma}
We note that \ref{hh2} and \ref{hh4} are slightly different from their analogues (2) and (4) of \cite[Lem.\  10.2]{MR819557} since we define $h_i$ by \cref{eqhi} rather than by Brin--Katok ``local entropies'' as in \cite[\S 9]{MR819557}.  Also \ref{hh3} is a simpler version of (3) of \cite[Lem.\  10.2]{MR819557} using \cref{lem:partadaptedtolocalentropy} rather than \cite[(10.1)]{MR819557}.
\begin{proof}[Proof of \cref{lem:ooo}]
\ref{hh1}, \ref{hh5}, and \ref{hh6} follow from definition. \ref{hh4} follows from Lemma \ref{lemma:pastfinitepart}.   \ref{hh3} follows from Lemma \ref{lem:partadaptedtolocalentropy}.

We have $F^n\xi^i\vee \xi^{i-1}=F^n\xi^{i-1}$ from \cref{cor:refinedynamically}.  Since $\xi^i$ and $\xi^{i-1}$ are sub-invariant, we thus have 
\begin{align*}
(\xi^i\vee  \partP \vee\zeta)  ^{(-n,0]}\vee \xi^{i-1}&=F^{n-1}\xi^i \vee (\partP \vee\zeta)  ^{(-n,0]} \vee  \xi^{i-1}  \\
&=F^{n-1}\xi^{i-1} \vee (\partP \vee\zeta)  ^{(-n,0]} \\
&=(\xi^{i-1}\vee  \partP \vee\zeta)  ^{(-n,0]}
\end{align*}
and so \ref{hh2} also follows from  Lemma \ref{lemma:pastfinitepart}.

   Note also that when $i-1=0$,  \ref{hh1} 
   and \ref{hh5} hold trivially  and \ref{hh2} follows from  Lemma \ref{lemma:pastfinitepart}.      
\end{proof}

\subsubsection{Proof of  inequality \ref{ineqI}}\label{sss:ineqI}
We now prove the first inequality of \cref{sec:outline}.
\begin{proof}[Proof of  \ref{ineqI}]Retain all notation (in particular, the choice of $\partP$) from Lemma \ref{lem:ooo}.  Fix a choice of $1\le i\le s$.   Recall  $h_i:= h_\mu(F, \xi^{i}\vee \zeta)$.

Fix $N_1>0$ and set $\Gamma := \{(\omega,x) : n_0(\omega, x) \le N_1\}$.
With $$g(\omega,x) := \mu^{ \xi ^{0}\vee \zeta}_{(\omega,x)} (\Gamma),$$ select $N_1$ large enough so that for a  positive measure subset of $(\omega,x)\in \Gamma$ satisfying Lemma \ref{lem:ooo}, we have 
 $g(\omega,x)>0.9$, and
for all $n\ge N_1$, 
\begin{align*} &\frac{  \mu^{ \xi ^{i-1}\vee \zeta}_{(\omega,x)} \left( \Gamma \cap  \wtd   B^{i-1}((\omega,x) , e^{n(\lambda_i+2\epsilon)} ) \right  )  }{  \mu^{ \xi ^{i-1}\vee \zeta}_{(\omega,x)} \left(  \wtd   B^{i-1}((\omega,x) , e^{n(\lambda_i+2\epsilon)} ) \right  )  }
\\
&\quad \quad =
\frac 1{  \mu^{ \xi ^{i-1}\vee \zeta}_{(\omega,x)} \left(  \wtd  B^{i-1}((\omega,x) , e^{n(\lambda_i+2\epsilon)} ) \right  )  }\int_{ \wtd  B^{i-1}((\omega,x) , e^{n(\lambda_i+2\epsilon)}) }
g(\cdot) \, d\mu^{ \xi ^{i-1}\vee \zeta}_{(\omega,x)} (\cdot)
\\
  &\quad \quad  \ge 0.87.
  \end{align*}
Such $(\omega,x)$ exist by the Lebesgue differentiation theorem applied to the pairs of measures $$ \tau_*
\left( \restrict{\mu^{ \xi ^{i-1}\vee \zeta}_{(\omega,x)}}{\Gamma}\right)\ll
 \tau_*   \mu^{ \xi ^{i-1}\vee \zeta}_{(\omega,x)} $$
  on the manifolds $W ^{i-1}({\tau(\omega,x)})$. 
Fix such  $(\omega,x)$.  
Observe  by Lemma \ref{lem:ooo}\ref{hh1} that for $n\ge N_1$, 
\begin{equation}
  \mu^{ \xi ^{i-1}\vee \zeta}_{(\omega,x)} \left( \Gamma \cap \wtd  B^{i-1}((\omega,x) , e^{n(\lambda_i+2\epsilon)} ) \right  )
\ge 0.87 e^{(n(\lambda_i +2\epsilon))  ( \overline{\dim}^{i-1}(  \mu\mid  \zeta) + \epsilon)}. 
\end{equation}

For $n\ge N_1$,  take $$L_n:=   \wtd B^{i-1}((\omega,x) , e^{n(\lambda_i+2\epsilon)} )\cap \zeta(\omega,x).$$  
The careful reader will observe that unlike in  other adaptations of Ledrappier--Young, our ``balls''  $\wtd B^{i-1}((\omega,x) , e^{n(\lambda_i+2\epsilon)} )$ are  fully saturated by $\wtd \sigA$-sets.  In particular, $L_n$ witnesses  no decay in the $\xi^0$-directions despite seeing decay in the fiberwise $W^{i-1}$-directions (if $i\ge 2$); this allows us to avoid defining a notion of pointwise dimension along the $\xi^0$-atoms (which we view as the ``strongest'' stable directions).
For $(\omega' ,y)\in \Gamma \cap L_n $, we have $\tau(\omega')=\tau(\omega)$;  for $n\ge N_1$,   Lemma \ref{lem:ooo}\ref{hh5} then implies   that $ \xi^{i-1}(\omega', y) = \xi^{i-1}(\omega,x)$ and so   Lemma \ref{lem:ooo}\ref{hh5} again implies
$  \wtd B^{i-1}\left((\omega',y), e^{n(\lambda_i+2\epsilon)} \right) \subset   \xi^{i-1}(\omega, x)$.  Thus, by Lemma \ref{lem:ooo}\ref{hh2}, for $(\omega' ,y)\in \Gamma \cap L_n $ and $n\ge N_1$, 
$$   \mu^{ \xi ^{i-1}\vee \zeta}_{(\omega,x)}  ( \xi^i\vee \partP\vee \zeta)  ^{(-n,0]}(\omega',y)
=   \mu^{ \xi ^{i-1}\vee \zeta}_{(\omega',y)} ( \xi^i\vee \partP\vee\zeta)  ^{(-n,0]}(\omega', y)\le e^{-n(h_{i-1} - \epsilon)}.$$

For $n\ge N_1$, we obtain a lower   bound  on the number of distinct $(  \xi^i \vee \partP \vee\zeta) ^{(-n,0]}$-atoms meeting   $\Gamma \cap L_n $ by 
\begin{align*} \#\{( \xi^i\vee \partP\vee\zeta)  ^{(-n,0]}(\omega',y)  & : (\omega',y)\in  \Gamma \cap L_n \} 
\ge   \mu^{ \xi ^{i-1}\vee\zeta}_{(\omega,x)} \left( \Gamma \cap L_n  \right  )  / e^{-n(h_{i-1} - \epsilon)} \\
&\ge  0.87 e^{(n(\lambda_i +2\epsilon))  ( \overline{\dim}^{i-1}(  \mu\mid \zeta) + \epsilon)}  e^{n(h_{i-1} - \epsilon)}.  
\end{align*}

For $(\omega',y)\in \Gamma \cap L_n  $ and $n\ge N_1$, we have by  Lemma \ref{lem:ooo}\ref{hh3} that 
$( \xi^i\vee \partP\vee\zeta)  ^{(-n,0]}(\omega',y)\subset   \wtd  B^i \left((\omega',y), e^{n(\lambda_i + 2\epsilon)}\right)$ whence $( \xi^i\vee \partP\vee\zeta)  ^{(-n,0]}(\omega',y)\subset 
 \wtd B^{i}\bigl((\omega,x) ,  2 e^{n(\lambda_i+2\epsilon)} \bigr).$
From Lemma  \ref{lem:ooo}\ref{hh4}, we have 
$$  \mu^{ \xi ^{i}\vee \zeta}_{(\omega',y)} ( \xi^i\vee  \partP\vee\zeta)  ^{(-n,0]}(\omega',y) \ge e^{-n( h_i  + \epsilon)}$$
and obtain inequalities  
\begin{align*}
  \mu^{ \xi ^{i}\vee\zeta}_{(\omega,x)}& \wtd   B^{i} \left((\omega,x) , {2} e^{n(\lambda_i+2\epsilon)} \right) \\
 &\ge  \#\{( \xi^i\vee \partP\vee\zeta)  ^{(-n,0]}(\omega', y)  : (\omega', y)\in  \Gamma \cap L_n \} \cdot   e^{-n( h_i + \epsilon)}\\
  &\ge   0.87 e^{(n(\lambda_i +2\epsilon))  ( \overline{\dim}^{i-1}( \mu\mid \zeta ) + \epsilon)}  e^{n(h_{i-1} - \epsilon)} e^{-n( h_i  + \epsilon)}.
\end{align*}

Comparing to Lemma \ref{lem:ooo}\ref{hh6},
we have for infinitely many $n$ that 
	$$(\lambda_i + 2\epsilon)  (\overline{\dim}^{i}(  \mu\mid\zeta) - \epsilon )
	\ge\left( \frac {\log (0.87)} n+
(\lambda_i + 2\epsilon)  ( \overline{\dim}^{i-1}( \mu\mid\zeta) + \epsilon)
	 - 
	h_i+h_{i-1} - 2\epsilon\right).
	$$
 Choosing such an $n$ sufficiently large, 
 $$h_i-h_{i-1} \ge -(\lambda_i + 2\epsilon) \left (\overline{\dim}^{i}(   \mu \mid \zeta)  -  \overline{\dim}^{i-1}(  \mu \mid \zeta) - 2 \epsilon\right) - 3 \epsilon.$$
Inequality  \ref{ineqI} then follows by the arbitrariness of  $\epsilon>0 $.   
\end{proof}

\subsubsection{Addendum: proof of \cref{thm:exactdimproj}}\label{sss:proj msr}
After slight reformulation,  the proof of \cref{thm:exactdimproj} is essentially the same as the proof of  inequality \ref{ineqI} in \cref{sss:ineqI}.  We outline the details.

Recall we assume the base dynamics $\hat \theta \colon (\what \Omega,\hat \nu) \to (\what \Omega,\hat \nu)$ admits a one-sided generator $\what \partQ$.  Then (cf.~\cite[Lem.~3.9]{MR0262464}) we have $$h_{\hat\nu}(\hat \theta) = H_{\hat \nu}(\ptpart\mid \hat \theta\inv \ptpart).$$
By \cref{prop:GenAR,prop:howtoNatExt}, we have $\maxent_{   \mu}(  F\mid    \sigC)=h_ \mu( F)$.  Indeed,
\begin{align}
\maxent_{   \mu}(  F\mid    \sigC) &= \maxent_{   \hat \mu}(  \what F\mid   \hat  \sigF)
=  h_{   \hat \mu}(  \what F\mid   \hat  \sigF) +  h_{   \hat \mu}(  \what F,  \Xi( \hat  \sigF)) 
\notag\\&
=  h_{   \hat \mu}(  \what F\mid   \hat  \sigF) +  H_{   \hat \nu}( \ptpart \mid \hat \theta\inv(\ptpart))
=  h_{   \hat \mu}(  \what F\mid   \hat  \sigF) +  h_{ \hat \nu}( \hat \theta)
\notag\\&
= h_{\hat \mu}(\what F)\label{eqAR15}
= h_ \mu( F).\end{align}
Here $\ptpart$ denotes the point partition on the base space $(\what \Omega, \hat \nu)$.  The penultimate equality in \eqref{eqAR15} is the classical Abramov--Rohlin formula.

\begin{proof}[Proof of \cref{thm:exactdimproj}]
The assumption that $\lambda_{\top}(\hat \mu)<0$ ensures stable manifolds are open in the fibers of $\what X$.  The  measures $\hat \mu_\omega$ are exact dimensional with an a.s.~common dimension $d$ by  \cref{prop:dimensions are defined} showing \cref{thm:exactdimproj}\ref{projdima}.

Given $x\in M$, write 
 $$\underline \delta (x)= \liminf_{r\to 0}\frac{\log \mu_0(B(x,r))}{\log r}, \quad \quad \overline \delta (x)= \limsup_{r\to 0}\frac{\log \mu_0(B(x,r))}{\log r}.$$
\cref{thm:exactdimproj}\ref{projdimb} follows if for $\hat \mu$-a.e.~$(\omega,x)$, we have  \begin{equation}
d\le \underline \delta (x)\le  \overline \delta (x)\le d
\end{equation}
Viewing every $\hat \mu_\omega$ as a measure on $M$, we have $\mu_0 = \int \hat \mu_\omega\, d \hat \nu(\omega)$. Then by \cite[Lem.~11.3.2]{MR819557}, for $\hat \mu$-a.e.~$(\omega,x)$ we obtain the lower bound $$d\le \underline \delta (x).$$

It remains to show that $\overline \delta (x)\le d$ for $\hat \mu$-a.e.~$(\omega,x)$.  The proof follows the proof of  inequality \ref{ineqI} in \cref{sss:ineqI}.
Fix $\epsilon>0$ sufficiently small; in particular, we assume $\lambda=\lambda_{\top}(\hat\mu) + 2\epsilon<0$.  
Let $r_n = e^{(n-1)\lambda}$ and let $$h_*=h_{\mu}(F)=\maxent_{   \mu}(  F\mid    \sigC).$$  Let $\xi^s$ be the pullback under $\tau$ of a  partition subordinate to fiber stable manifolds as in \cref{prop:entropyisunstable}.  Note that atoms of $\xi^s$ are $\wtd \sigA$-saturates of open subsets in the fibers of $X$.  In particular, atoms of $\xi^s$ have positive $\mu^\sigC_{(\omega,x)}$-measure a.s.

As in  Lemma \ref{lem:ooo}, there exists a partition $\partR$ of $  X$ with $H_\mu(\partR)<\infty$ and 
a measurable function $n_0\colon X\to \N$ such that for $  \mu$-a.e.\  $(\omega,x)$, the following hold for all $n\ge n_0(\omega,x)$:  
\begin{enumerate}[label= {(\alph*)}, ref= {(\alph*)}]
	
	\item \label{jhh4} $\displaystyle - \frac 1 n \log\mu^{ \sigC}_{(\omega,x)} ( \partR  ^{(-n,0]}(\omega,x)) \ge h_*  - \epsilon$;  

	\item \label{jhh3}
	$\displaystyle \diam \left(\pi_2\bigl( \partR  ^{(-n,0]}(\omega,x)\bigr)\right)\le r_n$
	\item\label{jhh2} { $\displaystyle - \frac 1 n \log   \mu  (   \partR  ^{(-n,0]}(\omega,x) )\le h_{*} +\epsilon$}.  
		\end{enumerate}Moreover, 
\begin{enumerate}[label= {(\alph*)}, ref= {(\alph*)}, resume]
	\item\label{jhh1}  $\displaystyle\lim_{n\to \infty}   \dfrac{ \log  \left(\hat \mu_{\tau(\omega)} \left(  B(x, r_n) \right) \right)} {\log r_n}=
	\displaystyle\lim_{n\to \infty} \dfrac{ \log  \left(\mu_{(\omega,x)}^{\sigC} \left(\wtd  B((\omega,x), r_n) \right) \right)} {\log r_n}	= d$;

	\item  \label{jhh6} $\displaystyle  \limsup_{n\to \infty} \dfrac{ \log  \left(\mu_0 \left( B(x, 2r_n) \right) \right)} {\log r_n}= \displaystyle \limsup_{n\to \infty} \dfrac{ \log  \left(\mu \left( \Omega\times B(x, 2r_n) \right) \right)} {\log r_n}= \overline{\delta}(x)$.
	\end{enumerate} 

\begin{proof}[Proof sketch]
Properties \ref{jhh1} and \ref{jhh6} follow from definition.  

For property \ref{jhh3}, first take  $\partP$ as in the construction of a partition adapted to suitable (two-sided) Lyapunov charts as in  \cref{lem:subpartfinent}.  
Recall the functions $T_\pm(\omega,x)$ in the proof of \cref{lem:subpartfinent}.  
Let $\partQ= \tau\inv \what \partQ$, let $\wtd \partQ= \pi_1\inv \partQ$ and let $\partR=\partP\vee \wtd \partQ$.  
Recall $\Lambda$ and $T(\omega,x)$ in the proof of \cref{lem:subpartfinent}. 
If  $$(\omega',y)\in \partR^{[0,n)}(\omega,x)\subset \wtd \partQ^{[0,n)}(\omega,x)$$ then \eqref{eq:pesinadapt} in the proof of  \cref{lem:subpartfinent} holds.  
If $n\ge T_+(\omega,x)+1$, then for $(\omega',y)\in \partR^{[0,n)}$, 
$$ \pi_2\bigl(\partR^{[0,n)}(\omega',y)\cap \sigF(\omega',y)\bigr)\subset \Phi_{(\omega,x)}(\R^d(\delta)).$$
From  \eqref{eq:pesinadapt}, the dynamics inside Lyapunov charts, and that Lyapunov charts are $1$-Lipschitz, we have 
$$\pi_2\bigl(F^{n-1}\bigl(\partR^{[0,n)}(\omega',y)\cap \sigF(\omega',y)\bigr) \bigr)\subset B(f_\omega^{(n-1)}(x), \delta r_n).$$
Now, consider $(\omega,x)\in X$.  If $n\ge T_-(\omega,x)+2$, we have $T_+(F^{-(n-1)}(\omega,x))\le n-1$.  
Since  $$\partR^{(-n,0]}(\omega,x)=F^{n-1}\partR^{[0,n)}(F^{-(n-1)}(\omega,x))$$  the above discussion applies at the point $F^{-(n-1)}(\omega,x)$ and
it follows $$\pi_2(\partR^{(-n,0]}(\omega,x))\subset B(x, \delta r_n)$$ and so $\diam( \pi_2(\partR^{(-n,0]}(\omega,x)))\le r_n$ assuming $\delta\le \frac 1 2 $ and  property \ref{jhh3} follows.

Property \ref{jhh3} ensures $\partR^+$ (and in fact $\partP^+$) distinguishes points in the same fiber.  Since $\wtd \partQ^\Z$ is the partition into fibers, $\partR^\Z$ is the point partition on $X$.  In particular, $\partR$ is a two-sided generator and, in particular, a generator relative to $\sigC$.  Since $\what \partQ$ is a one-sided generator, we have $\Xi(\sigC) = \wtd\partQ^{[0,\infty)}$.  Thus for $n\ge 1$ we have  $$ \Xi(\sigC)^{(-n,0]}= F^{n-1}\wtd\partQ^{[0,\infty)}= \wtd\partQ^{(-n,\infty)}= \wtd\partQ^{(-n,0]}\vee   \Xi(\sigC).$$ Since $\wtd \partQ\prec \partR$, we thus have $$\partR^{(-n,0]}\vee\sigC= (\partR\vee\sigC)^{(-n,0]}.$$
Property \ref{jhh4} then follows from \cref{cor:SMBgen} using that $\maxent_\mu(F\mid \sigC) = h_*$.  
For Property \ref{jhh2}, we may further refine  $\partR$ if needed to ensure  that $$h_\mu(F, \partR)\ge h_\mu(F)-\frac{\epsilon}{2}.$$  Then \ref{jhh2}  follows from the classical Shannon--McMillan--Breiman Theorem.      
 \end{proof}

We return to the proof of  \cref{thm:exactdimproj}.  Fix $N_1>0$ and set $\Gamma := \{(\omega,x) : n_0(\omega, x) \le N_1\}$.
Differentiation and density arguments as in \cref{sss:ineqI} imply we can select $N_1$ sufficiently large  so that for a set  of $(\omega,x)\in \Gamma$ satisfying properties \ref{jhh4}---\ref{jhh6} above with measure arbitrarily close to $\mu(\Gamma)$, we have
$$\mu_{(\omega,x)}^{\sigC} \left(\Gamma \cap \wtd  B((\omega,x), r_n) \right)\ge .87 \mu_{(\omega,x)}^{\sigC} \left(\wtd  B((\omega,x), r_n) \right)
=.87\hat \mu_{\tau(\omega)} (B(x, r_n)).$$

Fix such $(\omega,x)$.  
For $n\ge N_1$, using property \ref{jhh4} we  lower   bound the number of distinct $ \partR  ^{(-n,0]}$-atoms meeting   $\Gamma \cap  \wtd B((\omega,x) , r_n)  $ by 
\begin{align*} \#\{\partR  ^{(-n,0]}(\omega',y)  & : (\omega',y)\in  \Gamma \cap \wtd B((\omega,x) , r_n)  \} 
\\&\ge   \mu^{ \sigC}_{(\omega,x)} \left( \Gamma \cap \wtd B((\omega,x) , r_n)   \right  )  / e^{-n(h_{*} - \epsilon)} \\
&\ge .87 \hat \mu_{\tau(\omega)} \left(  B(x, r_n) \right) e^{n(h_{*} - \epsilon)}.  
\end{align*}

For $(\omega',y)\in \Gamma \cap \wtd B((\omega,x) , r_n)   $ and $n\ge N_1$, by  property \ref{jhh3} we have 
$\pi_2\left(\partR  ^{(-n,0]}(\omega',y)\right) \subset B(x,2r_n) $.
Combining the above bound with property \ref{jhh2}, we have 
\begin{align*}
 \mu_0(B(x,2r_n))&= \mu  \left(\Omega\times B(x, 2r_n)\right)  \\
 &\ge  \#\{ \partR  ^{(-n,0]}(\omega', y)  : (\omega',y)\in  \Gamma \cap \wtd B((\omega,x) , r_n)  \} \cdot   e^{-n( h_* + \epsilon)}\\
  &\ge   .87 \hat \mu_{\tau(\omega)} \left(B(x, r_n) \right)    e^{n(h_{*} - \epsilon)} e^{-n( h_*  + \epsilon)}\\
    &=  .87 \hat \mu_{\tau(\omega)} \left(   B(x, r_n) \right)    e^{-2n \epsilon} .
\end{align*}
Comparing to property \ref{jhh6},  we have
$$\overline \delta (x)  \le 
\limsup_{n\to \infty} \dfrac{\log\left(.87 \hat \mu_{\tau(\omega)} \left(  B(x, r_n) \right)    e^{-2n \epsilon} \right)}{\log r_n} =
d + \frac{2\epsilon}{-(\lambda_\top(\hat\mu)+2\epsilon)}.$$
The bound $\overline \delta (x)\le d$ then follows.  
\end{proof}

\subsection{Auxiliary partitions and their metric properties}\label{sec:thespartitions} 
We build a number of auxiliary objects used to prove inequalities \ref{ineqII}--\ref{ineqV}. 

We first note  that if the measure  $\mu$ contains atoms, then by ergodicity, the measure is finitely supported.  In this case all entropies and dimensions appearing in  \ref{ineqII}--\ref{ineqV} vanish.      
Thus we assume for the remainder that $\mu$ has no atoms.

Fix a choice of $0<\epsilon<\eps$ and a family of two-sided $(\epsilon, \ell)$-Lyapunov charts $\{\Phi_{(\omega,x)}\}$ as in \cref{sec:charts}.  
Consider a fixed  family of measurable partitions  $\{\hat  \xi^i\}_{ 0\le i\le s} $  of $(\what X, \hat \mu)$ as in \cref{good parts}.   We recall  also the choice of $ \Lambda_1\subset \what X$, $x_0\in M$, $ S_\rho$ 
in \cref{sec:sss}.   
Also recall the two-sided $(\epsilon, \ell)$-Lyapunov charts $\{\Phi_{(\omega,x)}\}$ and their properties outlined in \cref{sec:twosidecharts}.

\subsubsection{Auxiliary partitions}
Let $\delta_0$ be as in \cref{delta0}.  
 Write $\wtd S:= \tau\inv S_\rho\subset X$.    
 Fix $\ell_0>0$  sufficiently large so that $\wtd S'=\{(\omega,x)\in \wtd S: \ell(\omega,x)\le \ell_0\}$ has positive measure in $B(x_0, \rho)$.  Then select 
  $x_1\in B(x_0, \rho)$ with 
$$\mu \left(\left(\Omega\times B\left(x_1, \frac{ \delta_0 \ell_0 \inv} 2\right)\right) \cap \wtd S' \right)>0$$ 
 and $\ell(\omega,x)\le \ell_0$ for all $(\omega,x)\in \wtd S'$.
Such an $x_1$ can be taken by considering points where $$\frac{d(\pi_2)_*\restrict{\mu}{\wtd S'}}{d(\pi_2)_*\restrict{\mu}{\wtd S}}>0.$$

For $1\le i\le\min\{ s+1,p\}$, we select a $(m_{i}+\dots +m_p)$-dimensional embedded disc
 $$T_i \subset B\left(x_1, \frac{\delta_0\ell_0\inv}{2}\right). $$
We may assume the $T_i$ are chosen so that there exists a subset   $$E\subset \left(\Omega\times B\left(x_1, \frac{ \delta_0 \ell_0 \inv} 2\right)\right) \cap \wtd S'$$ with $\mu(E)>0$ 
  and such that for every $(\omega,x)\in E$,
   $$W^{i-1}_{ \loc, 2\delta_0}(\omega,x)\subset M$$ intersects $T_i$ in a singleton
and    $$W^{i}_{ \loc, 2\delta_0}(\omega, x)\subset M$$ intersects $T_i$ in an $m_i$-dimensional disc; moreover,  all such intersections are uniformly (over $(\omega,x)\in E$) transverse to $T_i$.

Let $\partP$ be a measurable partition of $(X,\mu)$ with the following properties: 
\begin{enumerate}
\item $  \partP$ is adapted    to $(\{\Phi_{(\omega, x)}\},\delta_0)$ (see  Definition \ref{def:adapted}); 
\item $  \partP$ refines $\{\wtd S,   X \sm \wtd S\}$; 
\item $  \partP$ refines $\{E,   X\sm E\}$;     
\item $  \partP$ refines $\{E',   X\sm E'\}$ where $E'\subset E$ is a subset with $\mu(E')>0$ to be specified below;
\item  $ H_\mu(  \partP)<\infty$.  
\end{enumerate}

Following the notation in \cref{little props about eta}, for $0\le i\le s$ take $$ \eta_i =    \partP^-\vee   \xi^i \vee\zeta.$$Recall also the notational convention
 $\eta_*= \eta_{s+1}:=(\partP\vee \zeta)^-$ and so if $\lambda_s=\lambda_p<0$ or if $\lambda_s<0< \lambda_{s+1}$ (in particular, if no exponent is zero), then  $  \eta_{s+1}=  \eta_s.$

\subsubsection{Preliminary choices of transverse metrics   on $\eta_i/\eta_{i-1}$}

Consider $(\omega,x)\in \Lambda_0$.  For $1\le i\le s+1$, let  $U^i= \bigoplus _{j\ge i} \R_j$ (with the convention the empty direct sum is $\{0\}$).
Given $(\omega',y),(\omega'',z)\in \eta_*{(\omega,x)}$ we have that 
$$\tau(\omega) = \tau(\omega')=\tau(\omega'').$$
From \cref{lem:localunstable}, the intersections 
$$  W^{i-1}_{{(\omega,x)},2\delta_0}(\Phi_{(\omega,x)}\inv(y))\cap U^i$$
$$  W^{i-1}_{{(\omega,x)},2\delta_0}(\Phi_{(\omega,x)}\inv(z))\cap U^i$$ are singletons.  
Let $$d^i_{(\omega,x)}((\omega',y),(\omega'',z))= \|   W^{i-1}_{(\omega,x),2\delta_0}(\Phi_{(\omega,x)}\inv(y))\cap U^i -   W^{i-1}_{(\omega,x),2\delta_0}(\Phi_{(\omega,x)}\inv(z))\cap U^i \|.$$

Recall $\pi_2\colon X\to M$ is the canonical projection.  By \cref{claim:neatfacts}\ref{yuff}, and for every $(\omega',y)\in \eta_*(\omega,x)$,  
 $$\Phi_{(\omega,x)}\inv \bigl(\pi_2(\eta_{i-1}(\omega',y))\bigr) \subset W^{i-1}_{{(\omega,x)},2\delta_0}(\Phi_{(\omega,x)}\inv(y)).$$ 
Combined with \cref{lemma:44}\ref{PP1}, $d^i_{(\omega,x)}$ defines a metric on $\eta_{i-1}$-equivalence classes in $\eta_*(\omega,x)$; below we will restrict   $d^i_{(\omega,x)}$ to a metric on $\eta_{i-1}$-equivalence classes in $\eta_i(\omega,x)$, hence the notation in the superscript.  
As in  \cite[Lem.\  2.3.2]{MR819556} and  \cite[Lem.\  8.3.2]{MR819557}, applying the dynamics inside charts, we have the following.  
\begin{claim}\label{growthonorbit} For $(\omega'',z)\in \eta_i(\omega, x)$
 $$d^i_{F(\omega,x)}(F(\omega,x),F(\omega'',z))\ge e^{\lambda_i - 3\epsilon} d^i_{(\omega,x)}((\omega,x),(\omega'',z)).$$
\end{claim}

We construct an alternative metric that is independent of the choice of $x$.   This is slightly different from the construction in  \cite[(4.2)]{MR819556} and \cite[(8.4)]{MR819557}; in particular, we use the Lipschitzness of holonomies in the setting of $C^{1+\beta}$ diffeomorphisms established in  \cite{MR4504676}.
By the choice of the set $E$  above, for each $1\le i \le \min\{s+1,p\}$, we have a  $(m_i+ \dots + m_p)$-dimensional embedded disc 
	$$T_i \subset B\left (x_1,  \tfrac 1 2\delta_0\ell_0\inv\right) $$
  that is uniformly transverse to each $W^{i-1}_{\loc,\delta_0}(\omega',y)$ for   $(\omega',y)\in E$  and 
such that for $1\le i\le s$, 	$W^{i}_{\loc, 2\delta_0}(\omega',y)\cap T_i $ is an embedded $m_i$-dimensional submanifold for each $(\omega',y)\in E$.

Identify each transversal $T_i$ with a subset of $\R^{m_i+ \dots + m_p}$ through a bi-Lipschitz identification.  Let $d_{T_i}$ denote the pullback to $T_i$ of the Euclidean distance under this identification.  (We fix a Euclidean geometry on $T_i$ to apply Besicovitch arguments later.)
Given $(\omega,x)\in E$ and $(\omega',y),(\omega'',z)\in \eta_*(\omega,x)$,  we  define $d^{T_i}((\omega',y),(\omega'',z))$ to be \begin{align*}d^{T_i}((\omega',y),(\omega'',z))&:= d_{T_i}(W^{i-1}_{\loc,2\delta_0}(\omega', y)\cap T_i, W^{i-1}_{\loc,2\delta_0}(\omega'',z)\cap T_i)\end{align*}
For $(\omega,x)\in E$, \cref{claim:insides} and \cref{lemma:44}\ref{PP1}  imply that $d^{T_i}$ induces a metric on $\eta_{i-1}$-equivalence classes in $\eta_*(\omega,x )$.

For $(\omega,x)\in E$, let $\wtd T_i(\omega,x) = \Phi_{(\omega,x)}\inv( T_i)\subset \R^d$.  View $\wtd T_i(\omega,x)$ as a subspace of Euclidean metric.  We have that $\wtd T_i(\omega,x)$ and $T_i$ are   bi-Lipschitz equivalent with Lipschitz constant uniform  (by the choice of $\ell_0$) over $(\omega, x)\in E$.   
 Given $(\omega,x)\in E,$ consider the holonomy map  $$U^i\cap   W^{i}_{(\omega,x),2\delta_0}(0) \to  \wtd T_i(\omega,x) \cap   W^{i}_{(\omega,x),2\delta_0}(0)$$  along $  W^{i-1}_{(\omega,x),2\delta_0}$-leaves.  The main result of \cite{MR4504676} ensures this map is uniformly   bi-Lipschitz (with Lipschitz constant uniform over $(\omega,x)\in E$).
 Given $(\omega,x)\in E$ and $(\omega',y),(\omega'',z)\in \eta_*(\omega,x)$, define  \begin{align*} &\td d^{ \wtd T_i}_{(\omega,x)}  ((\omega',y),(\omega'',z))\\&  \quad : =\left\|\left (  W^{i-1}_{(\omega,x),2\delta_0} (\Phi_{(\omega,x)}\inv (y))\cap \wtd T_i(\omega,x) \right)-\left(   W^{i-1}_{{(\omega,x)},2\delta_0}(\Phi_{(\omega,x)}\inv (z)) \cap \wtd T_i{(\omega,x)} \right)\right\|.\end{align*}  
 It follows that, restricted to $\eta_{i-1}$-equivalence classes in $\eta_i(\omega,x) $, 
  the metrics $d_{(\omega,x)}^i$ and $\td d^{\wtd T_i}_{(\omega,x)}$ are uniformly   bi-Lipschitz equivalent (with Lipschitz constant uniform over ${(\omega,x)}\in E$).
In particular, 
using that the charts $\{\Phi_{(\omega,x)} :{(\omega,x)}\in E\}$ are uniformly Lipschitz embeddings, we obtain the  following.   
\begin{claim}\label{claim:DontGiveMeNoLipBoy}
There is an  $N>1$ (depending only on  $\ell_0$ and $E$) such that for all ${(\omega,x)}\in E$ and ${(\omega',y)},{(\omega'',z)}\in \eta_i{(\omega,x)}$, 
$$\frac 1{   N} d^{T_i}({(\omega',y)},{(\omega'',z)}) \le d_{(\omega,x)}^i({(\omega',y)},{(\omega'',z)}) \le   N d^{T_i}({(\omega',y)},{(\omega'',z)}).$$ 
\end{claim}

We note the $N$ in Claim \ref{claim:DontGiveMeNoLipBoy} depends only on the choice of the set $E$ (on which $\ell_0$ is bounded from above and relative to which the transversals $T_i$ were chosen) but not on the choice of $E'\subset E$ in the choice of $\partP$ in \cref{sec:thespartitions} which we now specify.

\subsubsection{Choice of transverse metric on $\eta_i/\eta_{i-1}$ and transverse dimension}   
We now define the metric used for the proof of \cref{prop:key} below.  

We take $E'\subset E$ to be any subset with $$0<\mu(E')<\frac {\epsilon ^2}{8\log N}.$$
We will define a metric on $\eta_{i-1}$-equivalence classes in each $\eta_i$-atom relative to which all hyperbolicity is concentrated on returns to $E'$.  The choice of $E'$ with small measure will be used to control bi-Lipschitz changes of metric that arise over subsequent returns to $E'$.

For   ${(\omega,x)}\in \Lambda_0$, let $$n{(\omega,x)}:= \min \{n\ge 0 : F^{n}{(\omega,x)} \in E'\}.$$  By the choice of $\partP$ refining $\{E', X\sm E'\}$,  $n{(\omega,x)}$ is constant on elements of $\eta_*$. 
 For  $1\le i \le \min\{s+1,p\}$ and for $(\omega',y),(\omega'',z)\in \eta_*{(\omega,x)}$, 
define \begin{equation}\label{eq:metrics} \bar d^{i}((\omega',y),(\omega'',z)):= d^{T_i}(F^{n{(\omega,x)}}(\omega',y), F^{n{(\omega,x)}}(\omega'',z)).\end{equation}
As $E'\subset E$ and as $F\eta_{i-1}{(\omega,x)}\subset \eta_{i-1}(F{(\omega,x)})$, from \cref{lemma:44}\ref{PP2} we have that $ \bar d^{i}$ defines a metric on $\eta_{i-1}$-equivalence classes in every  $\eta_i{(\omega,x)}$.  In particular,  we view the space of equivalence classes as a topological quotient $\eta_i{(\omega,x)}/\eta_{i-1}$ endowed with this metric.

For $1\le i\le \min\{s+1,p\}$, let 
$$\wtd U^{i}({(\omega,x)}, \rho) := \{ (\omega',y) \in \eta_{i } {(\omega,x)}: \bar d^{i}((\omega',y),{(\omega,x)})< \rho\}$$ denote  the ball of $\eta_{i-1}$-equivalence classes in $\eta_i{(\omega,x)}$ of radius $\rho$ centered at ${(\omega,x)}$;  note, in particular, that $(\omega',y)\in \wtd U^{i}({(\omega,x)}, \rho) $ implies $\eta_{i-1}(\omega',y)\subset \wtd  U^{i}({(\omega,x)}, \rho) $.

For $1\le i \le \min\{s+1,p\}$, write  \begin{equation}\label{eq:dimmmm}\td \gamma_i{(\omega,x)}= \liminf _{\rho\to 0}\dfrac {\log \mu^{\eta_i}_{(\omega,x)} (\wtd  U^{i}({(\omega,x)}, \rho) )}{\log \rho}.\end{equation}

Standard lemmas in measure theory (cf. \cite[Lem. 11.1.4]{MR819557}) show the following.   
\begin{claim}\label{claim:dimbounds}
For $1\le i\le s$ and  a.e.\ ${(\omega,x)}$, $\td \gamma_i{(\omega,x)}\le m_i$ and if $s<p$ then $\td \gamma_{s+1}{(\omega,x)}\le \dim (M)$.  
\end{claim}

The following claim summarizes the discussion in  \cite[(11.3)--(11.4)]{MR819557} (using the Lipschitzness of holonomies established in \cite{MR4504676} rather than the coordinates of \cite[(8.4)]{MR819557}). 

\begin{claim}\label{claim:oink}For $1\le i\le s$, suppose $c\ge0$ is such that $\td \gamma_i{(\omega,x)}\ge c$\ for almost every ${(\omega,x)}$.  Then 
$$c \le \underline{\dim}^i( \mu\mid \zeta) - \underline{\dim}^{i-1}(\mu\mid \zeta) .$$
\end{claim}

 	\subsection{Key Proposition and proof of \ref{ineqII}, \ref{ineqIIII}, and  \ref{ineqV}}
	Recall we assume $\mu$ has no atoms and so may apply the discussion and constructions in \cref{sec:thespartitions}.   
	
Recall from discussion in \cref{finiteentass} that $h_\mu(F,\eta_i) <\infty$ for every $1\le i\le s+1$.  
	The following proposition is the analogue of \cite[Prop.\  5.1]{MR819556} and \cite[Prop.\  11.2]{MR819557}.

		\begin{proposition}\label{prop:key}
		For the family of partitions 
$$ \eta_{s+1} \prec \eta_s \prec \dots \prec \eta_2 \prec \eta_1\prec \eta_0=\xi^0\vee \zeta\vee\partP^{-}$$ 
		as constructed in \cref{sec:thespartitions}, for every $1\le i\le \min\{s+1,p\}$ with $\lambda_i\le 0$ and  for $\mu$-a.e.\ $(\omega,x) \in X$,  we have 
\begin{equation}(-\lambda_i + 3\epsilon)\td \gamma_i(\omega,x)\ge(1-\epsilon) \left(h_\mu(F,\eta_i) -  h_{\mu}(F,\eta_{i-1}) - 2\epsilon  \right)\end{equation}
	where $\td \gamma_i(\omega,x)$ is the transverse dimension from  \cref{eq:dimmmm} associated to $  \eta_i({\omega, x})/  \eta_{i-1}$.
	\end{proposition}

Inequality \ref{ineqII} now follows from  Proposition \ref{prop:key}, Claims \ref{claim:oink} and \ref{claim:dimbounds}, and the arbitrariness of $\epsilon>0$.   
Indeed,  using   \cref{cor:entropysame}, for $1\le i\le s$, the quantities $$h_\mu(F, \eta_i) = h_\mu(F,  \xi^i\vee \zeta)$$ are independent of the choice of the partition $\partP$.  Then the right-most inequality in \ref{ineqII} follows from Proposition \ref{prop:key} and Claim \ref{claim:dimbounds}.   
The first inequality in \ref{ineqII} follows from Claim \ref{claim:oink} with 
	$$c =\max\left\{0,\frac{(1-\epsilon) \left(h_\mu(F,\eta_i) -  h_{\mu}(F,\eta_{i-1}) - 2\epsilon  \right)}{(-\lambda_i + 3\epsilon)} \right\}$$ and Proposition \ref{prop:key}.

For \ref{ineqIIII}, 
first note by definition that  
$$h_\mu(F,  \xi^s\vee \zeta) = h_\mu(F,  \xi^s\vee \partP \vee \zeta) \le  \maxent_\mu(F\mid \sigma(\zeta)).$$Moreover, given  any 
$M< \maxent_\mu(F\mid \sigma(\zeta))$ we may  assume $\partP$ in \cref{sec:thespartitions} was further chosen so  that \begin{equation} M<h_\mu (F,  \partP \vee \zeta ):= h_\mu(F, \eta_{s+1} ).\end{equation}
{If $\lambda_{s+1} >0$ then $\eta_{s} = \eta_{s+1}$ and }
$$M\le h_\mu(F, \eta_{s+1}) =  h_\mu(F,   \xi^s\vee \partP\vee\zeta) = h_\mu(F,   \xi^s\vee\zeta)
 \le   \maxent_\mu(F\mid \sigma(\zeta)) $$
 and  the result follows.  If
$\lambda_{s+1} =0$ then  Proposition \ref{prop:key} gives $$3 \epsilon \dim(M)\ge (1-\epsilon) (M - h_\mu(F,   \xi^s\vee\zeta) - 2 \epsilon).$$ As  $\epsilon<\eps$  is arbitrary,  we have  $M\le  h_\mu(F,   \xi^s\vee\zeta)$ and \ref{ineqIIII} follows.  

For \ref{ineqV}, let $\zeta'$ be a sub-invariant measurable partition with $\zeta\prec \zeta'$.  
We may construct $\eta'_i = \eta_i\vee \zeta'$ and apply \cref{prop:key} with the family $\eta'_i$.   Let $\wtd \gamma_i'(\omega,x)$ denote the corresponding transverse dimension.   
Since $\eta_i\prec\eta_i'$,  \cite[Lem.~11.3.2]{MR819557} implies $\wtd \gamma_i'(\omega,x)\le \wtd \gamma_i(\omega,x)$ for almost every $(\omega,x)$.  
\cref{prop:key} combined with \cref{claim:oink} then implies 
\begin{align*}
 (-\lambda_i + 3\epsilon) \big(\underline{\dim}^i( \mu\mid \zeta) - &\underline{\dim}^{i-1}(\mu\mid \zeta)\big)
\ge (1-\epsilon)\left(h_\mu(F,\eta_i') -  h_{\mu}(F,\eta_{i-1}') - 2\epsilon  \right).
 \end{align*}
Combined with  \cref{cor:entropysame} and \ref{ineqI} applied to $h_\mu(F, \xi^i\vee \zeta')$ we then have
 \begin{align*}
 (-\lambda_i + 3\epsilon) \big(\underline{\dim}^i( \mu\mid \zeta) - &\underline{\dim}^{i-1}(\mu\mid \zeta)\big)
\\& \ge (1-\epsilon)\left(h_\mu(F,\eta_i') -  h_{\mu}(F,\eta_{i-1}') - 2\epsilon  \right)\\
 &=(1-\epsilon)\left(h_\mu(F, \xi^i\vee \zeta')- h_\mu(F, \xi^{i-1}\vee \zeta')-2\epsilon \right)\\
& \ge (1-\epsilon) \left(  -\lambda_i  \big( \overline{\dim}^i(   \mu\mid \zeta') - \overline{\dim}^{i-1}(    \mu\mid \zeta')\big)  - 2\epsilon  \right).
 \end{align*}
Arbitrariness of $\epsilon>0$ then implies  \ref{ineqV}.

\subsection{Proof of Proposition \ref{prop:key}}

The proof closely follows   the proof of \cite[Prop.\  5.1]{MR819556} (and \cite[Prop.~11.2]{MR819557}); we include it adapted to  our notation for completeness.

	For  $1\le i\le s+1$ with $\lambda_i\le 0$, define measurable functions  $g, g_\delta, g_*\colon   X \to \R$ as follows:
\begin{align*}
g(\omega, x) & :=   \mu^{\eta_{i-1}}_{(\omega, x)} (F \eta_{i}(\omega, x)) 
	 =   \mu^{\eta_{i-1}}_{(\omega, x)} (F \eta_{i-1}(\omega, x)) \\
 g_\delta(\omega, x) & :=\frac {1}{   \mu^{\eta_{i}}_{(\omega, x)} (\wtd U^{i}((\omega, x),\delta))}    \int _{\wtd U^{i}((\omega, x),\delta)}       \mu^{\eta_{i-1}}_{(\omega' ,z)} ((F \eta_{i})(\omega, x) )\, d   \mu^{\eta_{i}}_{(\omega, x)}(\omega', z) \\
 & = 	\frac {{  \mu^{\eta_{i}}_{(\omega, x)} \left[(\wtd U^{i}((\omega,x),\delta)) \cap (F \eta_i)(\omega,x)\right]}}{ {  \mu^{\eta_{i}}_{(\omega, x)} (\wtd U^{i}((\omega,x),\delta))}}   \\
	 g_*(\omega, x) & := \inf _{\delta>0} g_\delta(\omega,x)
	 \end{align*}
		The equality asserted after the definition of $g(\omega,x)$ follows from Lemma \ref{lemma:44}\ref{PP2}.

For $1\le i\le \min\{s+1,p\}$ and $\mu$-almost every $(\omega,x)$, the metric $\bar d^{ i}$  identifies the quotient $\left(\eta_{i}(\omega,x)/\eta_{i-1}\right)$ with a subset of  $\R^{\sum_{j\ge i} m_j}$. 
This identification gives a measurable map $\eta_i(\omega,x) \to \R^{\sum_{j\ge i} m_j}$ 
that is constant on elements of $\eta_{i-1}$.  
We may apply \cite[Lem.\  4.1.3]{MR819556} to this setting (taking $\alpha$ in the statement of \cite[Lem.\  4.1.3]{MR819556} to be the partition $\restrict{(F \eta_{i})}{\eta_i(\omega,x)}$ of $\eta_i(\omega,x)$) 
and obtain 
	$$\int -\log g_*(\omega', z) \, d   \mu^{\eta_{i}}_{(\omega, x)}(\omega', z) \le H_{  \mu^{\eta_{i}}_{(\omega, x)}} (F \eta_{i}) + \log c + 1$$ where $c$ is a constant depending only on the dimension of $ \R^{\sum_{j\ge i} m_j}$ 
	coming from the Besicovitch Covering Lemma (see \cite[(4.1)]{MR819556}.)  
	It follows that $g_\delta\to g$ $  \mu$-\ae  and that 
\begin{align*}
\int -\log g_*(\omega', z) \, d   \mu (\omega', z) & \le \int H_{  \mu^{\eta_{i}}_{(\omega, x)}} (F \eta_{i})  \  d   \mu(\omega, x)+ \log c + 1
\\& \le H_\mu(F \eta_i\mid \eta_i ) +\log c + 1\\& <\infty.\end{align*}

 Recall our choice of  $E'\subset E$ with $0<  \mu(E')<\frac{\epsilon^2}{8 \log N}$ relative to which our metric in \eqref{eq:metrics} is defined.  Recall that relative to the metric $ \bar d^{i}(\cdot, \cdot)$, all hyperbolicity is seen only on subsequent returns to $E'$.
 
 Consider $(\omega,x)\in X$ whose backwards orbit returns to $E'$ infinitely often and let $r_0\le 0 $ denote the largest $n\le 0$ with $F^{-n}(\omega,x)\in E'$  and let 
$r_0\le 0 < r_1< \dots$ denote the distinct times when $F^{-r_i}(\omega,x)\in E'$.  For $n$ sufficiently large and for $0\le k <n$ with $r_j \le k< r_{j+1}$ write 
	$$a((\omega,x);n, k) := \wtd U^{i}\left (F^{-k}(\omega,x),e^{(\lambda_i -3\epsilon)(n-r_j)} N^{2j} \right)\subset \eta_i(F^{-k}(\omega,x)).$$
We check from definition the following.
	\begin{claim}\label{claim:thisguy8908}
	$a((\omega,x);n, k) \cap( F \eta_i)(F^{-k}(\omega,x))\subset F(a((\omega,x); n , k+1)).$
	\end{claim} 
\begin{proof}  First consider $r_j \le k< r_{j+1}-1$ for some $j\ge 0$.  By definition of the metric $\bar d^{i}(\cdot, \cdot)$,
\begin{align*}F\inv\bigl(a((\omega,x); n,k)\bigr)\cap \eta_i(F^{-(k+1)}(\omega,x)) &= \wtd U^{i}\left (
F^{-(k+1)}(\omega,x),e^{(\lambda_i - 3\epsilon)(n-r_j)} N^{2j} \right)
\\&= a((\omega,x); n,k+1).\end{align*}
If $k=r_{j+1}-1$ for some $j\ge 0$  then 
\begin{align*}
F\inv&\bigl(a((\omega,x); n,k)\bigr)\cap \eta_i(F^{-(k+1)}(\omega,x)) \\
&= F^{-r_{j+1}+ r_j}\left(\wtd U^{i}\left (
F^{-r_j}(\omega,x),e^{(\lambda_i - 3\epsilon)(n-r_j)} N^{2j} \right)\right)\cap \eta_i(F^{-r_{j+1}}(\omega,x)) \\
&\subset \wtd U^{i}\left (
F^{-r_{j+1}}(\omega,x),e^{(\lambda_i - 3\epsilon)(n-r_j)}
e^{-(\lambda_i - 3\epsilon)(r_{j+1}- r_j)} N^{2j}N^2 \right)\\
&=a((\omega,x); n, k+1)
\end{align*}
where the growth by at most $e^{-(\lambda_i - 3\epsilon)(r_{j+1}- r_j)}$ along the backward orbit from $F^{-r_{j}}(\omega,x)$ to $F^{-r_{j+1}}(\omega,x)$ is bounded above by Claim \ref{growthonorbit} and the factor of $N^2$ comes by converting between the metrics $d^{T_i} $ and $d_{F^{-r_{j}}(\omega,x)}^i$ and $d_{F^{-r_{j+1}}(\omega,x)}^i$ using Claim \ref{claim:DontGiveMeNoLipBoy}.
\end{proof}

\begin{proof}[Proof of Proposition \ref{prop:key}]
Consider the quantity  $  \mu_{(\omega, x)} ^{\eta_{i}}(a((\omega,x); n, 0))$.   
Write $p = p(n) := \lfloor (1-\epsilon)n\rfloor$. 
Then
\begin{align*}   \mu_{(\omega, x)} ^{\eta_{i}}&(a((\omega,x); n, 0)) \\&=   \mu_{F^{-p}(\omega,x)}^{\eta_{i}}(a((\omega,x);n,p) )\prod_{k=0} ^{p-1} \dfrac{  \mu_{F^{-k}(\omega,x)} ^{\eta_{i}}(a((\omega,x);n,k))  }{{  \mu_{F^{-(k+1)}(\omega,x)} ^{\eta_{i}}(a((\omega,x);n,k+1))  }} \\ &  \le \prod_{k=0} ^{p-1} \dfrac{  \mu_{F^{-k}(\omega,x)} ^{\eta_{i}}(a((\omega,x);n,k))  }{{  \mu_{F^{-(k+1)}(\omega,x)} ^{\eta_{i}}(a((\omega,x);n,k+1))  }}.\end{align*}
 Renormalizing and applying Claim \ref{claim:thisguy8908},  
\begin{align*}
&\dfrac{  \mu_{F^{-k}(\omega,x)} ^{\eta_{i}}(a((\omega,x);n,k))  }{{  \mu_{F^{-(k+1)}(\omega,x)} ^{\eta_{i}}(a((\omega,x);n,k+1))  }} \\
&\quad\quad   
= 
  \mu_{F^{-k}(\omega,x)} ^{\eta_{i}}(a((\omega,x);n,k))  
\dfrac{  \mu_{F^{-k}(\omega,x)} ^{\eta_{i}}( (F \eta_i)(F^{-k}(\omega,x)))  }{{  \mu_{F^{-k}(\omega,x)} ^{\eta_{i}}(F(a((\omega,x);n,k+1)))  }} \\
&\quad\quad   \le
\dfrac{   \mu_{F^{-k}(\omega,x)} ^{\eta_{i}}(a((\omega,x);n,k))  }{{  \mu_{F^{-k}(\omega,x)} ^{\eta_{i}}\left[(F\eta_i)(F^{-k}(\omega,x))\cap(a((\omega,x);n,k))\right]  }}  \mu_{F^{-k}(\omega,x)} ^{\eta_{i}}( (F \eta_i)(F^{-k}(\omega,x)))\\
&\quad\quad   = [g_{\delta((\omega,x);n,k)}(F^{-k}(\omega,x)){ ] \inv} e^{-I(F^{-k}{ (\omega,x)})}
\end{align*}
where $$\delta((\omega,x);n,k):= e^{(\lambda_i - 3 \epsilon) (n-r_{j_k})}N^{2j_k}$$ for  $r_{j_k}\le k< r_{j_k +1}$ and $$I(\omega',y) := - \log   \mu^{\eta_i}_{(\omega',y)}((F \eta_i){(\omega',y)}).$$  
Then  
\begin{align}
(-&\lambda_i +3 \epsilon) \liminf_{r\to 0}  \dfrac {\log\left(   \mu_{(\omega, x)} ^{\eta_{i}}\wtd U^{i}((\omega,x), r)   \right)}{\log r} \notag \\
&= (-\lambda_i + 3 \epsilon) \liminf_{n\to \infty }  \dfrac {\log\left(   \mu_{(\omega, x)} ^{\eta_{i}}\wtd U^{i}((\omega,x),  e^{(\lambda_i - 3 \epsilon) (n-r_0)} )    \right)}{\log e^{(\lambda_i - 3 \epsilon) (n-r_0)}} \notag \\
&= (-\lambda_i + 3 \epsilon) \liminf_{n\to \infty }  \dfrac {\log\left(   \mu_{(\omega, x)} ^{\eta_{i}} a((\omega,x); n,0)   \right)}{\log e^{(\lambda_i - 3 \epsilon) n} }\notag \\
&= \liminf_{n\to \infty }  \dfrac {1}{n}  -\log\left(   \mu_{(\omega, x)} ^{\eta_{i}} a((\omega,x); n,0)   \right)\notag \\
&\ge \liminf_{n\to \infty} \frac 1 n \sum _{k = 0 } ^{p-1}  \log g_{\delta((\omega,x); n, k)} (F^{-k}(\omega,x)) 
+  \lim_{n\to \infty} \frac 1 n \sum _{k=0} ^{p-1} I(F^{-k}(\omega,x)).  \label{eq:999}
\end{align}
For almost every $(\omega,x)$,  $\lim_{n\to \infty} \frac 1 {{n}} \sum _{k=0} ^{{p-1}} I(F^{-k}(\omega,x))  $ exists and 
\begin{equation}\label{eq:998} \lim_{n\to \infty} \frac 1 n \sum _{k=0} ^{{p-1}} I(F^{-k}(\omega,x)) = (1-\epsilon) H_{  \mu} (F \eta_i \mid \eta_i) = (1-\epsilon)  h_{  \mu} (F, \eta_i).\end{equation}

Now, let $$A_{\delta'}:= \{(\omega,x): -\log g_\delta(\omega,x) \le -\log g(\omega,x) + \epsilon \text{ for all $\delta\le \delta'$}\}.$$
Since $\int -\log g_* <\infty$, taking $\delta'>0$ sufficiently small we can ensure that $$\int _{X\sm A_{\delta'}} -\log g_*  \, d \mu<\epsilon.$$

For a.e.\ $(\omega,x)$, we may  take $n$ sufficiently large so that  $$J_n= \#\{0\le j \le n : F^{-j}(\omega,x)\in E'\}\le 2   n   \mu(E').$$  By the choice of $  \mu(E')$, for large enough $n$  we have $$\delta((\omega,x); n,k):= 
e^{(\lambda_i - 3 \epsilon) (n-r_{j_k})}N^{2j_k} \le 
e^{-\epsilon (n-p(n)) + 2 \log N J_n} \le e^{(-\epsilon^2+ \epsilon^2/2)n} \le \delta'$$ for all $0\le k\le p(n) = \lfloor (1-\epsilon)n\rfloor$.

Then, for such $(\omega,x)$ and $n$, 
\begin{align*}
\sum _{k = 0 } ^ {p-1} & - \log g_{\delta((\omega,x); n, k)} (F^{-k}(\omega,x)) \\
&\le 
\sum _{\{0\le k  \le {p-1} :  F^{-k}(\omega,x)\in A_{\delta'} \}}  - \log g_{\delta((\omega,x); n, k)} (F^{-k}(\omega,x)) \\&\quad \quad \quad+
\sum _{\{0\le k  \le {p-1} :  F^{-k}(\omega,x)\notin A_{\delta'} \}}  - \log g_{\delta((\omega,x); n, k)} (F^{-k}(\omega,x)) \\
&\le \sum _{\{0\le k  \le {p-1} :  F^{-k}(\omega,x)\in A_{\delta'} \}}  \left(- \log g (F^{-k}(\omega,x))+ \epsilon\right) \\&\quad \quad \quad+
\sum _{\{0\le k  \le {p-1} :  F^{-k}(\omega,x)\notin A_{\delta'} \}}  - \log g_{*} (F^{-k}(\omega,x)) 
\end{align*}
whence (from the pointwise ergodic theorem),
\begin{align*}\limsup_{n\to \infty}\frac 1  { n} \sum _{k = 0 } ^ {p-1} &  - \log g_{\delta((\omega,x); n, k)} (F^{-k}(\omega,x)) \\
&\le (1-\epsilon)\left[\int _{A_{\delta'}}\left(- \log g + \epsilon\right) \, d   \mu  + \int _{X\sm A_{\delta'}} - \log g_{*}\, d   \mu \right]\\
&\le (1-\epsilon)\left[ \int -\log g \, d   \mu + 2\epsilon \right]\\
&=(1-\epsilon)\left[  H_{  \mu} (F \eta_{i-1}\mid \eta_{i-1})+ 2\epsilon \right]\\
&= (1-\epsilon)\left[ h_{  \mu} (F, \eta_{i-1})+ 2\epsilon \right].
\end{align*}
Combining the above with \cref{eq:999} and \cref{eq:998}, the proposition follows.  \end{proof}

\end{document}